\documentclass{amsart}
\usepackage{etex}
\usepackage{fixltx2e}

\usepackage[english,french]{babel}

\usepackage[usenames,dvipsnames]{xcolor}
\usepackage[bookmarksopen,bookmarksdepth=2]{hyperref}
\hypersetup{colorlinks=true,citecolor=NavyBlue,linkcolor=BrickRed,urlcolor=Green} 
\usepackage[mathscr]{eucal}
\usepackage{tikz}
\usepackage{tikz-cd}
\usepackage{tikzsymbols}
\usepackage{adjustbox}
\usepackage{enumitem}
\usepackage{dsfont}

\usepackage{microtype}

\renewcommand{\mathbb}{\mathbf}

\usepackage{amssymb,amsmath,amsfonts,amsthm,epsfig,amscd}
\usepackage{stmaryrd}
\usepackage[all,cmtip,poly]{xy}
\usepackage{color}

\newcommand{\hyp}{\operatorname{hyp}}
\newcommand{\Iw}{\operatorname{Iw}}
\newcommand{\gL}{\mathfrak{L}}

\newcommand{\tx}{\tilde{x}}

\newcommand{\red}{\operatorname{red}}

\newcommand{\To}{\longrightarrow}
\newcommand{\isoto}{\stackrel{\sim}{\To}}
\newcommand{\cotimes}{\, \widehat{\otimes}}

\newcommand{\Rep}{\operatorname{Rep}}

\newcommand{\WBT}{W^{\operatorname{BT}}}
\newcommand{\BT}{\operatorname{BT}}
\newcommand{\ocZ}{\overline{\cZ}}

\newcommand{\Kbar}{\bar{K}}

\newcommand{\ghat}{\hat{g}}

\newcommand{\JH}{\operatorname{JH}}
 \newcommand{\sigmabar   }{\overline{\sigma}}   
\def\iso{\buildrel \sim \over \longrightarrow}

\newcommand{\id}{\operatorname{id}}

\newtheorem{thm}[subsubsection]{Theorem}
\newtheorem{lemma}[subsubsection]{Lemma}
\newtheorem{lem}[subsubsection]{Lemma}

\newtheorem{cor}[subsubsection]{Corollary}

\newtheorem{prop}[subsubsection]{Proposition}

\newtheorem{alemma}[subsection]{Lemma}

\newtheorem{adefn}[subsection]{Definition}

\theoremstyle{definition}
\newtheorem{df}[subsubsection]{Definition}
\newtheorem{defn}[subsubsection]{Definition}

\theoremstyle{remark}
\newtheorem{remark}[subsubsection]{Remark}
\newtheorem{rem}[subsubsection]{Remark}

\def\numequation{\addtocounter{subsubsection}{1}\begin{equation}}
\def\nummultline{\addtocounter{subsubsection}{1}\begin{multline}}
\def\anumequation{\addtocounter{subsection}{1}\begin{equation}}
\def\anummultline{\addtocounter{subsection}{1}\begin{multline}}
\renewcommand{\theequation}{\arabic{section}.\arabic{subsection}.\arabic{subsubsection}}

\newif\iffinalrun
\finalruntrue
\iffinalrun
\else
\fi

\iffinalrun
  \newcommand{\need}[1]{}
  \newcommand{\mar}[1]{}
\else
  \newcommand{\need}[1]{{\tiny *** #1}}
  \newcommand{\mar}[1]{\marginpar{\raggedright\tiny fixme #1}}
\fi

\newcommand{\F}{\FF}

\newcommand{\Q}{\QQ}

\newcommand{\Z}{\ZZ}

\newcommand{\FF}{{\mathbb F}}
\newcommand{\GG}{{\mathbb G}}

\newcommand{\QQ}{{\mathbb Q}}

\newcommand{\ZZ}{{\mathbb Z}}

\renewcommand{\bf}{\ensuremath{\mathbf{f}}}

\newcommand{\cC}{{\mathcal C}}
\newcommand{\cD}{{\mathcal D}}
\newcommand{\cE}{{\mathcal E}}
\newcommand{\cF}{{\mathcal F}}
\newcommand{\cG}{{\mathcal G}}

\newcommand{\cL}{{\mathcal L}}

\newcommand{\cM}{{\mathcal M}}

\newcommand{\cO}{{\mathcal O}}
\newcommand{\cP}{{\mathcal P}}

\newcommand{\cR}{{\mathcal R}}

\newcommand{\cW}{{\mathcal W}}
\newcommand{\cX}{{\mathcal X}}
\newcommand{\cY}{{\mathcal Y}}
\newcommand{\cZ}{{\mathcal Z}}

\newcommand{\gM}{\mathfrak{M}}
\newcommand{\gN}{\mathfrak{N}}
\newcommand{\gS}{\mathfrak{S}}

\newcommand{\tB}{\mathrm{B}}

\newcommand{\Fbar}{\overline{\F}}
\newcommand{\Qbar}{\overline{\Q}}
\newcommand{\Zbar}{\overline{\Z}}

\newcommand{\Fp}{\F_p}

\newcommand{\Fpbar}{\Fbar_p}

\newcommand{\Zp}{\Z_p}

\newcommand{\Zpbar}{\Zbar_p}

\newcommand{\Qp}{\Q_p}

\newcommand{\Qpbar}{\Qbar_p}

\DeclareMathOperator{\Aut}{Aut}
\DeclareMathOperator{\coker}{coker}

\DeclareMathOperator{\End}{End}
\DeclareMathOperator{\Ext}{Ext}

\DeclareMathOperator{\Gal}{Gal}
\DeclareMathOperator{\GL}{GL}

\DeclareMathOperator{\Hom}{Hom}

\DeclareMathOperator{\im}{im}

\DeclareMathOperator{\Spec}{Spec}

\DeclareMathOperator{\Spf}{Spf}

\DeclareMathOperator{\Sym}{Sym}

\newcommand{\cris}{\mathrm{cris}}

\newcommand{\Bcris}{\tB_{\cris}}

\newcommand{\Dpcris}{\operatorname{D_{pcris}}}

\newcommand{\rhobar}{\overline{\rho}}

\newcommand*{\invlim}{\varprojlim} 
\newcommand{\into}{\hookrightarrow}
\newcommand{\onto}{\twoheadrightarrow}
\newcommand{\toisom}{\buildrel\sim\over\to}

\newcommand{\longnto}[1]{\buildrel#1\over\longrightarrow}

\newcommand{\Gm}{\GG_m}

\newcommand{\dd}{\mathrm{dd}}
\newcommand{\loc}{\mathrm{loc}}

\newcommand{\varepsilonbar}{\overline{\varepsilon}}
\newcommand{\rbar}{\overline{r}}
\newcommand{\Res}{\operatorname{Res}}

\newcommand{\czero}{c}

\newcommand{\overphi}{/_{\hskip -2pt \phi \hskip 2pt}}

\begin{document}
\selectlanguage{english}
\title[Local geometry of moduli stacks]{Local geometry of moduli stacks of two-dimensional Galois
representations}

\author[A. Caraiani]{Ana Caraiani}\email{a.caraiani@imperial.ac.uk}
\address{Department of Mathematics, Imperial College London, London SW7 2AZ, UK}

\author[M. Emerton]{Matthew Emerton}\email{emerton@math.uchicago.edu}
\address{Department of Mathematics, University of Chicago,
5734 S.\ University Ave., Chicago, IL 60637, USA}

\author[T. Gee]{Toby Gee} \email{toby.gee@imperial.ac.uk} \address{Department of
  Mathematics, Imperial College London,
  London SW7 2AZ, UK}

\author[D. Savitt]{David Savitt} \email{savitt@math.jhu.edu}
\address{Department of Mathematics, Johns Hopkins University}
\thanks{The first author was supported in part by NSF grant DMS-1501064, 
  by a Royal Society University Research Fellowship,
  and by ERC Starting Grant 804176. The second author was supported in part by the
  NSF grants DMS-1303450, DMS-1601871,
  DMS-1902307, DMS-1952705, and DMS-2201242.  The third author was 
  supported in part by a Leverhulme Prize, EPSRC grant EP/L025485/1, Marie Curie Career
  Integration Grant 303605, 
  ERC Starting Grant 306326, a Royal Society Wolfson Research
  Merit Award, and ERC Advanced Grant 884596. The fourth author was supported in part by NSF CAREER
  grant DMS-1564367 and  NSF grant 
  DMS-1702161.}

\begin{abstract}
 We construct moduli stacks of two-dimensional mod~$p$ representations of the
 absolute Galois group of a $p$-adic local field, as well as their
 resolutions by moduli stacks of two-dimensional Breuil--Kisin modules
 with tame descent data. We study  the
 local geometry of these moduli stacks by comparing them with local
 models of Shimura varieties  at hyperspecial and Iwahori level.
\end{abstract}
\maketitle

\selectlanguage{english}

\setcounter{tocdepth}{1}
\tableofcontents

\section{Introduction}\subsection{Moduli of Galois
  representations}Let~$K/\Qp$ be a finite extension,
let~$\overline{K}$ be an algebraic closure of~$K$, and let
$\rbar:\Gal(\overline{K}/K)\to\GL_d(\Fpbar)$ be a continuous
representation. 
The
theory of deformations of $\rbar$ --- that is, liftings of $\rbar$ to
continuous representations $r:\Gal(\overline{K}/K)\to\GL_d(A)$, where $A$ is
a complete local ring with residue field $\Fpbar$ --- is extremely
important in the Langlands program, and in particular is crucial for
proving automorphy lifting theorems via the Taylor--Wiles method.
Proving such theorems often comes down to studying the moduli spaces
of those deformations which satisfy various $p$-adic Hodge-theoretic
conditions. 

From the point of view of algebraic geometry, it seems unnatural to
study only formal deformations of this kind, and Kisin observed about
fifteen years ago that results on the reduction modulo $p$ of
two-dimensional crystalline representations suggested that there
should be moduli spaces of $p$-adic representations in which the
residual representations $\rbar$ should be allowed to vary. In
particular, the special fibres of these moduli spaces 
would be moduli
spaces of mod~$p$ representations of 
$\Gal(\Kbar/K)$. 

In this paper and its sequels~\cite{cegsC,cegsA} we construct such a
space (or rather stack) $\cZ$ of mod $p$ representations in the
case~$d=2$, and describe its geometry. In particular, over the course
of these three papers we show that their irreducible components are
naturally labelled by Serre weights, and that our spaces give a
geometrisation of the weight part of Serre's conjecture.

In the present paper we prove in particular the following
theorem (see Theorem~\ref{thm: existence of picture with descent data and its basic
    properties} and Proposition~\ref{prop: dimensions of the Z stacks},
  as well as Definition~\ref{defn: tres ram} for the notion of a tr\`es
  ramifi\'ee representation). 
\begin{thm}
  \label{thm:intro thm on Z}There is  an algebraic stack~$\cZ$ of
  finite type over~$\Fpbar$, which
  is equidimensional of
  dimension~$[K:\Qp]$, and whose $\Fpbar$-points are in natural
  bijection with the continuous representations
  $\rhobar:G_K\to\GL_2(\Fpbar)$ which are not tr\`es ramifi\'ee. 
\end{thm}

\subsection{The construction}The reason that we restrict to the case
of two-dimensional representations is that in this case one knows that
most mod~$p$ representations are ``tamely potentially finite flat'';
that is, after restriction to a finite tamely ramified extension, they
come from the generic fibres of finite flat group schemes.  Indeed,
the only representations not of this form are the so-called tr\`es
ramifi\'ee representations, which are twists of certain extensions of the
trivial character by the mod~$p$ cyclotomic character, and can be
described explicitly in terms of Kummer theory. (This is a local
Galois-theoretic analogue of the well-known fact that, up to twist,
modular forms of any weight and level $\Gamma_1(N)$, with $N$ prime to
$p$, are congruent modulo~$p$ to modular forms of weight two and level
$\Gamma_1(Np)$;\ the corresponding modular curves acquire semistable
reduction over a tamely ramified extension of~$\Qp$.)

These Galois representations, and the corresponding finite flat group
schemes, can be described in terms of semilinear algebra data. Such
descriptions also exist for more general $p$-adic Hodge theoretic
conditions (such as being crystalline of prescribed Hodge--Tate
weights), although they are more complicated,
and can be used to construct 
analogues
for higher dimensional representations
of the moduli stacks we construct here;
this construction is the subject of 
~\cite{EGmoduli}.

The semilinear algebra data that we use in the present paper are Breuil--Kisin modules
and \'etale $\varphi$-modules. A Breuil--Kisin module is a module with
Frobenius over a power series ring, satisfying a condition on the
cokernel of the Frobenius which depends on a fixed integer, called the
\emph{height} of the Breuil--Kisin module. Inverting the formal variable in the power
series ring gives a functor from the category of Breuil--Kisin modules to the
category of 
\'etale $\varphi$-modules. By Fontaine's
theory~\cite{MR1106901}, these \'etale $\varphi$-modules correspond to
representations of~$\Gal(\overline{K}/K_\infty)$, where~$K_\infty$ is
an infinite non-Galois extension of~$K$ obtained by
extracting $p$-power roots of a uniformiser. By work of Breuil and
Kisin (in particular~\cite{kis04}), for \'etale
$\varphi$-modules that arise from a Breuil--Kisin module of height at most~$1$ 
the corresponding
representations admit a natural extension to~$\Gal(\overline{K}/K)$,
and in this way one obtains precisely the finite flat
representations. 
This is the case that we will consider throughout this paper, extended
slightly to
incorporate descent data from a finite tamely ramified
extension~$K'/K$ and thereby allowing us to study tamely potentially finite
flat representations.

Following Pappas and Rapoport~\cite{MR2562795}, we then consider the moduli stack~$\cC^{\dd}$ of rank two
projective Breuil--Kisin modules with descent data, and the moduli stack~$\cR^{\dd}$ of \'etale
$\varphi$-modules with descent data, together with the natural
map~$\cC^{\dd} \to\cR^{\dd}$. We
deduce from the results of~\cite{MR2562795} that the stack~$\cC^{\dd}$ is algebraic (that is, it
is an Artin stack);\ however~$\cR$ is not algebraic, and indeed is infinite-dimensional. (In fact, we consider
versions of these stacks with $p$-adic coefficients, in which
case~$\cC^{\dd}$ is a $p$-adic formal algebraic stack, but we suppress
this for the purpose of this introduction.) The analogous construction
without tame descent data was considered
in~\cite{EGstacktheoreticimages}, where it was shown that one can
define a notion of the ``scheme-theoretic image'' of the morphism $\cC^{\dd}\to\cR^{\dd}$,
and that the scheme-theoretic image is algebraic. Using similar
arguments, we define our moduli stack~$\cZ^{\dd}$ of two-dimensional Galois
representations to be the scheme-theoretic image, under the 
morphism~$\cC^{\dd}\to\cR^{\dd}$, of a substack of $\cC^{\dd}$
cut out by a condition modeling that of being tamely potentially Barsotti--Tate.

By construction, we know that the finite type points of~$\cZ^{\dd}$ are in
bijection with the (non-tr\`es ramifi\'ee) representations
$\Gal(\overline{K}/K)\to\GL_2(\Fpbar)$, and by using standard results on
the corresponding formal deformation problems, we deduce that~$\cZ^{\dd}$ is
equidimensional of dimension~$[K:\Qp]$. The finite type points of~$\cC^{\dd}$
correspond to potentially finite flat models of these Galois
representations, and we are able to deduce that~$\cC^{\dd}$ is also
equidimensional of dimension~$[K:\Qp]$ (at least morally, this is by
Tate's theorem on the uniqueness of prolongations of $p$-divisible
groups).

An important tool in our proofs is that~$\cC^{\dd}$ has rather mild
singularities, and in particular is Cohen--Macaulay and
reduced. 
We show this by relating the singularities of (certain substacks
of)~$\cC^{\dd}$
to the local
models at Iwahori level of Shimura varieties of $\GL_2$-type;\ such a
relationship was first found in~\cite{kis04} (in the context of formal
deformation spaces, with no descent data) and~\cite{MR2562795} (in the
context of stacks of Breuil--Kisin modules, although again
without descent data) and developed further
by AC and Levin in~\cite{2015arXiv151007503C}.

\subsection{Acknowledgements}We would like to thank Ulrich G\"ortz,
Kalyani Kansal, Wansu Kim and Brandon Levin
for helpful comments, conversations, and correspondence. 

AC and TG would like to thank the organisers (Najmuddin Fakhruddin, Eknath
Ghate, Arvind Nair, C.\ S.\ Rajan, and Sandeep Varma) of the
International Colloquium on Arithmetic Geometry at TIFR in January
2020, for the invitation to speak at the conference, for their
hospitality, and for the invitation to submit this article to the
proceedings. We are grateful to the anonymous referee for their careful reading of the paper and valuable suggestions.
 
This project has received funding from the European Research Council (ERC) under the European Union’s Horizon 2020 research and innovation programme (grant agreements No.\ 804176 and No.\ 884596).

\subsection{Notation and conventions}\label{subsec: notation}
\subsubsection*{Topological groups} If~$M$ is an abelian
topological group with a linear topology, then as
in~\cite[\href{https://stacks.math.columbia.edu/tag/07E7}{Tag
  07E7}]{stacks-project} we say that~$M$ is {\em complete} if the
natural morphism $M\to \varinjlim_i M/U_i$ is an isomorphism,
where~$\{U_i\}_{i \in I}$ is some (equivalently any) fundamental
system of neighbourhoods of~$0$ consisting of subgroups. Note that in
some other references this would be referred to as being~{\em complete
  and separated}. In particular, any $p$-adically complete ring~$A$ is
by definition $p$-adically separated.

\subsubsection*{Galois theory and local class field theory} If $M$ is a field, we let $G_M$ denote its
absolute Galois group.
If~$M$ is a global field and $v$ is a place of $M$, let $M_v$ denote
the completion of $M$ at $v$. If~$M$ is a local field, we write~$I_M$
for the inertia subgroup of~$G_M$. 

 Let $p$ be a prime number. 
 Fix a finite extension $K/\Qp$, with
 ring of integers $\cO_K$ and residue field $k$.  Let $e$ and $f$
 be the  ramification and inertial degrees of $K$, respectively, and
 write $\# k=p^f$ for the cardinality of~$k$.   
Let $K'/K$ be a finite
tamely ramified Galois extension. Let $k'$ be the residue field of $K'$, and let $e',f'$ be the
ramification and inertial degrees of $K'$ respectively.

Our representations of $G_K$ will often have coefficients in $\Qpbar$,
a fixed algebraic closure of $\Qp$ whose residue field we denote by~$\Fpbar$. Let $E$ be a finite
extension of $\Qp$ contained in $\Qpbar$ and containing the image of every
embedding of $K'$ into $\Qpbar$. Let $\cO$ be the ring of integers in
$E$, with uniformiser $\varpi$ and residue field $\F \subset
\Fpbar$.  

Fix an embedding $\sigma_0:k'\into\F$, and recursively define
$\sigma_i:k'\into\F$ for all $i\in\Z$ so that
$\sigma_{i+1}^p=\sigma_i$;\ of course, we have $\sigma_{i+f'}=\sigma_i$
for all~$i$. We let $e_i\in k'\otimes_{\Fp} \F$ denote the idempotent
satisfying $(x\otimes 1)e_i=(1\otimes\sigma_i(x))e_i$ for all $x\in
k'$;\ note that $\varphi(e_i)=e_{i+1}$. We also denote by $e_i$ the
natural lift of $e_i$ to an idempotent in
$W(k')\otimes_{\Zp}\cO$. If $M$ is an
$W(k')\otimes_{\Zp}\cO$-module, then we write $M_i$ for
$e_iM$.

\subsubsection*{$p$-adic Hodge theory} We normalise Hodge--Tate weights so that all Hodge--Tate weights of
the cyclotomic character are equal to $-1$. We say that a potentially
crystalline representation $\rho:G_K\to\GL_2(\Qpbar)$ has \emph{Hodge
  type} $0$, or is \emph{potentially Barsotti--Tate}, if for each
$\varsigma :K\into \Qpbar$, the Hodge--Tate weights of $\rho$ with
respect to $\varsigma$ are $0$ and $1$.  (Note that this is a more
restrictive definition of potentially Barsotti--Tate than is sometimes
used;\ however, we will have no reason to deal with representations
with non-regular Hodge-Tate weights, and so we exclude them from
consideration. Note also that it is more usual in the literature to
say that $\rho$ is potentially Barsotti--Tate if it is potentially
crystalline, and $\rho^\vee$ has Hodge type $0$.) 

We say
that a potentially crystalline representation
$\rho:G_K\to\GL_2(\Qpbar)$ has  \emph{inertial type} $\tau$ if the traces of
elements of $I_K$ acting on~$\tau$ and on
\[\Dpcris(\rho)=\varinjlim_{K'/K}(\Bcris\otimes_{\Qp}V_\rho)^{G_{K'}}\] are
equal (here~$V_\rho$ is the underlying vector space
of~$V_\rho$). 
A representation $\rbar:G_K\to\GL_2(\Fpbar)$ \emph{has a potentially
    Barsotti--Tate lift of
    type~$\tau$} if and
  only if $\rbar$ admits a lift to a representation
  $r:G_K\to\GL_2(\Zpbar)$ of Hodge type~$0$ and inertial type~$\tau$.

  \subsubsection*{Stacks}We will use the terminology and results
  theory of algebraic stacks as developed in~\cite{stacks-project},
  and that of formal algebraic stacks developed
  in~\cite{Emertonformalstacks}. For an overview of this material we
  refer the reader to~\cite[App.\ A]{EGmoduli}. 
For a commutative ring $A$, an \emph{fppf stack over $A$} (or
\emph{fppf} $A$-stack) is a stack fibred in groupoids over the big \emph{fppf}
site of $\Spec A$. 

\subsubsection*{Scheme-theoretic images}We briefly remind the reader of some 
definitions from~\cite[\S3.2]{EGstacktheoreticimages}. 
Let
 $\cX \to \cF$ be a proper morphism of stacks over a locally
Noetherian base-scheme~$S$, 
where $\cX$ is an algebraic stack which is locally of finite presentation over~$S$,
and the diagonal of $\cF$ is representable by algebraic spaces and locally of
finite presentation.

We refer to~\cite[Defn.\ 3.2.8]{EGstacktheoreticimages} for the
definition of the \emph{scheme-theoretic image}~$\cZ$ of the proper morphism $\cX \to
\cF$. By definition, it is a full subcategory in groupoids of~$\cF$, and in fact
by~\cite[Lem.\ 3.2.9]{EGstacktheoreticimages} it is a Zariski substack
of~$\cF$. By~\cite[Lem.\ 3.2.14]{EGstacktheoreticimages}, the finite type points
of~$\cZ$ are precisely the finite type points of~$\cF$ for which the
corresponding fibre of~$\cX$ is nonzero. 

 The results of~\cite[\S3.2]{EGstacktheoreticimages} give criteria
for~$\cZ$ to be an algebraic stack, and prove a number of associated results
(such as universal properties of the morphism $\cZ\to\cF$, and a description of
versal deformation rings for~$\cZ$);\ rather than recalling these results in detail
here, we will refer to them as needed in the body of the paper.

\section{Integral $p$-adic Hodge theory with tame descent data}\label{sec:kisin modules with descent
  data} In this section we introduce various objects in semilinear algebra which
arise in the study of potentially Barsotti--Tate Galois representations with tame
descent data. Much of this material is standard, and none of it will surprise an
expert, but we do not know of a treatment in the literature in the level of
generality that we require;\ in particular, we are not aware of a
treatment of the theory of tame descent data for Breuil--Kisin modules
with coefficients. However, the arguments are almost identical
to those for strongly divisible
modules and Breuil modules, so we will be brief.

The various equivalences of categories between the objects we consider and
finite flat group schemes or $p$-divisible groups will not be
relevant to our main arguments, except at a motivational level, so we
largely ignore them.

\subsection{Breuil--Kisin modules and \texorpdfstring{$\varphi$}{phi}-modules
  with descent data}\label{subsec: kisin modules with dd} 
Recall that we have a finite
tamely ramified Galois extension $K'/K$. Suppose further that there exists a uniformiser $\pi'$ of
$\cO_{K'}$ such that $\pi:=(\pi')^{e(K'/K)}$ is an element of~$K$, 
where $e(K'/K)$ is the ramification index of
$K'/K$. 
Recall that $k'$ is the residue field of $K'$, while $e',f'$ are the
ramification and inertial degrees of $K'$ respectively.
 Let $E(u)$ be the minimal polynomial of $\pi'$ over $W(k')[1/p]$. 

Let $\varphi$ denote the arithmetic Frobenius automorphism of $k'$, which lifts uniquely
to an automorphism of $W(k')$ that we also denote by $\varphi$. Define
$\gS:=W(k')[[u]]$, and extend $\varphi$ to $\gS$ by \[\varphi\left(\sum a_iu^i\right)=\sum
\varphi(a_i)u^{pi}.\] By our assumptions that $(\pi')^{e(K'/K)} \in K$ 
and that $K'/K$ is Galois, for each
$g\in\Gal(K'/K)$ we can write $g(\pi')/\pi'=h(g)$ with $h(g)\in
\mu_{e(K'/K)}(K') \subset W(k')$,
and we
let $\Gal(K'/K)$ act on $\gS$ via \[g\left(\sum a_iu^i\right)=\sum g(a_i)h(g)^iu^i.\]

Let $A$ be a $p$-adically complete $\Zp$-algebra, set $\gS_A:=(W(k')\otimes_{\Zp} A)[[u]]$, and extend
the actions of $\varphi$ and $\Gal(K'/K)$ on $\gS$ to actions on $\gS_A$ in the
obvious ($A$-linear) fashion.

The actions of $\varphi$ and $\Gal(K'/K)$ on $\gS_A$
extend to actions on $\gS_A[1/u]=(W(k')\otimes_{\Zp} A)((u))$ in the obvious
way.  It will sometimes be necessary to consider the subring $\gS_A^0
:=(W(k)\otimes_{\Zp} A)[[v]]$ of $\gS_A$
  consisting of power series in 
$v:=u^{e(K'/K)}$, on which  $\Gal(K'/K)$ acts
  trivially. 

\begin{defn}\label{defn: Kisin module with descent data}
Fix a $p$-adically complete $\Zp$-algebra~$A$. A \emph{weak Breuil--Kisin module with
  $A$-coefficients and descent data from $K'$ to $K$}
  is a triple $(\gM,\varphi_{\gM},\{\hat{g}\}_{g\in\Gal(K'/K)})$ consisting of
  a 
  $\gS_A$-module~$\gM$ and a $\varphi$-semilinear map
  $\varphi_{\gM}:\gM\to\gM$ 
   such that:
  \begin{itemize}
  \item the $\gS_A$-module $\gM$ is finitely generated and $u$-torsion
    free, and 
  \item the  induced
  map $\Phi_{\gM} = 1 \otimes \varphi_{\gM} :\varphi^*\gM\to\gM$ is an isomorphism after
  inverting $E(u)$  (here as usual we write $\varphi^*\gM:=\gS_A \otimes_{\varphi,\gS_A}\gM$), 
  \end{itemize}
together with
  additive bijections  $\hat{g}:\gM\to\gM$, satisfying the further
  properties that 
    the maps $\hat{g}$ commute with $\varphi_\gM$, satisfy 
    $\hat{g_1}\circ\hat{g_2}=\widehat{g_1\circ g_2}$, and have
    $\hat{g}(sm)=g(s)\hat{g}(m)$ for all $s\in\gS_A$, $m\in\gM$. 
   We say that $\gM$ is has \emph{height at most $h$} if the cokernel of 
   $\Phi_{\gM}$ is killed by $E(u)^h$. 

 If $\gM$ as above is projective as an $\gS_A$-module (equivalently, if
   the condition that  the $\gM$ is $u$-torsion free is replaced with the condition
   that $\gM$ is projective) then we say 
 that $\gM$ is a \emph{Breuil--Kisin
   module  with $A$-coefficients and descent data from $K'$ to $K$},
 or even simply that $\gM$ is a \emph{Breuil--Kisin module}. 

The Breuil--Kisin module $\gM$ is said to be  of rank~$d$ if the underlying
finitely generated projective $\gS_A$-module has constant rank~$d$. It
is said to be free if the underlying $\gS_A$-module is free.
\end{defn}

A morphism of (weak) Breuil--Kisin modules with descent data is
a  morphism 
of~$\gS_A$-modules 
that commutes with $\varphi$ and with the~$\hat{g}$. 
In the case that $K'=K$ the data of the $\hat{g}$ is trivial, so
it can be forgotten, giving the category of (\emph{weak}) \emph{Breuil--Kisin modules with
  $A$-coefficients.} In this case it will sometimes be convenient to elide the difference between a
Breuil--Kisin module with trivial descent data, and a Breuil--Kisin module without
descent data, in order to avoid making separate definitions in the
case of Breuil--Kisin modules without descent data;\ the same convention will
apply to the \'etale $\varphi$-modules considered below.

\begin{lem}\label{lem:kisin injective} Suppose either that $A$ is a $\Z/p^a\Z$-algebra for some
  $a \ge 1$, or that $A$ is $p$-adically complete and $\gM$ is projective. Then in
  Definition~\ref{defn: Kisin module with descent data} the condition
  that $\Phi_{\gM}$ is an isomorphism after inverting $E(u)$ may
  equivalently be replaced with the condition that $\Phi_\gM$ is
  injective and its cokernel is killed by a power
  of $E(u)$.
\end{lem}

  \label{rem: Kisin modules are etale} 

  \begin{proof}
   
  If $A$ is a $\Z/p^a\Z$-algebra for some $a\ge 0$, then $E(u)^h$
  divides $u^{e(a+h-1)}$ in $\gS_A$ (see ~\cite[Lem.\
  5.2.6]{EGstacktheoreticimages} and its proof), so that $\gM[1/u]$ is
  \'etale in the sense that the induced
  map \[\Phi_{\gM}[1/u]:\varphi^*\gM[1/u]\to\gM[1/u]\]
  is an isomorphism. The injectivity of $\Phi_{\gM}$ now follows
  because $\gM$, and therefore $\varphi^*\gM$, is $u$-torsion free.

If instead  $A$ is $p$-adically complete, then no Eisenstein
polynomial over $W(k')$ is a zero divisor in $\gS_A$:\ this is plainly
true if $p$ is nilpotent in $A$, from which one deduces the same for
$p$-adically complete $A$. Assuming that $\gM$ is projective, it
follows that the maps $\gM \to \gM[1/E(u)]$ and $\varphi^*\gM \to
(\varphi^*\gM)[1/E(u)]$ are injective, and we are done.
  \end{proof}

\begin{df}
\label{def:completed tensor}
If $Q$ is any (not necessarily finitely generated) $A$-module,
and $\gM$ is an $A[[u]]$-module,
then  we let $\gM\cotimes_A Q$ denote the $u$-adic completion
of $\gM\otimes_A Q$.
\end{df} 

\begin{lem}
  \label{rem: base change of locally free Kisin module is a
    locally free Kisin module}
If $\gM$ is a Breuil--Kisin module and $B$ is an $A$-algebra, then the base
change $\gM \cotimes_A B$ is a Breuil--Kisin module.
\end{lem}

\begin{proof}

We claim that $\gM \cotimes_A B \cong \gM \otimes_{A[[u]]} B[[u]]$ for
any finitely generated projective $A[[u]]$-module;\ the lemma then follows
immediately from Definition~\ref{defn: Kisin module with descent
  data}.

To check the claim,
we must see that
the finitely generated $B[[u]]$-module
$\gM\otimes_{A[[u]]} B[[u]]$ is $u$-adically complete.
But $\gM$ is a direct summand of a free $A[[u]]$-module of finite rank,
in which case
$\gM\otimes_{A[[u]]} B[[u]]$ is a direct summand of a free $B[[u]]$-module
of finite rank and hence is $u$-adically complete.
\end{proof}

\begin{remark}\label{rem:Kisin-mod-ideal}
If $I \subset A$ is a finitely generated ideal then  $A[[u]] \otimes_A
A/I \cong (A/I)[[u]]$, and $\gM \otimes_A A/I
\cong \gM \otimes_{A[[u]]} (A/I)[[u]] \cong \gM \cotimes_A A/I$;
so in this case $\gM \otimes_A A/I$ itself is a Breuil--Kisin module.
\end{remark}

Note that the base change (in the sense of Definition~\ref{def:completed
  tensor})  of a weak Breuil--Kisin
  module may not be a weak
  Breuil--Kisin module, because  the property of being $u$-torsion free is not
  always preserved by base change.

\begin{defn}\label{defn: etale phi module} 
Let $A$ be a $\Z/p^a\Z$-algebra for some $a\ge 1$.  A \emph{weak \'etale
  $\varphi$-module} with $A$-coefficients
  and descent data from $K'$ to $K$ is a triple
  $(M,\varphi_M,\{\hat{g}\})$ consisting of: 
  \begin{itemize}
\item 
a finitely generated
  $\gS_A[1/u]$-module $M$;\ 
\item a  $\varphi$-semilinear map $\varphi_M:M\to M$ with the
  property that the induced
  map   \[\Phi_M = 1 \otimes \varphi_M:\varphi^*M:=\gS_A[1/u]\otimes_{\varphi,\gS_A[1/u]}M\to M\]is an
  isomorphism,
  \end{itemize}
together with   additive bijections  $\hat{g}:M\to M$ for $g\in\Gal(K'/K)$, satisfying the further
  properties that the maps $\hat{g}$ commute with $\varphi_M$, satisfy
    $\hat{g_1}\circ\hat{g_2}=\widehat{g_1\circ g_2}$, and have
    $\hat{g}(sm)=g(s)\hat{g}(m)$ for all $s\in\gS_A[1/u]$, $m\in M$.

If $M$ as above is projective as an $\gS_A[1/u]$-module then we say
simply that $M$ is an \'etale $\varphi$-module. The \'etale $\varphi$-module $M$ is said to be  of rank~$d$ if the underlying
finitely generated projective $\gS_A[1/u]$-module has constant rank~$d$.

\end{defn}

\begin{rem} 
  \label{rem: completed version if $p$ not nilpotent}We could also
  consider \'etale $\varphi$-modules for general $p$-adically complete $\Zp$-algebras~$A$, but
  we would need to replace $\gS_A[1/u]$ by its $p$-adic completion. As
  we will not need to consider these modules in this paper, we do not
  do so here, but we refer the interested reader to~\cite{EGmoduli}. 
\end{rem}
A morphism
of weak \'etale
$\varphi$-modules with $A$-coefficients and descent data from $K'$ to
$K$ 
is a morphism of~$\gS_A[1/u]$-modules 
that commutes with $\varphi$ and with the
$\hat{g}$. Again, in the case $K'=K$ the descent data is trivial, and we
obtain the  usual category of \'etale $\varphi$-modules with
$A$-coefficients. 

Note that
if $A$ is a $\Z/p^a\Z$-algebra, and $\gM$ is a  Breuil--Kisin module (resp.,
weak Breuil--Kisin module)  with descent data, then $\gM[1/u]$ naturally has the
structure of an \'etale $\varphi$-module (resp., weak \'etale $\varphi$-module) with descent data.

Suppose that $A$ is an $\cO$-algebra (where $\cO$ is as in
Section~\ref{subsec: notation}). In making calculations, it is often
convenient to use the idempotents~$e_i$ (again as in
Section~\ref{subsec: notation}). In particular if $\gM$ is a Breuil--Kisin
module, then writing as usual
$\gM_i:=e_i\gM$, we write $\Phi_{\gM,i}:\varphi^*(\gM_{i-1})\to\gM_{i}$ for
the morphism induced by~$\Phi_{\gM}$.  Similarly if $M$ is an
\'etale $\varphi$-module  then we write
$M_i:=e_iM$, and we write $\Phi_{M,i}:\varphi^*(M_{i-1}) \to M_{i}$ for
the morphism induced by~$\Phi_{M}$.

\subsection{Dieudonn\'e modules}\label{subsec:Dieudonne modules} 
Let $A$ be a $\Zp$-algebra. We define a \emph{Dieudonn\'e module of rank $d$ with $A$-coefficients and
  descent data from $K'$ to $K$} to be a finitely generated projective
$W(k')\otimes_{\Zp}A$-module $D$ of constant rank 
$d$ on $\Spec A$, together with:

\begin{itemize}
\item $A$-linear endomorphisms $F,V$ satisfying $FV = VF = p$ such that $F$ is $\varphi$-semilinear and $V$ is
  $\varphi^{-1}$-semilinear for the action of $W(k')$, and
\item a $W(k')\otimes_{\Zp}A$-semilinear action of $\Gal(K'/K)$
which commutes with $F$ and $V$. 
\end{itemize}

\begin{defn}\label{def: Dieudonne module formulas}
If $\gM$ is a Breuil--Kisin module of height at most~$1$ and rank~$d$ with descent data, 
then there is a
corresponding Dieudonn\'e module $D=D(\gM)$ of rank~$d$ defined as follows. We set
$D:=\gM/u\gM$
with the induced action of $\Gal(K'/K)$, and $F$ given by the induced
action of $\varphi$. 
The endomorphism $V$ is determined as follows.  Write $E(0) = \czero p$,
so that we have $p \equiv \czero^{-1}E(u) \pmod{u}$. The
condition that the cokernel of $\varphi^*\gM\to\gM$ is killed by $E(u)$
allows us to factor the multiplication-by-$E(u)$ map on $\gM$ uniquely
as $\mathfrak{V} \circ \varphi$,  and $V$ is
defined to be
$\czero^{-1} \mathfrak{V}$ modulo~$u$. 
\end{defn}

\subsection{Galois representations}\label{subsec: etale phi modules
  and Galois representations}
The theory of fields of norms~\cite{MR526137} 
was used in~\cite{MR1106901} to relate \'etale $\varphi$-modules with descent data to representations
of a certain absolute Galois group;\ not the group $G_K$,
but rather the group $G_{K_{\infty}}$,
where $K_{\infty}$ is a certain infinite extension of~$K$ (whose
definition is recalled below).  Breuil--Kisin modules of height $h \leq 1$ are closely related to finite
flat group schemes (defined over $\cO_{K'}$, but with descent data to $K$
on their generic fibre). Passage from a Breuil--Kisin module to its associated
\'etale $\varphi$-module can morally be interpreted as the
passage from a finite flat group scheme (with descent data) to
its corresponding Galois representation (restricted to $G_{K_{\infty}}$).
Since the generic fibre of a finite flat group scheme over $\cO_{K'}$,
when equipped with descent data to $K$,
in fact gives rise to a representation of $G_K$,
in the case $h = 1$ we may relate Breuil--Kisin modules with descent data
(or, more precisely, their associated \'etale $\varphi$-modules),
not only to representations of $G_{K_{\infty}}$, but to representations
of~$G_K$.

 In 
this subsection, we recall
some results coming from this connection,
and draw some conclusions for Galois deformation rings.

\subsubsection{From \'etale $\varphi$-modules to $G_{K_{\infty}}$-representations}
We begin by recalling from \cite{kis04} some constructions arising in $p$-adic Hodge theory
and the theory of fields of norms, which go back to~\cite{MR1106901}. 
Following Fontaine,
we write $R:=\varprojlim_{x\mapsto
  x^p}\cO_{\Kbar}/p$. 
Fix a compatible system $(\! \sqrt[p^n]{\pi}\,)_{n\ge 0}$ of
$p^n$th roots of $\pi$ in $\Kbar$ (compatible in the obvious sense that 
$\bigl(\! \sqrt[p^{n+1}]{\pi}\,\bigr)^p = \sqrt[p^n]{\pi}\,$),
and let
$K_{\infty}:=\cup_{n}K(\sqrt[p^n]{\pi})$, and
also $K'_\infty:=\cup_{n}K'(\sqrt[p^n]{\pi})$. Since $(e(K'/K),p)=1$, the compatible system
$(\! \sqrt[p^n]{\pi}\,)_{n\ge 0}$ determines a unique compatible system $(\!
\sqrt[p^n]{\pi'}\,)_{n\ge 0}$ of $p^n$th roots of~$\pi'$ such that $(\!
\sqrt[p^n]{\pi'}\,)^{e(K'/K)} =\sqrt[p^n]{\pi}$.
Write
$\underline{\pi}'=(\sqrt[p^n]{\pi'})_{n\ge 0}\in R$, and $[\underline{\pi}']\in
W(R)$ for its image under the natural multiplicative map $R \to W(R)$. We have a Frobenius-equivariant
inclusion $\gS\into W(R)$ by sending $u\mapsto[\underline{\pi}']$.  We can naturally identify
$\Gal(K'_\infty/K_\infty)$ with $\Gal(K'/K)$, and doing this we see that the
action of $g\in G_{K_\infty}$ on $u$ is via $g(u)=h(g)u$.

We let $\cO_{\cE}$ denote the $p$-adic completion of $\gS[1/u]$, and let $\cE$ be the
field of fractions of~$\cO_\cE$. The inclusion $\gS\into W(R)$ extends to an
inclusion $\cE\into W(\operatorname{Frac}(R))[1/p]$. Let $\cE^{\text{nr}}$ be
the maximal unramified extension of $\cE$ in $ W(\operatorname{Frac}(R))[1/p]$,
and let $\cO_{\cE^{\text{nr}}}\subset W(\operatorname{Frac}(R))$ denote
its ring of
integers. Let $\cO_{\widehat{\cE^{\text{nr}}}}$ be the $p$-adic completion of
$\cO_{\cE^{\text{nr}}}$. Note that $\cO_{\widehat{\cE^{\text{nr}}}}$ is stable
under the action of $G_{K_\infty}$. 

 \begin{defn} 
Suppose that $A$ is a $\Z/p^a\Z$-algebra for some $a \ge 1$.  If $M$
is a weak \'etale $\varphi$-module with $A$-coefficients and descent data, set
  $T_A(M):=\left(\cO_{\widehat{\cE^{\text{nr}}}}\otimes_{\gS[1/u]}M
  \right)^{\varphi=1}$, an $A$-module with a
  $G_{K_\infty}$-action (via the diagonal action on
  $\cO_{\widehat{\cE^{\text{nr}}}}$ and $M$, the latter given by
  the~$\hat{g}$). If $\gM$ is a weak Breuil--Kisin module with
  $A$-coefficients and descent data,
  set 
  $T_A(\gM):=T_A(\gM[1/u])$. 
\end{defn}

\begin{lem}
  \label{lem: Galois rep is a functor if A is actually finite local} Suppose
  that $A$ is a local $\Zp$-algebra and that $|A|<\infty$. Then~$T_A$ induces an 
  equivalence of categories from the category of weak \'etale
  $\varphi$-modules with $A$-coefficients and descent data to the category of continuous
  representations of $G_{K_\infty}$ on finite $A$-modules. If $A\to A'$
  is finite, then there is a natural isomorphism $T_A(M)\otimes_A A'\iso
  T_{A'}(M\otimes_A A')$. A weak \'etale $\varphi$-module with
  $A$-coefficients and descent
  data~$M$ is free of rank~$d$ if and only if $T_A(M)$ is a free
  $A$-module of rank~$d$. 
\end{lem}
\begin{proof}
  This is due to Fontaine~\cite{MR1106901}, and can be proved in exactly the same way as~\cite[Lem.\ 1.2.7]{kis04}.
\end{proof}

We will frequently simply write $T$ for $T_A$. Note that if we let
$M'$ be the \'etale $\varphi$-module obtained from $M$ by forgetting the
descent data, then by definition we have
$T(M')=T(M)|_{G_{K'_\infty}}$.

\subsubsection{Relationships between $G_K$-representations
and $G_{K_{\infty}}$-representations}\label{subsubsec: deformation
rings and Kisin modules}
We will later need to study deformation rings for representations
of~$G_K$ in terms of the deformation rings for the restrictions of
these representations to~$G_{K_\infty}$. Note that the representations
of~$G_{K_\infty}$ coming from Breuil--Kisin modules of height at most~$1$
admit canonical extensions to~$G_K$ by~\cite[Prop.\ 1.1.13]{kis04}.

Let $\rbar:G_K\to\GL_2(\F)$
be a continuous representation, let $R_{\rbar}$ denote the universal
framed deformation $\cO$-algebra for~$\rbar$, and let
$R_{\rbar}^{[0,1]}$ 
be the quotient with the property that if~$A$ is
an Artinian local~$\cO$-algebra with residue field~$\F$, 
then
a local $\cO$-morphism $R_{\rbar}\to A$ factors
through $R_{\rbar}^{[0,1]}$ if and only if the corresponding
$G_K$-module (ignoring the $A$-action) admits a $G_K$-equivariant surjection from a potentially crystalline
$\cO$-representation  all of whose Hodge--Tate weights are equal to $0$ or
$1$, and whose restriction to~$G_{K'}$ is crystalline. (The existence of this quotient follows as in~\cite[\S2.1]{MR2782840}.)

Let $R_{\rbar|_{G_{K_\infty}}}$ be the universal framed deformation
$\cO$-algebra for~$\rbar|_{G_{K_\infty}}$, and let
$R_{\rbar|_{G_{K_\infty}}}^{\le 1}$ denote the quotient with the property that if~$A$ is
an Artinian local~$\cO$-algebra with residue field~$\F$, then a morphism $R_{\rbar|_{G_{K_\infty}}}\to A$ factors
through $R_{\rbar|_{G_{K_\infty}}}^{\le 1}$ if and only if the corresponding
$G_{K_\infty}$-module is isomorphic to  $T(\gM)$ for some weak Breuil--Kisin
module~$\gM$ of height at most one with $A$-coefficients and descent
data from~$K'$ to~$K$. 
(The existence of this quotient
follows exactly as for~\cite[Thm.~1.3]{MR2782840}.)

\begin{prop}
  \label{prop: restriction gives an equivalence of height at most 1
    deformation}The natural map induced by restriction from $G_K$ to
  $G_{K_\infty}$ induces an isomorphism $\Spec R_{\rbar}^{[0,1]}\to
  \Spec R_{\rbar|_{G_{K_\infty}}}^{\le 1}$. 
\end{prop}
\begin{proof}
  This can be proved in exactly the same way as~\cite[Cor.\
  2.2.1]{MR2782840} (which is the case that $E=\Qp$ and $K'=K$).
\end{proof}

\section{Moduli stacks of Breuil--Kisin modules with descent data}\label{sec:moduli stacks}

In this section we introduce certain moduli stacks of Breuil--Kisin modules and \'etale
$\varphi$-modules with tame
descent data, following~\cite{MR2562795,EGstacktheoreticimages} (which
consider the case without descent data).

The definitions and basic properties of our moduli stacks of
Breuil--Kisin modules with tame descent
data, and of \'etale $\varphi$-modules  with tame descent data,  are contained in
Section~\ref{subsec:moduli}.  In Section~\ref{subsec:tame inertial types,
  general height} we explain how to decompose the Breuil--Kisin
moduli stacks 
 as a disjoint union of substacks according to the \emph{type} of the
 descent data, making use of results in Section~\ref{subsec: tame
   groups} on the stacks classifying representations of tame groups.

\subsection{Moduli stacks of Breuil--Kisin modules and 
  \texorpdfstring{$\varphi$}{phi}-modules with descent data}
\label{subsec:moduli}
We begin by defining the 
moduli stacks of Breuil--Kisin modules, with and without descent data. 

\begin{defn}
  \label{defn: C^dd,a }For each integer $a\ge 1$, we let $\cC_{d,h,K'}^{\dd,a}$ be
  the {\em fppf} stack over~$\cO/\varpi^a$ which associates to any $\cO/\varpi^a$-algebra $A$
the groupoid $\cC_{d,h,K'}^{\dd,a}(A)$ 
of rank~$d$ Breuil--Kisin modules of height at most~$h$ with $A$-coefficients and descent data from
$K'$ to $K$. 

By~\cite[\href{http://stacks.math.columbia.edu/tag/04WV}{Tag 04WV}]{stacks-project}, 
we may also regard each of the stacks $\cC_{d,h,K'}^{\dd,a}$ as an {\em fppf}
stack over $\cO$,
and we then write $\cC_{d,h,K'}^{\dd}:=\varinjlim_{a}\cC_{d,h,K'}^{\dd,a}$;\ this
is again an {\em fppf} stack over $\cO$.

We will frequently omit any (or all) of the subscripts $d,h,K'$ from this notation
when doing so will not cause confusion. We write $\cC^a$ for $\cC^{\dd,a}$ and $\cC$
for $\cC^{\dd}$ for the stacks with descent data for the trivial extension $K'/K'$.
\end{defn} 

The natural morphism $\cC^{\dd} \to \Spec \cO$ factors through $\Spf \cO$,
and by construction, there is an isomorphism $\cC^{\dd,a} \iso \cC^{\dd}\times_{\Spf \cO}
\Spec \cO/\varpi^a,$ for each $a \geq 1$;\ in particular, each of the morphisms
$\cC^{\dd,a} \to \cC^{\dd,a+1}$
is a thickening
(in the sense that its pullback under any test morphism $\Spec A \to \cC^{\dd,a+1}$
becomes a thickening of schemes, as defined
in~\cite[\href{http://stacks.math.columbia.edu/tag/04EX}{Tag
  04EX}]{stacks-project}\footnote{Note that for morphisms of algebraic
stacks --- and we will see below that $\cC^{\dd,a}$ and $\cC^{\dd,a+1}$ {\em are}
algebraic stacks --- this notion of thickening coincides with the notion defined
in~\cite[\href{http://stacks.math.columbia.edu/tag/0BPP}{Tag
  0BPP}]{stacks-project},
by~\cite[\href{http://stacks.math.columbia.edu/tag/0CJ7}{Tag
  0CJ7}]{stacks-project}.}).
  In Corollary~\ref{cor: C^dd,a is an algebraic stack} below we show that for each integer $a\ge 1$,
  $\cC^{\dd,a}$ is in fact an algebraic stack of finite type over
  $\Spec\cO/\varpi^a$. The stack $\cC^{\dd}$ is then
  {\em a priori} an Ind-algebraic
 stack, endowed with a morphism to $\Spf \cO$ which is representable
by algebraic stacks;\ we show that it 
is in fact a 
  \emph{$\varpi$-adic formal algebraic stack}, in the sense of the following definition whose basic theory has been 
 developed by Emerton~\cite{Emertonformalstacks}. 

\begin{defn}
  \label{df:formal stack}
An {\em fppf} stack in groupoids $\cX$ over a scheme $S$ is called
a {\em formal algebraic stack} if
there is a morphism $U \to \cX$,
  whose domain $U$ is a formal algebraic space over $S$
  (in the sense 
of~\cite[\href{http://stacks.math.columbia.edu/tag/0AIL}{Tag 0AIL}]{stacks-project}),
and which is representable by algebraic spaces,
smooth, and surjective.

  A formal algebraic stack $\cX$ over $\Spec \cO$
  is called \emph{$\varpi$-adic} if the canonical map
  $\cX \to \Spec \cO$ factors 
  through $\Spf \cO$,
  and if the induced map $\cX \to \Spf \cO$
  is algebraic, i.e.\ representable by algebraic
  stacks (in the sense 
of~\cite[\href{http://stacks.math.columbia.edu/tag/06CF}{Tag 06CF}]{stacks-project} and \cite[Def.~3.1]{Emertonformalstacks}).
\end{defn}

Our approach will be to deduce the statements in the case with descent 
data from the corresponding statements in the case with no descent
data, which follow from the methods of Pappas and Rapoport~\cite{MR2562795}.  
More precisely, in that reference it is proved that each
$\cC^a$ is an algebraic stack over~$\cO/\varpi^a$~\cite[Thm.~0.1 (i)]{MR2562795},
and thus that $\cC := \varinjlim_a \cC^a$ is a $\varpi$-adic
Ind-algebraic stack (in the sense that it is an Ind-algebraic stack
with a morphism to~$\Spf\cO$ that is representable by algebraic
stacks). 
(In \cite{MR2562795} the stack
	$\cC$ is described as being a $p$-adic formal algebraic stack. 
	However, in that reference, this
	term is used synonomously with our notion of a $p$-adic 
	Ind-algebraic stack;\ the question of the existence of a smooth cover
	of $\cC$ by a $p$-adic formal algebraic space is not discussed.
	As we will see, though, the existence of such a cover is easily
	deduced from the set-up of \cite{MR2562795}.)

We thank Brandon Levin for pointing out the following
result to us.  The proof is essentially (but somewhat 
implicitly) contained in the proof of 
\cite[Thm.~3.5]{2015arXiv151007503C},
but we take the opportunity to make it explicit. Note that it could
also be directly deduced from the results of~\cite{MR2562795}
using~\cite[Prop.\ A.10]{EGmoduli}, but the proof that we give has the
advantage of giving an explicit cover by a formal algebraic space. 

\begin{prop}
  \label{prop: C without dd is a $p$-adic formal  algebraic
    stack}For any choice of $d,h$, $\cC$ is a 
$\varpi$-adic formal algebraic stack of finite type over~$\Spf\cO$ with affine diagonal.
\end{prop}
\begin{proof}
  We begin by recalling some results from~\cite[\S 3.b]{MR2562795} (which
  is where the
  proof that each~$\cC^a$ is an algebraic stack of finite type
  over~$\cO/\varpi^a$ is given). If~$A$ is an $\cO/\varpi^a$-algebra for some $a\ge 1$,
  then we set\[L^+G(A):=GL_d(\gS_A),\] \[LG^{h,K'}(A):=\{X\in
    M_d(\gS_A) \mid  X^{-1}\in E(u)^{-h}M_d(\gS_A)\},\]and let $g\in L^+G(A)$ act
  on the right on $LG^{h,K'}(A)$ by $\varphi$-conjugation as $g^{-1}\cdot
  X\cdot \varphi(g)$. Then we may write 
  \[\cC=[LG^{h,K'} \overphi L^+G].\]

  For each $n\ge 1$ we have the principal congruence subgroup $U_n$ of
  $L^+G$ given by $U_n(A)=I+u^n\cdot M_d(\gS_A)$. As in~\cite[\S
  3.b.2]{MR2562795}, for any integer $n(a)>eah/(p-1)$ we have a
  natural identification
  \numequation\label{eqn:U_n phi conjugacy}
[LG^{h,K'}\overphi U_{n(a)}]_{\cO/\varpi^a}\cong [LG^{h,K'}/
    U_{n(a)}]_{\cO/\varpi^a}\end{equation}where
  the $U_{n(a)}$-action on the right hand side is by left
  translation by the inverse;\ 
  moreover this quotient stack is represented by a finite type
  scheme $(X^{h,K'}_{n(a)})_{\cO/\varpi^a}$, and we find that
  \[\cC^a\cong
    [(X^{h,K'}_{n(a)})_{\cO/\varpi^a}\overphi(\cG_{{n(a)}})_{\cO/\varpi^a}],\]
  where $(\cG_{{n(a)}})_{\cO/\varpi^a}=(L^+G/U_{n(a)})_{\cO/\varpi^a}$ is a
  smooth finite type group scheme over $\cO/\varpi^a$.

Now define $Y_a := [(X_{n(a)}^{h,K'})_{\cO/\varpi^a}\overphi (U_{n(1)})_{\cO/\varpi^a}].$
If $a \geq b,$ then there is a natural isomorphism
$(Y_a)_{\cO/\varpi^b} \cong Y_b.$ 
Thus we may form the $\varpi$-adic Ind-algebraic stack
$Y := \varinjlim_a Y_a.$
Since $Y_1 := (X_{n(1)}^{h,K'})_{\F}$ is a scheme,
each $Y_a$ is in fact a
scheme~\cite[\href{http://stacks.math.columbia.edu/tag/0BPW}{Tag 0BPW}]{stacks-project}, 
and thus $Y$ is a $\varpi$-adic formal scheme.
(In fact, it is easy to check directly that $U_{n(1)}$ acts freely
on $X_{n(a)}^{h,K'},$ and thus to see that $Y_a$ is an algebraic space.)
The natural morphism $Y \to \cC$ is then representable by algebraic
spaces;\ indeed, any morphism from an affine scheme to~$\cC$ factors
through some~$\cC^a$, and representability by algebraic spaces then follows from the
representability by algebraic spaces of~$Y_a\to\cC^a$, and the
Cartesianness of the diagram \[\xymatrix{Y_a\ar[r]\ar[d]&Y\ar[d]\\ \cC^a\ar[r]&\cC}\]
Similarly, the morphism $Y\to\cC$ is smooth and surjective, and so witnesses the claim that
$\cC$ is a $\varpi$-adic formal algebraic stack.

To check that $\cC$ has affine diagonal, it suffices to check
that each $\cC^a$ has affine diagonal, which follows from the fact that
$(\cG_{n(a)})_{\cO/\varpi^a}$ is in fact an affine group scheme
over $\cO/\varpi^a$ (indeed, as in~\cite[\S 2.b.1]{MR2562795}, it is a
Weil restriction of~$\GL_d$). 
\end{proof}

We next introduce the moduli stack of \'etale $\varphi$-modules, again both with
and without descent data.

\begin{defn}\label{defn: R^dd}
For each integer~$a\ge 1$, we let 
 $\cR^{\dd,a}_{d,K'}$ be the \emph{fppf} $\cO/\varpi^a$-stack which
  associates
  to any $\cO/\varpi^a$-algebra $A$ the groupoid $\cR^{\dd,a}_{d,K'}(A)$ of rank $d$ \'etale
  $\varphi$-modules  with $A$-coefficients and  descent data from $K'$ to
  $K$.

By~\cite[\href{http://stacks.math.columbia.edu/tag/04WV}{Tag 04WV}]{stacks-project}, 
we may also regard each of the stacks $\cR_{d,K'}^{\dd,a}$ as an {\em fppf}
$\cO$-stack, 
and we then write 
  $\cR^{\dd}:=\varinjlim_{a}\cR^{\dd,a}$,
  which is again an {\em fppf} $\cO$-stack.

  We will omit $d, K'$ from the
  notation wherever doing so will not cause confusion.
 We write $\cR$ and $\cR^a$ for $\cR^{\dd}$ and $\cR^{\dd,a}$ when the descent data is for the trivial extension $K'/K'$.
 \end{defn}

Just as in the case of $\cC^{\dd}$, the morphism $\cR^{\dd} \to \Spec \cO$ factors through
$\Spf \cO$, and for each $a \geq 1$, there is a natural isomorphism
$\cR^{\dd,a} \iso \cR^{\dd}\times_{\Spf \cO} \Spec \cO/\varpi^a.$
Thus each of the morphisms $\cR^{\dd,a} \to \cR^{\dd,a+1}$ is a thickening.

  There is a natural morphism $\cC_{d,h,K'}^{\dd}\to\cR_d^{\dd}$,
defined via
 $$ (\gM,\varphi,\{\hat{g}\}_{g\in\Gal(K'/K)}) \mapsto
(\gM[1/u],\varphi,\{\hat{g}\}_{g\in\Gal(K'/K)}),$$ and natural morphisms
$\cC^{\dd}\to\cC$ and $\cR^{\dd}\to\cR$ given by forgetting the descent
data.  
In the optic of Section~\ref{subsec: etale phi modules
  and Galois representations}, the stack~$\cR^{\dd}_d$ may
morally be thought of as a moduli of $G_{K_{\infty}}$-representations,
and the morphisms
$\cC^{\dd}_{d,h,K'}\to\cR^{\dd}_d$ correspond to passage from a Breuil--Kisin module to
its underlying Galois representation.

\begin{prop}
  \label{prop: relative representability of descent data for R} For
  each $a\ge 1$, the natural morphism
  $\cR^{\dd,a}\to\cR^a$ is representable by algebraic spaces,
  affine, and of finite presentation.
\end{prop}
\begin{proof} To see this, consider the pullback along some morphism $\Spec
  A\to\cR^a$ (where $A$ is a $\cO/\varpi^a$-algebra);\ we must show
  that given an \'etale $\varphi$-module $M$ of rank $d$ without descent data,
  the
data of giving 
additive bijections $\hat{g}:M\to M$, satisfying the further
  property that:
  \begin{itemize}
  \item the maps $\hat{g}$ commute with $\varphi$, satisfy
    $\hat{g_1}\circ\hat{g_2}=\widehat{g_1\circ g_2}$, and we have
    $\hat{g}(sm)=g(s)\hat{g}(m)$ for all $s\in\gS_A[1/u]$, $m\in M$
  \end{itemize}is represented by an affine algebraic space (i.e.\ an affine scheme!)
of finite presentation over~$A$.

To see this, note first that such maps $\hat{g}$ are by definition
$\gS_A^0[1/v]$-linear. 
The data of giving
an $\gS^0_A[1/v]$-linear automorphism of $M$ which commutes with~$\varphi$ is representable
by an affine scheme of finite presentation over~$A$ by~\cite[Prop.\ 5.4.8]{EGstacktheoreticimages} 
and so the
data of a finite collection of automorphisms is also representable by
a finitely presented affine scheme over $A$. The further commutation and composition
conditions on the $\hat{g}$ cut out a closed subscheme, as does the condition of
$\gS_A[1/u]$-semi-linearity, so the result follows.
\end{proof}

\begin{cor}
\label{cor:diagonal of R}
The diagonal of $\cR^{\dd}$ is representable by algebraic spaces, affine,
and of finite presentation.
\end{cor}
\begin{proof}
Since $\cR^{\dd} = \varinjlim_a\cR^{\dd,a} \iso \varinjlim_a \cR^{\dd}\times_{\Spf \cO}
\Spec \cO/\varpi^a$,
and since the transition morphisms are closed immersions (and hence monomorphisms),
we have a Cartesian diagram
$$\xymatrix{ \cR^{\dd,a} \ar[r]\ar[d] & \cR^{\dd,a}\times_{\cO/\varpi^a} \cR^{\dd,a} 
\ar[d] \\
\cR^\dd \ar[r] & \cR^{\dd} \times_{\cO} \cR^{\dd}}$$
for each~$a~\geq~1$,
and the diagonal morphism of $\cR^{\dd}$ is the inductive limit of the diagonal
morphisms of the various $\cR^{\dd,a}$.
Any morphism from an affine scheme $T$ to $\cR^{\dd} \times_{\cO} \cR^{\dd}$
thus factors through
one of the $\cR^{\dd,a}\times_{\cO/\varpi^a} \cR^{\dd,a}$, and the fibre product
$\cR^{\dd} \times_{\cR^{\dd} \times \cR^{\dd}} T$
may be identified with
$\cR^{\dd,a} \times_{\cR^{\dd,a} \times \cR^{\dd,a}} T$.
It is thus equivalent to prove that each of the diagonal morphisms
$\cR^{\dd,a} \to \cR^{\dd,a} \times_{\cO/\varpi^a} \cR^{\dd,a}$
is representable by algebraic spaces, affine, and of finite presentation.

The diagonal of $\cR^{\dd,a}$ may be obtained
by composing the pullback over $\cR^{\dd,a}\times_{\cO}
\cR^{\dd,a}$ of the diagonal $\cR^a \to \cR^a \times_{\cO} \cR^a$ with
the relative diagonal of the morphism $\cR^{\dd,a} \to \cR^a$.
The first of these morphisms is representable by algebraic spaces,
affine, and of finite presentation, by \cite[Thm.~5.4.11~(2)]{EGstacktheoreticimages},
and the second is also representable by algebraic
spaces, affine, and of finite presentation,
since it is the relative diagonal of a morphism which has these properties,
by~Proposition~\ref{prop: relative representability of descent data for R}.
\end{proof}

\begin{cor}
  \label{cor: C^dd,a is an algebraic stack}
\begin{enumerate}
\item
  For each $a\ge 1$,
  $\cC^{\dd,a}$ is an algebraic stack of finite presentation over
  $\Spec\cO/\varpi^a$, with affine diagonal.
\item
  The Ind-algebraic stack $\cC^{\dd} := \varinjlim_a \cC^{\dd,a}$ 
  is furthermore a $\varpi$-adic formal algebraic stack.
\item
The morphism $\cC^{\dd}_h \to \cR^{\dd}$ is representable by algebraic spaces and proper.
\end{enumerate}
\end{cor}
\begin{proof}By Proposition~\ref{prop: C without dd is a $p$-adic formal  algebraic
    stack}, $\cC^a$ is an algebraic stack of
finite type over $\Spec\cO/\varpi^a$ with affine diagonal.
In particular it has quasi-compact diagonal, and so is quasi-separated.
Since $\cO/\varpi^a$ is Noetherian, it follows from~\cite[\href{http://stacks.math.columbia.edu/tag/0DQJ}{Tag 0DQJ}]{stacks-project} that $\cC^a$ is 
in fact of finite presentation over $\Spec \cO/\varpi^a$.

By Proposition~\ref{prop: relative representability of descent data for
  R}, the  morphism $\cR^{\dd,a}\times_{\cR^a}\cC^a\to\cC^a$ is
representable by algebraic spaces and of finite presentation, so it follows from~\cite[\href{http://stacks.math.columbia.edu/tag/05UM}{Tag
    05UM}]{stacks-project} that $\cR^{\dd,a}\times_{\cR^a}\cC^a$ is
  an algebraic stack of finite presentation over $\Spec\cO/\varpi^a$. In order to show
that $\cC^{\dd,a}$ is an algebraic stack of finite presentation over $\Spec\cO/\varpi^a$,
it therefore suffices 
  to show 
  that the natural monomorphism 
\numequation
\label{eqn:closed}
\cC^{\dd,a}\to \cR^{\dd,a}\times_{\cR^a}\cC^a
\end{equation}
 is representable by
algebraic spaces and of finite presentation.  We will in fact show that it is
a closed immersion (in the sense that its pull-back under 
any morphism from a scheme to its target becomes a closed immersion
of schemes);\ since the target is locally Noetherian, and closed
immersions are automatically of finite type and quasi-separated, 
it follows
from~\cite[\href{http://stacks.math.columbia.edu/tag/0DQJ}{Tag
  0DQJ}]{stacks-project} that this closed immersion is of finite
presentation, as required.

By~\cite[\href{http://stacks.math.columbia.edu/tag/0420}{Tag
  0420}]{stacks-project}, the property of being a closed immersion can  be checked
after pulling back to an affine scheme, and then working
\emph{fpqc}-locally. 
The claim then follows easily from the proof of~\cite[Prop.\ 5.4.8]{EGstacktheoreticimages}, as
  \emph{fpqc}-locally the condition that a lattice in an \'etale
  $\varphi$-module of rank $d$ with descent data is preserved by the action
  of the $\hat{g}$ is determined by the vanishing of the coefficients of
  negative powers of~$u$ in a matrix.

To complete the proof of~(1), it suffices to show that the diagonal of
$\cC^{\dd,a}$ is affine.  Since (as we have shown)
the morphism~\eqref{eqn:closed} is a closed immersion, and thus a monomorphism,
it is equivalent to show that the diagonal 
of $\cR^{\dd,a}\times_{\cR^a} \cC^a$ is affine.
To ease notation, we denote this fibre product by $\cY$.
We may then factor the diagonal of $\cY$ as the composite 
of the pull-back over $\cY \times_{\cO/\varpi^a} \cY$
of the diagonal morphism $\cC^a \to \cC^a \times_{\cO/\varpi^a} \cC^a$
and the relative diagonal $\cY \to \cY\times_{\cC^a} \cY$.
The former morphism is affine, by~\cite[Thm.~5.4.9~(1)]{EGstacktheoreticimages},
and the latter morphism is also affine, since it is the pullback via $\cC^a \to
\cR^a$ 
of the relative diagonal morphism $\cR^{\dd,a} \to \cR^{\dd,a}\times_{\cR^a} \cR^{\dd,a}$,
which is affine (as already observed in the proof of Corollary~\ref{cor:diagonal of R}).

To prove~(2),
consider the morphism $\cC^{\dd} \to \cC.$  
This is a morphism of $\varpi$-adic Ind-algebraic stacks,
and by what we have already proved, it is representable 
by algebraic spaces.    Since the target is a $\varpi$-adic formal
algebraic stack, it follows from ~\cite[Lem.\
7.9]{Emertonformalstacks} that the source is also a $\varpi$-adic formal
algebraic stack, as required. 

To prove~(3), since each of $\cC^{\dd,a}$ and $\cR^{\dd,a}$ is obtained from
$\cC^{\dd}$ and $\cR^{\dd}$ via pull-back over $\cO/\varpi^a$,
it suffices to prove that each of the morphisms $\cC^{\dd,a} \to \cR^{\dd,a}$
is representable by algebraic spaces and proper.   Each of these morphisms
factors as
$$\cC^{\dd,a} \buildrel \text{\eqref{eqn:closed}} \over
\longrightarrow \cR^{\dd,a} \times_{\cR^a} \cC^a \buildrel \text{proj.} \over
\longrightarrow  \cR^{\dd,a}.$$
We have already shown that the first of these morphisms is a closed immersion,
and hence representable by algebraic spaces and proper.  The second morphism is also 
representable by algebraic spaces and proper, since it is a base-change 
of the morphism $\cC^a \to \cR^a$, which has these properties
by~\cite[Thm.~5.4.11~(1)]{EGstacktheoreticimages}.
\end{proof}

The next lemma gives a concrete interpretation of the points of $\cC^{\dd}$ over 
$\varpi$-adically complete $\cO$-algebras, extending the tautological interpretation
of the points of each $\cC^{\dd,a}$ prescribed by Definition~\ref{defn: C^dd,a }.

\begin{lemma}\label{lem:points-of-Cdd}
 If $A$ is a $\varpi$-adically complete $\cO$-algebra then the
 $\Spf(A)$-points of $\cC^{\dd}$ are the Breuil--Kisin modules
 of rank $d$ and height $h$ with $A$-coefficients and descent
 data. 
\end{lemma} 

\begin{proof} 

Let $\gM$ be a Breuil--Kisin module of rank $d$ and height $h$ with
$A$-coefficients and descent data. Then the sequence
$\{\gM/\varpi^a\gM\}_{a\ge 1}$ defines a $\Spf(A)$-point of $\cC^{\dd}$ 
 (\emph{cf.}\ Remark~\ref{rem:Kisin-mod-ideal}), and since $\gM$ is
 $\varpi$-adically complete  it is recoverable from the sequence
 $\{\gM/\varpi^a\gM\}_{a \ge 1}$. 

In the other direction, suppose that $\{\gM_a\}$ is a $\Spf(A)$-point
of $\cC^{\dd}$, so that $\gM_a \in \cC^{\dd,a}(A/\varpi^a)$. Define
$\gM = \invlim_{a} \gM_a$, and similarly define $\varphi_\gM$ and
$\{\ghat\}$ as inverse limits. Observe that $\varphi^* \gM = \invlim_a
\varphi^*\gM_a$ (since $\varphi : \gS_A \to \gS_A$ makes $\gS_A$ into
a free $\gS_A$-module). Since each $\Phi_{\gM_a}$ is injective with
cokernel killed by $E(u)^h$ the same holds for $\Phi_\gM$. 

Since the required properties of the descent data are immediate, to complete the proof it remains to check that $\gM$ is a projective
$\gS_A$-module (necessarily of rank~$d$, since its rank will equal
that of $\gM_1$), which is a consequence of \cite[Prop.\ 0.7.2.10(ii)]{MR3075000}.
\end{proof}

We now temporarily reintroduce~$h$ to the notation. 

\begin{df} For each $h\ge 0$,  write~$\cR^a_h$ for the
scheme-theoretic image of $\cC^a_h\to\cR^a$ in the sense of~\cite[Defn.\
3.2.8]{EGstacktheoreticimages};\ then by~\cite[Thms.\ 5.4.19,
5.4.20]{EGstacktheoreticimages}, $\cR^a_h$ is an algebraic stack of
finite presentation over~$\Spec\cO/\varpi^a$, the morphism 
$\cC^a_h\to\cR^a$ factors through $\cR^a_h$, 
and we may write
$\cR^a \cong \varinjlim_h\cR^a_h$ as an inductive limit of closed substacks,
the natural transition morphisms being closed immersions. 

We similarly write $\cR^{\dd,a}_h$ for the scheme-theoretic image of the morphism
$\cC_h^{\dd,a}\to\cR^{\dd,a}$ in the sense of~\cite[Defn.\
3.2.8]{EGstacktheoreticimages}.
\end{df}

  \begin{thm}
    \label{thm: R^dd is an ind algebraic stack}For each $a\ge 1$,
    $\cR^{\dd,a}$ is an Ind-algebraic stack. Indeed, we can write
    $\cR^{\dd,a}=\varinjlim_h\cX_h$ as an inductive limit of
    algebraic stacks of finite presentation over $\Spec\cO/\varpi^a$, the
    transition morphisms being closed immersions.
  \end{thm}
  \begin{proof}
    As we have just recalled, by~\cite[Thm.\
    5.4.20]{EGstacktheoreticimages} 
    we can write
    $\cR^{a}=\varinjlim_h\cR^{a}_h$, so that if we set
    $\cX^{\dd,a}_h:=\cR^{\dd,a}\times_{\cR^a}\cR^a_h$, then
    $\cR^{\dd,a}=\varinjlim_h\cX^{\dd,a}_h$, and the transition morphisms are
    closed immersions.  Since~$\cR^a_h$ is of finite presentation
    over~$\Spec\cO/\varpi^a$, and a composite of morphisms of finite
    presentation is of finite presentation, it follows from
    Proposition~\ref{prop: relative representability of descent data
      for R} and
    \cite[\href{http://stacks.math.columbia.edu/tag/05UM}{Tag
      05UM}]{stacks-project} that $\cX^{\dd,a}_h$ is an algebraic stack of
    finite presentation over $\Spec\cO/\varpi^a$, as required.
  \end{proof}

\begin{thm}
  \label{thm: R^dd_h is an algebraic stack}$\cR^{\dd,a}_h$ is an algebraic stack of finite presentation over
    $\Spec\cO/\varpi^a$. It is a closed substack of $\cR^{\dd,a}$, and the
    morphism $\cC_h^{\dd,a}\to\cR^{\dd,a}$ factors through
a morphism 
    $\cC_h^{\dd,a} \to \cR^{\dd,a}_h$ which is representable by algebraic
spaces, proper, and scheme-theoretically dominant in the sense of \cite[Def.~3.1.3]{EGstacktheoreticimages}. 
\end{thm}
\begin{proof}
As in the proof of Theorem~\ref{thm: R^dd is an ind algebraic stack},
if we set $\cX^{\dd,a}_h:=\cR^{\dd,a}\times_{\cR^a}\cR^a_h$, then~$\cX^{\dd,a}_h$ is an algebraic stack of
    finite presentation over $\Spec\cO/\varpi^a$, and the natural
    morphism $\cX^{\dd,a}_h\to\cR^{\dd,a}$ is a closed immersion. The
    morphism~$\cC_h^{\dd,a}\to\cR^{\dd,a}$ factors through~$\cX^{\dd,a}_h$
    (because the morphism $\cC_h^a\to\cR^{a}$ factors through its
    scheme-theoretic image~$\cR^a_h$), so by~\cite[Prop.\
    3.2.31]{EGstacktheoreticimages}, 
    $\cR^{\dd,a}_h$ is the
    scheme-theoretic image of the morphism of algebraic stacks
    $\cC_h^{\dd,a}\to\cX^{\dd,a}_h$. The required properties now follow
    from \cite[Lem.\ 3.2.29]{EGstacktheoreticimages}
    (using representability by algebraic spaces and properness of the morphism
$\cC_h^{\dd,a} \to \cR^{\dd,a}$,
as proved in Corollary~\ref{cor: C^dd,a is an algebraic stack}~(3),
to see that the induced morphism $\cC_h^{\dd,a} \to \cR_h^{\dd,a}$ is representable
by algebraic spaces and proper,
along with~\cite[\href{http://stacks.math.columbia.edu/tag/0DQJ}{Tag
      0DQJ}]{stacks-project}, and the fact that~$\cX^{\dd,a}_h$ is of finite presentation
    over~$\Spec\cO/\varpi^a$, to see that~$\cR^{\dd,a}_h$ is of finite presentation). 
\end{proof}

\subsection{Representations of tame groups}\label{subsec: tame groups}
Let $G$ be a finite group.

\begin{df}
We let $\Rep_d(G)$ denote the algebraic stack
classifying $d$-dimensional
representations of $G$ over $\cO$: if $X$ is any $\cO$-scheme,
then $\Rep_d(G)(X)$ is the groupoid consisting 
of locally free sheaves of rank $d$ over $X$ endowed with an $\cO_X$-linear
action of $G$ (rank $d$ locally free $G$-sheaves, for short);
morphisms are $G$-equivariant isomorphisms of vector bundles.
\end{df}

We now suppose that $G$ is tame, i.e.\ that it has prime-to-$p$ order.
In this case (taking into account the fact that $\F$ has characteristic~$p$,
and that $\cO$ is Henselian),
the isomorphism classes of $d$-dimensional $G$-representations
of $G$ over $E$ and over $\F$ are in natural bijection.  Indeed, any 
finite-dimensional representation $\tau$
of $G$ over $E$ contains a $G$-invariant
$\cO$-lattice $\tau^{\circ}$,
and the associated representation of $G$ over $\F$ is given
by forming $\overline{\tau}:= \F\otimes_{\cO} \tau^{\circ}$.

\begin{lemma}
\label{lem:tame reps}
Suppose that~ $G$ is tame, and that~ $E$ is chosen large enough so
 that each irreducible
representation of $G$ over $E$ is absolutely irreducible {\em (}or,
equivalently, so that each irreducible representation of $G$ over $\F$
is absolutely irreducible{\em )}, and so that each irreducible
representation of~$G$ over~$\Qpbar$ is defined over~$E$
\emph{(}equivalently, so that each irreducible representation of~$G$
over~$\Fpbar$ is defined over~$\F$\emph{)}. 
\begin{enumerate}
\item $\Rep_d(G)$ is the disjoint union of irreducible
  components 
$\Rep_d(G)_{\tau}$, where $\tau$ ranges over the finite set of isomorphism
classes of $d$-dimensional representations of $G$ over $E$.
\item A morphism $X \to \Rep_d(G)$
factors through $\Rep_d(G)_{\tau}$
if and only if the associated locally free $G$-sheaf on $X$
is Zariski locally isomorphic to $\tau^{\circ}\otimes_{\cO} \cO_X$.
\item If we write $G_{\tau} := \Aut_{\cO[G]}(\tau^{\circ})$, 
then $G_{\tau}$ is a smooth
\emph{(}indeed reductive\emph{)} group scheme over $\cO$,
and
$\Rep_d(G)_{\tau}$ is isomorphic to the classifying space
$[\Spec \cO/G_{\tau}]$. 
\end{enumerate}
\end{lemma} 
\begin{proof}
Since~$G$ has order prime to~$p$, the
representation~$P:= \oplus_\sigma\sigma^\circ$ is a projective generator of
the category of $\cO[G]$-modules, where~$\sigma$ runs over a set of
representatives for the isomorphism classes of irreducible
$E$-representations of~$G$. (Indeed, each~$\sigma^\circ$ is
projective, because the fact that~$G$ has order prime to~$p$ means
that all of the~$\Ext^1$s against~$\sigma^\circ$ vanish. To see that
$\oplus_\sigma\sigma^\circ$ is a generator, we need to show that every
$\cO[G]$-module admits a non-zero map from some~$\sigma^\circ$. We can
reduce to the case of a finitely generated module~$M$, and it is
then enough (by projectivity) to prove that~$M\otimes_\cO\F$ admits
such a map, which is clear.)
Our assumption that each $\sigma$ is absolutely irreducible furthermore
shows that $\End_G(\sigma^{\circ}) = \cO$ for each $\sigma$,
so that $\End_G(P) = \prod_{\sigma} \cO$.  

Standard Morita theory then shows that the functor
$M \mapsto \Hom_G(P,M)$ induces an equivalence between the category
of $\cO[G]$-modules and the category of $\prod_{\sigma} \cO$-modules.
Of course, a $\prod_{\sigma}\cO$-module is just given by a tuple
$(N_{\sigma})_{\sigma}$ of $\cO$-modules,
and in this optic, the functor $\Hom_G(P,\text{--})$ 
can be written as $M \mapsto \bigl(\Hom_G(\sigma^{\circ},M)\bigr)_{\sigma}$,
with a quasi-inverse functor being given by $(N_{\sigma}) \mapsto
\bigoplus_{\sigma} \sigma^{\circ} \otimes_{\cO} N_{\sigma}.$
It is easily seen (just using the fact that $\Hom_G(P,\text{--})$
induces an equivalence of categories)
that $M$ is a finitely generated projective $A$-module,
for some $\cO$-algebra $A$,
if and only if each $\Hom_G(\sigma^{\circ},M)$
is a finitely generated projective $A$-module.

The preceding discussion shows that
giving a rank $d$ representation of $G$ over an $\cO$-algebra $A$
amounts to giving a tuple $(N_{\sigma})_{\sigma}$ of projective
$A$-modules, of ranks~$n_{\sigma}$, such that $\sum_{\sigma} n_{\sigma}
\dim \sigma = d.$
For each such tuple of ranks $(n_\sigma)$, we obtain a corresponding
moduli stack
$\Rep_{(n_{\sigma})}(G)$ classifying rank $d$ representations of $G$
which decompose in this manner,
and $\Rep_d(G)$ is isomorphic to the disjoint union
of the various stacks~$\Rep_{(n_{\sigma})}(G).$

If we write $\tau = \oplus_{\sigma} \sigma^{n_{\sigma}}$, 
then we may relabel $\Rep_{(n_{\sigma})}(G)$ as~$\Rep_{\tau}(G)$;\ 
statements~(1) and~(2) are then proved.
By construction, there is an isomorphism  
$$\Rep_{\tau}(G) = \Rep_{(n_{\sigma})}(G) \iso \prod_{\sigma}
[\Spec \cO/\GL_{n_{\sigma}}].$$
Noting that $G_{\tau} := \Aut(\tau) = 
\prod_\sigma\GL_{n_\sigma/\cO}$, 
we find that statement~(3) follows as well.
\end{proof}

For each~$\tau$, it follows from the identification of $\Rep_d(G)_{\tau}$ with
$[\Spec \cO/G_{\tau}]$ that there is a natural map
$\Rep_d(G)_{\tau}\to\Spec\cO$. We let~$\pi_0(\Rep_{d}(G))$ denote the
disjoint union of copies of~$\Spec\cO$, one for each isomorphism
class~$\tau$;\ then there is a natural map
$\Rep_d(G)\to\pi_0(\Rep_d(G))$. While we do not want to develop a
general theory of the \'etale~$\pi_0$ groups of algebraic stacks, we
note that it is natural to regard~$\pi_0(\Rep_d(G))$ as the
\'etale~$\pi_0$ of~$\Rep_d(G)$.
\subsection{Tame inertial types}
\label{subsec:tame inertial types,
  general height}Write $I(K'/K)$ for the inertia subgroup of~$\Gal(K'/K)$.
Since we are assuming that~$E$ is large enough that it contains the
image of every embedding $K'\into\Qpbar$, it follows in particular that every
$\Qpbar$-character of~$I(K'/K)$ is defined
over~$E$.

Recall from Subsection~\ref{subsec: notation} that if $A$ is an~$\cO$-algebra, and~$\gM$ is a Breuil--Kisin module with
$A$-coefficients, then we write $\gM_i$ for  $e_i\gM$. Since $I(K'/K)$
acts trivially on~$W(k')$, the $\hat{g}$ for~$g\in I(K'/K)$ stabilise
each~$\gM_i$, inducing an action of $I(K'/K)$ on $\gM_i/u\gM_i$. 

Write \[\Rep_{d,I(K'/K)}:=\prod_{i = 0}^{f'-1}
\Rep_d\bigl(I(K'/K)\bigr),\]the fibre product being taken over~$\cO$.  If~$\{\tau_i\}$ is an $f'$-tuple of isomorphism classes 
of $d$-dimensional representations of $I(K'/K)$, we write \[\Rep_{d,I(K'/K),\{\tau_i\}}:=\prod_{i = 0}^{f'-1}
\Rep_d\bigl(I(K'/K)\bigr)_{\tau_i}.\] Lemma~\ref{lem:tame reps} shows that we may write
$$\Rep_{d,I(K'/K)} = \coprod_{\{\tau_i\}}
\Rep_{d,I(K'/K),\tau_i}.$$
 Note that since $K'/K$ is tamely ramified,
$I(K'/K)$ is abelian of prime-to-$p$ order, and each~$\tau_i$ is just
a sum of characters. If all of the~$\tau_i$ are equal to some
fixed~$\tau$, then we write $\Rep_{d,I(K'/K),\tau}$ for
$\Rep_{d,I(K'/K),\{\tau_i\}}$. We have corresponding
stacks~$\pi_0(\Rep_{d,I(K'/K)})$,
~$\pi_0(\Rep_{d,I(K'/K),\{\tau_i\}})$
and~$\pi_0(\Rep_{d,I(K'/K),\tau})$, defined in the obvious way.

 If $\gM$ is a Breuil--Kisin module of rank~$d$ with descent data and $A$-coefficients,
then $\gM_i/u\gM_i$ is projective $A$-module of rank $d$, endowed
with an $A$-linear action of $I(K'/K)$,
and so is an $A$-valued point of $\Rep_d\bigl(I(K'/K)\bigr).$ Thus we obtain a morphism 
\numequation
\label{eqn:mixed type morphism}
\cC_d^{\dd} \to \Rep_{d,I(K'/K)},
\end{equation}
defined via $\gM \mapsto ( \gM_0/u\gM_0, \ldots, \gM_{f'-1}/u\gM_{f'-1}).$

\begin{defn}
  \label{defn: Kisin module of type tau general rank}Let $A$ be
  an~$\cO$-algebra, and let $\gM$ be a 
  Breuil--Kisin module of rank~$d$ with $A$-coefficients. We say that $\gM$ has
  \emph{mixed type}~$(\tau_i)_i$
  if the composite $\Spec A \to \cC_d^{\dd} \to \Rep_{d,I(K'/K)}$ (the first arrow being the morphism
that classifies $\gM$, and the second arrow being~(\ref{eqn:mixed 
type morphism}))
factors through $\Rep_{d,I(K'/K),\{\tau_i\}}$.
Concretely, this is equivalent to requiring that,
Zariski locally on
  $\Spec A$, 
there is an $I(K'/K)$-equivariant isomorphism $\gM_i/u\gM_i \cong
A\otimes_{\cO}\tau_i^{\circ}$ for each $i$.

If each $\tau_i=\tau$ for some fixed~$\tau$, then we say that the type
of $\gM$ is unmixed, or simply that
$\gM$ has \emph{type}~$\tau$.
\end{defn}
\begin{rem}
  \label{rem: Kisin modules don't always have a type} If $A=\cO$ then
  a Breuil--Kisin module necessarily has some (unmixed) type~$\tau$, since
after 
  inverting $E(u)$ and reducing modulo $u$ the map   $\Phi_{\gM,i}$
  gives an
  $I(K'/K)$-equivariant $E$-vector space isomorphism
  $\varphi^*(\gM_{i-1}/u\gM_{i-1})[\frac 1p] \toisom (\gM_{i}/u\gM_{i})[\frac 1p]$.
  However if $A=\F$
  there are Breuil--Kisin modules which have a genuinely mixed type;\ indeed, it is easy to write down examples of free Breuil--Kisin
  modules of rank one of any mixed type (see
  also~\cite[Rem. 3.7]{2015arXiv151007503C}), which necessarily cannot
  lift to characteristic zero. This shows that~$\cC_d^\dd$ is not flat
  over~$\cO$. In the following sections, when $d=2$ and $h=1$ we define a
  closed substack~$\cC_d^{\dd,\BT}$ of~$\cC^{\dd}$ which \emph{is} flat
  over~$\cO$, and can be thought of as taking the Zariski closure of
   $\Qpbar$-valued Galois representations that become Barsotti--Tate
   over $K'$ and such that all pairs of labeled Hodge--Tate
   weights are $\{0,1\}$ (see Remark~\ref{rem: strong determinant condition motivated by
      flat closure} below). 
\end{rem}

\begin{defn} 
  \label{defn: Ctau in general}Let $\cC_d^{(\tau_i)}$ be the
  \'etale substack of $\cC_d^{\dd}$ which associates to
  each $\cO$-algebra $A$ the subgroupoid $\cC_d^{(\tau_i)}(A)$ of
  $\cC_d^{\dd}(A)$ consisting of those Breuil--Kisin modules which are of
  mixed type~$(\tau_i)$. If each $\tau_i = \tau$ for some fixed
  $\tau$, we write $\cC_d^\tau$ for $\cC_d^{(\tau_i)}$.
\end{defn}

\begin{prop}
  \label{prop: Ctau is an algebraic stack}Each $\cC_d^{(\tau_i)}$ is an open
and closed
  substack of $\cC_d^\dd$, and $\cC_d^\dd$ is the disjoint union of its substacks~$\cC_d^{(\tau_i)}$. 
\end{prop}
\begin{proof}  
By Lemma~\ref{lem:tame reps}, $\Rep_{d,I(K'/K)}$ is the disjoint
union of its open and closed substacks $\Rep_{d,I(K'/K),\{\tau_i\}}$.
By definition $\cC_d^{(\tau_i)}$ is the preimage of 
$\Rep_{d,I(K'/K),\{\tau_i\}}$
under the morphism~(\ref{eqn:mixed type morphism});
the lemma follows.
\end{proof}

\section{Local geometry of Breuil--Kisin moduli stacks}
\label{sec:local-geom}

Our aim in this section is to  use 
local models of Shimura varieties at hyperspecial and Iwahori level to 
study the geometry of the moduli stacks of  Breuil--Kisin modules
described in Section~\ref{sec:moduli stacks}, in the case where the
height $h$ is equal to $1$. 

We begin in Section~\ref{subsec: local models
  generalities} by explaining how to construct a smooth morphism from
the stack $\cC_{1,K'}^{\dd}$ (and thus  its closed  substacks
$\cC^{(\tau_i)}$) to a suitable local model stack. When $d=2$, the main case
of interest in this paper,  we want to consider stacks
of Breuil--Kisin modules
of whose $\Spf(\cO)$-points
correspond to crystalline representations that
are not just of height $1$, but are actually potentially
Barsotti--Tate (i.e., whose Hodge--Tate weights are regular). This is done
by imposing Kottwitz-type determinant conditions, which we discuss in Section~\ref{subsec: local
  models rank two}. The resulting stacks will be denoted
$\cC^{\tau,\BT}$. 

In Section~\ref{subsec: explicit local models} we explain how to
relate our local models to the naive local models of~\cite{MR3135437}
at hyperspecial and Iwahori level, and finally in Section~\ref{subsec: local models
  relation to explicit PR ones} we deduce our main results on the
local geometry of the stacks $\cC^{\tau,\BT}$ by comparison with
known results on the geometry of these local models. Section~\ref{subsec: change of extension K
  prime over K} contains a brief interlude on the behaviour of our
constructions under change of base, which will allow us  to unify the
discussion of principal series and cuspidal types in Section~\ref{sec:galois-moduli}.

Finally, suppose that  $d=2$ and that the tame  type $\tau = \eta
\oplus \eta'$ has $\eta \neq \eta'$. In Section~\ref{subsec: Dieudonne
  stack} we construct a morphism from $\cC^{\tau,\BT}$ to an auxiliary
\emph{Dieudonn\'e stack} $\cD_{\eta}$. This construction is not needed
elsewhere in the present paper, but will be used
in \cite{cegsA,cegsC}. 

\subsection{Local models:\ generalities}\label{subsec: local models
  generalities}Throughout this section we allow~$d$ to be arbitrary;\ in Section~\ref{subsec: local
  models rank two} we specialise to the case~$d=2$, where we relate
the local models considered in this section to the local models
considered in the theory of Shimura varieties. We will usually
omit~$d$ from our notation, writing for example~$\cC$ for~$\cC_d$,
without any further comment. We begin with the
following lemma, for which we allow~$h$ to be arbitrary.

\begin{lem}
  \label{lem:cokernel is projective} Let $\gM$ be a rank~$d$ Breuil--Kisin
  module of height~$h$ with
  descent data over an $\cO$-algebra $A$. Assume further 
  either that~$A$ is $p^n$-torsion for some~$n$, or that~$A$ is Noetherian
  and $p$-adically complete. 
Then
  $\im \Phi_{\gM}/E(u)^h\gM$ is a finite projective
  $A$-module, and is a direct summand of~$\gM/E(u)^h\gM$ as an
  $A$-module. 
\end{lem}
\begin{proof}
  We follow the proof of~\cite[Lem.\ 1.2.2]{kis04}. We have a short
  exact sequence
  \[0\to \im \Phi_{\gM}/E(u)^h\gM\to
    \gM/E(u)^h\gM\to\gM/\im \Phi_{\gM}\to 0 \] in which the
  second term is a finite projective $A$-module (since it is a finite
  projective $\cO_{K'}\otimes_{\Zp}A$-module), so it is enough to show
  that the third term is a projective $A$-module. It is therefore enough to show
  that the finitely generated $A$-module $\gM/\im \Phi_{\gM}$
 is finitely presented and
  flat.

To see that it is finitely presented, note that we have the
equality \[\gM/\im \Phi_{\gM}=(\gM/E(u)^h)/(\im \Phi_{\gM}/E(u)^h), \]
and the right hand side admits a presentation
by finitely generated projective
$A$-modules  \[\varphi^*(\gM/E(u)^h)\to\gM/E(u)^h\to(\gM/E(u)^h) /(\im \Phi_{\gM}/E(u)^h)\to 0. \]

To see that it is flat, it is enough to show that for every finitely generated ideal~$I$
of~$A$, the map \[I\otimes_A\gM/\im \Phi_{\gM} \to
  \gM/\im \Phi_{\gM}\] is injective. It follows easily (for
example, from the snake lemma)
that it is enough to check that the complex \[0\to
  \varphi^*\gM\to \gM\to \gM/\im \Phi_{\gM}\to 0, \] 
which is exact by Lemma~\ref{lem:kisin injective},
remains
exact after tensoring with~$A/I$.  Since~$I$ is finitely generated we
have $\gM \otimes_A A/I \cong \gM \cotimes_A A/I$ by
Remark~\ref{rem:Kisin-mod-ideal}, and the desired exactness amounts
 to the injectivity of
$\Phi_{\gM \cotimes_A A/I}$ for the Breuil--Kisin module $\gM \cotimes_A A/I$.
This  follows immediately  from
  Lemma~\ref{lem:kisin injective} in the case that~$A$ is killed by
  $p^n$, and otherwise follows from the same lemma once we check that $A/I$ is
  $p$-adically complete, which 
follows from the
Artin--Rees lemma (as $A$ is assumed Noetherian and $p$-adically
complete). 
\end{proof}
We assume from now on that~$h=1$, but 
we continue to allow arbitrary~$d$. We allow~$K'/K$ to be any 
Galois extension such that~$[K':K]$ is prime to~$p$ (so in particular~$K'/K$
is tame). 

\begin{defn}\label{defn: }
  We let~$\cM_{\loc}^{K'/K,a}$ be the algebraic stack of finite presentation
over~$\Spec\cO/\varpi^a$ defined as follows:
  if~$A$ is an $\cO/\varpi^a$-algebra, then $\cM_{\loc}^{K'/K}(A)$ is
  the groupoid of tuples~$(\gL,\gL^+)$, where:
  \begin{itemize}
  \item $\gL$ is a rank $d$ projective $\cO_{K'}\otimes_{\Zp} A$-module,
    with a $\Gal(K'/K)$-semilinear, $A$-linear action of~$\Gal(K'/K)$;
  \item $\gL^+$ is an $\cO_{K'}\otimes_{\Zp} A$-submodule of~$\gL$, which is
    locally on~$\Spec A$ a direct summand of~$\gL$ as an $A$-module
    (or equivalently, for which $\gL/\gL^+$ is  projective as an $A$-module),
    and is preserved by~$\Gal(K'/K)$.
  \end{itemize}
  We set~$\cM_{\loc}^{K'/K}:=\varinjlim_a\cM_{\loc}^{K'/K,a}$, so
  that~$\cM_{\loc}^{K'/K}$ is a $\varpi$-adic formal algebraic stack,
of finite presentation over $\Spf \cO$
  (indeed, it is easily seen to be the $\varpi$-adic completion of an algebraic stack
of finite presentation over $\Spec \cO$).
\end{defn} 

In Section~\ref{subsec: explicit local models} we will explain how to relate these local model stacks (or more precisely, their refinements satisfying a certain determinant condition  that will be developed in Section~\ref{subsec: local
  models rank two}) to local models of Shimura varieties. Lemma~\ref{lem:cokernel is projective} shows that they are also connected to our stacks of Breuil--Kisin modules, as follows.

\begin{defn}By Lemma~\ref{lem:cokernel is projective}, we have a natural morphism
  $\Psi:\cC_{1,K'}^{\dd}\to\cM_{\loc}^{K'/K}$, which takes a Breuil--Kisin module
  with descent data~$\gM$ of height~$1$ to the pair\[(\gM/E(u)\gM, \im\Phi_{\gM}/E(u)\gM).\]
\end{defn}

\begin{remark}\label{rem:Mloc-is-too-big}
  The definition of the stack $\cM_{\loc}^{K'/K,a}$ does not include
  any condition that mirrors the commutativity between the Frobenius
  and the   descent data on a Breuil--Kisin module,   and so in general the
  morphism $\Psi_A:\cC_{1,K'}^{\dd}(A)\to\cM_{\loc}^{K'/K}(A)$ cannot be essentially surjective.
\end{remark}

It will be convenient to consider
the twisted group rings~$\gS_A[\Gal(K'/K)]$
and~$(\cO_{K'}\otimes_{\Zp}A)[\Gal(K'/K)]$, 
in which the elements $g\in\Gal(K'/K)$ obey the following
commutation relation with elements $s\in \gS_A$ (resp.\
$s\in\cO_{K'}$):
\[g\cdot s=g(s)\cdot g. \]
(In the literature these twisted group rings are more often written as
$\gS_A * \Gal(K'/K)$, $(\cO_{K'}\otimes_{\Zp}A) * \Gal(K'/K)$, in order
to distinguish them from the usual (untwisted) group rings, but as we
will only use the twisted versions in this paper, we prefer to use
this notation for them.)

By definition, endowing a finitely generated ~$\gS_A$-module~$P$ with a semilinear
$\Gal(K'/K)$-action is equivalent to giving it the structure of a
left~$\gS_A[\Gal(K'/K)]$-module. If~$P$ is projective as an
$\gS_A$-module, then it is also projective as
an $\gS_A[\Gal(K'/K)]$-module. Indeed, $\gS_A$ is a direct summand
of~$\gS_A[\Gal(K'/K)]$ as a $\gS_A[\Gal(K'/K)]$-module, given by the
central idempotent~$\frac{1}{|\Gal(K'/K)|}\sum_{g\in\Gal(K'/K)}g$,
so~$P$ is a direct summand of the projective module
$\gS_A[\Gal(K'/K)]\otimes_{\gS_A}P$. Similar remarks apply to the case
of  $(\cO_{K'}\otimes_{\Zp}A)[\Gal(K'/K)]$-modules.

\begin{thm}\label{thm: map from Kisin to general local model is
    formally smooth}The
  morphism~$\Psi:\cC_{1,K'}^{\dd}\to\cM_{\loc}^{K'/K}$ is representable by
algebraic spaces and smooth.  
\end{thm} 
\begin{proof}
We first show that the morphism~$\Psi$ is formally smooth,
in the sense that it satisfies the infinitesimal lifting criterion for nilpotent
thickenings of affine test objects \cite[Defn.\
2.4.2]{EGstacktheoreticimages}.
For this, 
we follow the proof of~\cite[Prop.\ 2.2.11]{kis04} (see also the proof
  of~\cite[Thm.\ 4.9]{2015arXiv151007503C}). Let $A$ be an $\cO/\varpi^a$-algebra and $I\subset A$ be a nilpotent ideal. 
  Suppose that we are given  $\gM_{A/I}\in \cC_{1,K'}^{\dd}(A/I)$
  and $(\gL_A,\gL^+_A)\in \cM_{\loc}^{K'/K}(A)$ 
together with an isomorphism
  \[
  \Psi(\gM_{A/I})\isoto (\gL_{A},\gL^+_{A})\otimes_{A}A/I=:(\gL_{A/I},\gL^+_{A/I}).
  \]
  We must show that there exists
  $\gM_{A}\in \cC_{1,K'}^{\dd}(A)$ together with an isomorphism 
  $\Psi(\gM_{A})\isoto
  (\gL_{A},\gL^+_{A})$ lifting the
  given isomorphism. 
  
As explained above, we can and do think of~$\gM_{A/I}$ as a finite
projective $\gS_{A/I}[\Gal(K'/K)]$-module, and~$\gL_A$ as a finite
projective $(\cO_{K'}\otimes_{\Zp}A)[\Gal(K'/K)]$-module. Since the
closed 2-sided ideal 
of~$\gS_{A}[\Gal(K'/K)]$ generated by~$I$ consists of
nilpotent elements, we may lift
$\gM_{A/I}$ to a finite
projective~$\gS_{A}[\Gal(K'/K)]$-module~$\gM_A$. This is presumably
standard, but for lack of a reference we indicate a proof. In fact the
proof of ~\cite[\href{https://stacks.math.columbia.edu/tag/0D47}{Tag
  0D47}]{stacks-project} goes over unchanged to our setting. Writing
~$\gM_{A/I}$ as a direct summand of a finite free
$\gS_{A/I}[\Gal(K'/K)]$-module $F$, it is enough to lift the
corresponding idempotent in~$\End_{\gS_{A/I}[\Gal(K'/K)]}(F)$, which
is possible
by~\cite[\href{http://stacks.math.columbia.edu/tag/05BU}{Tag
  05BU}]{stacks-project}. (See also~\cite[Thm.\ 21.28]{MR1125071} for another proof of the existence
of lifts of idempotents in this generality.) Note that since $\gM_{A/I}$ is of
rank~$d$ as a projective $\gS_{A/I}$-module, $\gM_A$ is of rank~$d$ as
a projective~$\gS_A$-module.

Since~$\gM_A/E(u)\gM_A$ is a projective
$(\cO_{K'}\otimes_{\Zp}A)[\Gal(K'/K)]$-module, we may lift the
composite \[\gM_A/E(u)\gM_A \onto \gM_{A/I}/E(u)\gM_{A/I} \isoto
\gL_{A/I}\] to a morphism $\theta : \gM_A/E(u)\gM_A \to \gL_A$.  
Since the composite of
$\theta$ with $\gL_A \onto \gL_{A/I}$ is
surjective, it follows by Nakayama's lemma that $\theta$ is
surjective. But a surjective map of projective modules of the same
rank is an isomorphism, and so $\theta$ is
an isomorphism lifting the given isomorphism $\gM_{A/I}/E(u)\gM_{A/I} \isoto
\gL_{A/I}$. 

We let~$\gM_A^+$ denote the
preimage in~$\gM_A$ of~$\theta^{-1}(\gL^+_A)$. The image of the induced
map $f : \gM_A^+ \subset \gM_A \onto \gM_A \cotimes_A A/I \cong  \gM_{A/I}$ is precisely $\im
\Phi_{\gM_{A/I}}$, since the same is true modulo $E(u)$ and because
$\gM^+_A$, $\im \Phi_{\gM_{A/I}}$
contain $E(u)\gM_A$, $E(u) \gM_{A/I}$ respectively. Observing that 
\[ \gM_A / \gM_A^+ \cong (\gM_A/E(u)\gM_A)/(\gM_A^+/E(u)\gM_A) \cong
\gL_A/\gL^+_A \] 
we deduce that $\gM_A/\gM_A^+$ is projective as an $A$-module, whence
$\gM_A^+$ is an $A$-module direct summand of $\gM_A$. By the same
argument $\im \Phi_{\gM_{A/I}}$ is an $A$-module direct summand of
$\gM_{A/I}$, and we conclude that the map $\gM_A^+ \cotimes_{A} A/I \to
\im \Phi_{\gM_{A/I}}$ induced by $f$ is an isomorphism.

Finally, we have the diagram
\[\xymatrix{\varphi^*\gM_A\ar[d]&\gM^+_A\ar[d]\\
    \varphi^*\gM_{A/I}\ar[r]& \im \Phi_{\gM_{A/I}}} \]
where the horizontal arrow is given by~$\Phi_{\gM_{A/I}}$, and the
right hand vertical arrow is~$f$. 
Since~$\varphi^*\gM_A$ is a
projective $\gS_{A}[\Gal(K'/K)]$-module, we may find a morphism of
$\gS_{A}[\Gal(K'/K)]$-modules $\varphi^*\gM_A\to \gM_A^+$ which fills in
the commutative square. Since the composite $\varphi^*\gM_A\to
\gM_A^+\to \im \Phi_{\gM_{A/I}} \cong \gM_A^+ \cotimes_{A} A/I $ is surjective, it follows by Nakayama's
lemma 
that
$\varphi^*\gM_A\to \gM_A^+$ is also surjective, and the composite
$\varphi^*\gM_A\to \gM_A^+\subset\gM_A$ gives a map
$\Phi_{\gM_A}$. Since $\Phi_{\gM_A}[1/E(u)]$ is a surjective
map of projective modules of the same rank, it is an isomorphism, and
we see that $\gM_A$ together with $\Phi_{\gM_A}$ is our required lifting to a
Breuil--Kisin module of rank~$d$ with descent data.

Since the source and target of $\Psi$ are of finite presentation
over $\Spf \cO$, and $\varpi$-adic, we see that $\Psi$ is representable
by algebraic spaces
(by \cite[Lem.~7.10]{Emertonformalstacks}) and locally of finite presentation
(by \cite[Cor.~2.1.8]{EGstacktheoreticimages}
and~\cite[\href{http://stacks.math.columbia.edu/tag/06CX}{Tag
  06CX}]{stacks-project}).  Thus $\Psi$ is in fact smooth (being formally
smooth and locally of finite presentation).
 \end{proof}

We now show that the inertial type of a Breuil--Kisin module is visible on the
local model.

\begin{lem}\label{lem: local model types map}There is a natural
  morphism  $\cM_{\loc}^{K'/K}\to \pi_0(\Rep_{I(K'/K)})$. 
\end{lem}
\begin{proof} 
  The morphism $\cM_{\loc}^{K'/K}\to \pi_0(\Rep_{I(K'/K)})$ is defined by sending
$(\gL,\gL^+)\mapsto\gL/\pi'$. More precisely, $\gL/\pi'$ is a rank~$d$
projective $k'\otimes_{\Zp}A$-module with a linear action of
$I(K'/K)$, so determines an $A/p$-point of~$\pi_0(\Rep_{I(K'/K)})=\prod_{i = 0}^{f'-1}
\pi_0(\Rep_d\bigl(I(K'/K)\bigr))$. Since the target is a disjoint
union of copies of~$\Spec\cO$, the morphism $\Spec A/p\to
\pi_0(\Rep_{I(K'/K)})$ lifts uniquely to a morphism $\Spec A \to
\pi_0(\Rep_{I(K'/K)})$, as required. 
\end{proof}

\begin{defn}\label{defn: local model stack fixed type}We let \[\cM_{\loc}^{(\tau_i)}:=\cM_{\loc}^{K'/K}\otimes_{\pi_0(\Rep_{I(K'/K)})}\pi_0(\Rep_{I(K'/K),\{\tau_i\}}).\]  If each $\tau_i = \tau$ for some fixed
  $\tau$, we write $\cM_{\loc}^\tau$ for $\cM_{\loc}^{(\tau_i)}$. By
  Lemma~\ref{lem:tame reps}, $\cM_{\loc}^{K'/K}$ is the disjoint union
  of the open and closed substacks~$\cM_{\loc}^{(\tau_i)}$.
  \end{defn}

\begin{lem} 
 We have~$\cC^{(\tau_i)}=\cC^{\dd}\times_{\cM_{\loc}^{K'/K}}\cM_{\loc}^{(\tau_i)}$. 
\end{lem}
\begin{proof}
  This is immediate from the definitions.
\end{proof}

In particular, $\cC^{(\tau_i)}$ is a closed substack of $\cC^{\dd}$.

\subsection{Local models:\ determinant conditions}\label{subsec: local
  models rank two}
Write~$N=K\cdot W(k')[1/p]$, so that~$K'/N$ is
totally ramified. Since~$I(K'/K)$ is cyclic of order prime to~$p$  and
acts trivially on~$\cO_N$, we may
write  \numequation\label{eqn: breaking up local model over
  xi}(\gL,\gL^+)=\oplus_{\xi}(\gL_{\xi},\gL^+_\xi)\end{equation}where the sum
is over all characters~$\xi:I(K'/K)\to\cO^\times$, and  
$\gL_\xi$ (resp.\ $\gL^+_\xi$) is the $\cO_N\otimes A$-submodule of~$\gL$
(resp.\ of~$\gL^+$) on which~$I(K'/K)$ acts through~$\xi$.
\begin{defn}\label{defn: strong determinant condition general local
    model} 
We say that an object $(\gL,\gL^+)$  of~$\cM_{\loc}^{K'/K}(A)$
\emph{satisfies the strong determinant condition} if Zariski
  locally on $\Spec A$ the
  following condition holds:\ for all $a\in\cO_N$ and all $\xi$, we
  have \numequation\label{eqn: strong det condn}
  \det{\!}_A(a|\gL^+_\xi)
  =\prod_{\psi:N\into E}\psi(a) 
  \end{equation}
  as polynomial functions on $\cO_N$ 
  in the sense of~\cite[\S5]{MR1124982}. 
  \end{defn}

  \begin{remark}\label{rem: explicit version of determinant condition
      local models version} 
An explicit version of this determinant condition is stated, in this generality, in~\cite[\S 2.2]{kis04}, specifically in the proof of \cite[Prop.\ 2.2.5]{kis04}. 
  We recall this here, with our notation. 
  We have a direct sum decomposition
  \[
  \cO_N\otimes_{\Z_p}A\toisom \bigoplus_{\sigma_i:k'\hookrightarrow \F}\cO_{N}\otimes_{W(k'),\sigma_i}A
  \]
  Recall that $e_{i}\in \cO_{N}\otimes_{\Z_p}\cO$ denotes the idempotent that identifies $e_i\cdot  \cO_N\otimes_{\Z_p}A$
  with the summand $\cO_{N}\otimes_{W(k'),\sigma_i}A$. For $j=0,1,\dots, e-1$, let $X_{j,\sigma_i}$ be an indeterminate. 
  Then the strong determinant condition on~$(\gL,\gL^+)$ is that for
  all~$\xi$, we have
  
\numequation\label{eqn: explicit strong det condn}  
\det{\!}_{A}\left(\sum_{j,\sigma_i}e_i\pi^jX_{j,\sigma_i}
\mid \gL^+_{\xi}\right) = 
\prod_{\psi}\sum_{j,\sigma_i}\left(\psi(e_i\pi^j)X_{j,\sigma_i}\right),
 \end{equation}
  where $j$ runs over $0,1,\dots,e-1$, $\sigma_i$ over embeddings $k'\hookrightarrow \F$, and $\psi$ over
  embeddings $\cO_N\hookrightarrow \cO$. Note that $\psi(e_i)=1$ if $\psi|_{W(k')}$ lifts $\sigma_i$
  and is equal to $0$ otherwise. 
  \end{remark}

\begin{defn}\label{defn: C^tau BT}
We write~$\cM_{\loc}^{K'/K,\BT}$ for the substack of~$\cM_{\loc}^{K'/K}$
given by those~$(\gL,\gL^+)$ which satisfy the strong determinant
condition. For each (possibly mixed) type~$(\tau_i)$, we write
$\cM_{\loc}^{(\tau_i),\BT}:=\cM_{\loc}^{(\tau_i)}\times_{\cM_{\loc}^{K'/K}}\cM_{\loc}^{K'/K,\BT}$.

Suppose for the remainder of this section that $d=2$ and $h=1$, so
that $\cC^{\dd}$ consists of Breuil--Kisin modules of rank two and height at
most $1$.
We then
  set~$\cC^{\dd,\BT}:=\cC^{\dd}\times_{\cM_{\loc}^{K'/K}}\cM_{\loc}^{K'/K,\BT}$,
  and
  $\cC^{(\tau_i),\BT}:=\cC^{(\tau_i)}\times_{\cM_{\loc}^{(\tau_i)}}\cM_{\loc}^{(\tau_i),\BT}$;\ 
  we also  write $\cC^{\dd,\BT,a} := \cC^{\dd,\BT} \times_{\cO}
  \cO/\varpi^a$ and $\cC^{(\tau_i),\BT,a} := \cC^{(\tau_i),\BT} \times_{\cO}
  \cO/\varpi^a$.

A Breuil--Kisin module $\gM \in
  \cC^{\dd}(A)$ is said to satisfy the strong determinant
  condition if and only if its image $\Psi(\gM) \in
  \cM_{\loc}^{K'/K}(A)$ does, i.e.\ if and only if it lies in $\cC^{\dd,\BT}$.

\end{defn}
\begin{prop}
  \label{prop:C tau BT is algebraic}$\cC^{(\tau_i),\BT}$ \emph{(}resp.\ $\cC^{\dd,\BT}$\emph{)} is a closed
  substack of $\cC^{(\tau_i)}$ \emph{(}resp.\ $\cC^{\dd}$\emph{)};\ in particular, it is a $\varpi$-adic formal
   algebraic
  stack of finite presentation over~$\cO$. 
\end{prop}
\begin{proof}
  This is immediate from Corollary~\ref{cor: C^dd,a is an algebraic
    stack} and the definition of the strong determinant
  condition as an equality of polynomial functions.
\end{proof}

  \begin{rem}\label{rem: strong determinant condition motivated by
      flat closure} 
The motivation for imposing the strong determinant condition is as
follows. One can take the flat part (in the sense of~\cite[Ex.\
9.11]{Emertonformalstacks}) of the $\varpi$-adic formal
stack~$\cC^{\dd}$, and on this flat part, one can impose the condition
that the corresponding Galois representations have all pairs of labelled
Hodge--Tate weights equal to~$\{0,1\}$;\ that is, we can consider the
substack of~$\cC^{\dd}$ corresponding to the Zariski closure of the
these Galois representations. 

We will soon see that $\cC^{\dd,\BT}$ is flat (Corollary~\ref{cor:
  Kisin moduli consequences of local models}).  By Lemma~\ref{lem: O points of moduli
  stacks} below, it follows that the substack of the previous paragraph is equal
to~$\cC^{\dd,\BT}$;\ so we may think of the strong determinant
condition as being precisely the condition which imposes this condition on the
labelled Hodge--Tate weights, and results in a formal stack which is
flat over~$\Spf\cO$. Since the inertial
types of $p$-adic Galois representations are unmixed, it is natural
from this perspective to expect
that $\cC^{\dd,\BT}$ should be the disjoint union of the stacks
$\cC^{\tau,\BT}$ for \emph{unmixed} types, and indeed this will be
proved shortly at Corollary~\ref{cor: BT C is union of tau C}.
  \end{rem}

To compare the strong determinant condition to the condition that the
type of a Breuil--Kisin module is unmixed, we make some observations about these conditions in the case of finite
field coefficients.

\begin{lem}
  \label{lem: determinant condition explicit finite field local models
  version}Let $\F'/\F$
  be a finite extension, and let~$(\gL,\gL^+)$ be an object
  of~$\cM_{\loc}^{K'/K}(\F')$. Then~$(\gL,\gL^+)$
satisfies the strong determinant condition if and only if the
  following property holds:\ for each~$i$ and for each $\xi: I(K'/K)
  \to \cO^\times$ we have $\dim_{\F'} (\gL^+_i)_{\xi} = e$. 
\end{lem}
\begin{proof} This is proved in a similar way to~\cite[Lemma
  2.5.1]{kis04},  using the explicit formulation of the strong
    determinant condition from Remark~\ref{rem: explicit version
      of determinant condition local models version}. 
    In the
    notation of that remark, we see that the strong determinant
    condition holds at~$\xi$ if and only if for each embedding
    $\sigma_i:k'\into\F$ we have
\numequation\label{eqn: explicit strong det condn slightly rewritten}  
\det{\!}_{A}\left(\sum_{j}\pi^jX_{j,\sigma_i}
\mid (\gL^+_i)_{\xi}\right) = 
\prod_{\psi}\sum_{j}\left(\psi(\pi)^jX_{j,\sigma_i}\right),
 \end{equation} where the product runs over the embeddings~$\psi:\cO_N\into\cO$
    with the property that~$\psi|_{W(k')}$ lifts~$\sigma_i$.
Since $\pi$ induces a nilpotent 
endomorphism of $(\gL^+_i)_{\xi}$ the left-hand side of \eqref{eqn: explicit strong det condn slightly rewritten} evaluates to $X_{0,\sigma_i}^{\dim_{\F'} (\gL^+_i)_{\xi}}$
while the right-hand side, which can be viewed as a norm
from $\cO_N \otimes_{\Zp} \F'$ down to $W(k') \otimes_{\Zp} \F'$, is equal to
$X_{0,\sigma_i}^e$.
     \end{proof}

\begin{lem} 
  \label{lem: determinant condition explicit finite field}Let $\F'/\F$
  be a finite extension, and let $\gM$ be a Breuil--Kisin module of rank~2 and
  height at most one with $\F'$-coefficients and descent data.
  
\begin{enumerate}
\item  $\gM$ satisfies the strong determinant condition if and only if the
  following property holds:\ for each~$i$ and for each $\xi
  : I(K'/K) \to \cO^\times$ we have $\dim_{\F'} (\im
  \Phi_{\gM,i}/E(u)\gM_i)_{\xi} = e$.
\item If~$\gM$ satisfies the strong determinant condition, then the determinant of~$\Phi_{\gM,i}$
  with respect to some choice of basis has $u$-adic
  valuation~$e'$. 
  \end{enumerate}
\end{lem}
\begin{proof} 
  The first part is immediate from Lemma~\ref{lem: determinant condition explicit finite field local models
  version}. 
  For the second part, let $\Phi_{\gM,i,\xi}$ be the restriction
  of $\Phi_{\gM, i}$ to $\varphi^*(\gM_{i-1})_{\xi}$. We think of
  $\gM_{i}$ and $\varphi^*(\gM_{i-1})$ as free $\F'[\![v]\!]$-modules of rank $2e(K'/K)$,
  where $v=u^{e(K'/K)}$. 
  We have
  \[
  \det{\!}_{\F'[\![v]\!]}(\Phi_{\gM,i}) =\left(\det{\!}_{\F'[\![u]\!]}(\Phi_{\gM,i})\right)^{e(K'/K)}.
  \] 
  Since $\Phi_{\gM,i}$ commutes with the descent datum, we also have 
  \[
  \det{\!}_{\F'[\![v]\!]}(\Phi_{\gM, i}) =\prod_{\xi} \det{\!}_{\F'[\![v]\!]}(\Phi_{\gM,i,\xi}),
  \]
  where $\xi$ runs over the $e(K'/K)$ characters $I(K'/K)\to \cO^\times$. 
  
  The proof of the 
  second part of~\cite[Lemma 2.5.1]{kis04} implies
  that, for each $\xi$, $\det_{\F'[\![v]\!]}(\Phi_{\gM,i,\xi})$ is $v^e=u^{e'}$ times a unit. 
  Indeed, each $\gM_{i,\xi}$ is a free $\F'[\![v]\!]$-module of rank $2$. It admits
  a basis $\{e_{1,\xi},e_{2,\xi}\}$ such that $\im\Phi_{\gM,i,\xi} = \langle v^{i}e_{1,\xi}, v^je_{2,\xi}\rangle$ 
  for some non-negative integers $i,j$. The strong determinant condition on 
  $(\im \Phi_{\gM,\i,\xi}/v^e\gM_{i,\xi})$
  implies that $i+j = 2e-e =e$, and this is precisely the $v$-adic valuation of $\det_{\F'[\![v]\!]}(\Phi_{\gM,i,\xi})$.
  We deduce that the $u$-adic valuation of  
  $\left(\det{\!}_{\F'[\![u]\!]}(\Phi_{\gM,i})\right)^{e(K'/K)}$ is $e(K'/K)\cdot e'$, which implies
  the second part of the lemma.  
   \end{proof}

By contrast, we have the following criterion for the type of a
Breuil--Kisin module to be unmixed.

 \begin{prop}\label{prop:types are unmixed}
  Let $\F'/\F$ be a finite extension, and let $\gM$ be a Breuil--Kisin module
  of rank~2 and height at most one with $\F'$-coefficients and descent
  data.  Then
the type of 
  $\gM$ is unmixed if and only if $\dim_{\F'} (\im
  \Phi_{\gM,i}/E(u)\gM_i)_{\xi}$ is independent of $\xi$ for each
  fixed $i$.  In particular, if
  $\gM$ satisfies the strong determinant condition, then the type of
  $\gM$ is unmixed. 
\end{prop}

\begin{proof} We begin the proof of the first part with the following observation. Let $\Lambda$ be a rank two free $\F'[[u]]$-module with an action of
 $I(K'/K)$ that is $\F'$-linear and $u$-semilinear with respect to a
 character $\chi$  (i.e., such that $g(u\lambda) = \chi(g)u
 g(\lambda)$ for $\lambda \in \Lambda$). In particular $I(K'/K)$ acts on $\Lambda/u\Lambda$ through a
 pair of characters which we call $\eta$ and $\eta'$. Let $\Lambda' \subset \Lambda$ be a rank two $I(K'/K)$-sublattice. We claim that there are integers $m, m'
 \ge 0$ such that the multiset of characters of $I(K'/K)$ occurring
in $\Lambda/\Lambda'$ has the shape $$\{ \eta \chi^i : 0 \le i < m\} \cup \{ \eta'
\chi^i : 0 \le j < m'\}$$ and the multiset of
characters occurring in $\Lambda'/u\Lambda'$ is $\{\eta \chi^m, \eta' \chi^{m'}\}$.

To check the claim we proceed by induction on $\dim_{\F'}
\Lambda/\Lambda'$, the case $\Lambda=\Lambda'$ being
trivial. Suppose $\dim_{\F'} \Lambda/\Lambda'=1$, so that $\Lambda'$ lies between $\Lambda$ and
$u\Lambda$. Consider the chain of containments $\Lambda \supset \Lambda' \supset u\Lambda
\supset u\Lambda'$. If without loss of generality $I(K'/K)$ acts on $\Lambda/\Lambda'$
via $\eta$, then it acts on $\Lambda'/u\Lambda$ by $\eta'$ and $u\Lambda/u\Lambda'$ by $\chi
\eta$, proving the claim with $m=1$ and $m'=0$. The general case follows by iterated application of the case
$\dim_k \Lambda/\Lambda'=1$, noting that since $I(K'/K)$ is abelian
the quotient $\Lambda/\Lambda'$ has a filtration by $I(K/K')$-submodules
whose graded pieces have dimension $1$.

Now return to the statement of the proposition. Let $(\tau_i)$ be the mixed type of $\gM$ and write $\tau_{i-1} = \eta
\oplus \eta'$. We apply the preceding observation with $\Lambda = \im
  \Phi_{\gM,i}$ and $\Lambda' = E(u)\gM_i = u^{e'} \gM_i$. Note that
  $\chi$ is a generator of the cyclic group $I(K'/K)$ of order $e'/e$. 
  Since
  $\Phi_{\gM}$ commutes with descent data, the group $I(K'/K)$ acts
  on $\Lambda/u\Lambda$ via $\eta$ and $\eta'$.   Then the
  the multiset 
\[\{ \eta \chi^i : 0 \le i < m\} \cup \{ \eta'
\chi^i : 0 \le j < m'\}\]
contains each character of $I(K'/K)$ with equal multiplicity if and
only if  one of $\eta,\eta'$  is the successor to $\eta\chi^{m-1}$ in the list
$\eta,\eta\chi,\eta\chi^2,\ldots$, and the other is the successor to
$\eta'\chi^{m'-1}$ in the list $\eta',\eta'\chi,\eta'\chi^2,\ldots$,
i.e., if and only if  $\{\eta\chi^m,\eta'\chi^{m'}\} = \{ \eta,\eta'\}$. Since $\gM_i/u\gM_i
\cong u^{e'}\gM_i/u^{e'+1}\gM_i = \Lambda'/u\Lambda'$, this occurs if
and only if
that $\tau_i = \tau_{i-1}$.

Finally, the last part of the proposition follows immediately from the
  first part and Lemma~\ref{lem: determinant condition explicit finite
    field}.\end{proof}

\begin{cor} 
  \label{cor: BT C is union of tau C}$\cC^{\dd,\BT}$ is the disjoint
  union of its closed substacks $\cC^{\tau,\BT}$. 
\end{cor}
\begin{proof}
  This follows from Propositions~\ref{prop: Ctau is
    an algebraic stack} and~\ref{prop:types are unmixed}. Indeed, from Proposition~\ref{prop: Ctau is
    an algebraic stack}, it suffices to show that if~$(\tau_i)$ is a
  mixed type,
  and~$\cC^{(\tau_i),\BT}$ is
  nonzero, then $(\tau_i)$ is in fact an unmixed type. Indeed,
  note that if $\cC^{(\tau_i),\BT}$ is nonzero, then it contains a dense set of
  finite type points, so in particular contains an $\F'$-point for
  some finite extension $\F'/\F$. It follows from
  Proposition~\ref{prop:types are unmixed} that the type is unmixed,
  as required.
\end{proof}

  \begin{rem}
    \label{rem: we don't consider mixed type local models} Since our
    primary interest is in Breuil--Kisin modules, we will have no further need
    to consider the stacks $\cM^{(\tau_i),\BT}_{\loc}$ or
    $\cC^{(\tau_i),\BT}$ for types that are not unmixed.
  \end{rem}

Let~$\tau$ be a tame
type;\ since $I(K'/K)$ is cyclic, we can write~$\tau=\eta\oplus\eta'$
for (possibly equal) characters $\eta,\eta':I(K'/K)\to \cO^\times$. 
Let~$(\gL,\gL^+)$ be an object of $\cM_{\loc}^{\tau}(A)$. Suppose that~$\xi\ne\eta,\eta'$. Then elements of $\gL_{\xi}$ are
divisible by $\pi'$ in $\gL$, and so multiplication by~$\pi'$
induces an isomorphism 
 of
projective $e_i(\cO_N \otimes A)$-modules of equal rank  \[p_{i,\xi} :
  e_i \gL_{\xi\chi_i^{-1} }
\isoto e_i\gL_{\xi} \]
where~$\chi_i:I(K'/K)\to\cO^{\times}$ is the
character sending $g\mapsto \sigma_i(h(g))$.  The 
induced map  \[p_{i,\xi}^+ : e_i \gL^+_{\xi\chi_i^{-1} }
\longrightarrow  e_i\gL^+_{\xi} \] is in particular an injection.
 The following lemma will be useful for checking the strong
 determinant condition when comparing various different stacks of local model
 stacks. 

\begin{lem} \label{lem: local model of type only have to check the strong
      determinant condition on that type}Let~$(\gL,\gL^+)$ be an object of $\cM_{\loc}^{\tau}(A)$. Then
  $(\gL,\gL^+)$ is an object of $\cM_{\loc}^{\tau,\BT}(A)$ if and only
  if both
  \begin{enumerate}
\item   the condition~\eqref{eqn: strong det condn} holds 
  for~$\xi=\eta,\eta'$, and 
  \item the injections $p_{i,\xi}^+ : e_i \gL^+_{\xi\chi_i^{-1} }
\longnto{\pi'}  e_i\gL^+_{\xi} $ are isomorphisms for all $\xi \neq
\eta,\eta'$ and for all $i$.
  \end{enumerate}
The second condition is equivalent to
\begin{enumerate}
\item[($\textrm{2}\,'$)] we have $(\gL^+/\pi' \gL^+)_{\xi} = 0$ for all $\xi \neq \eta,\eta'$.
\end{enumerate}
   \end{lem}
     \begin{proof}

The equivalence between (2) and ($2'$) is straightforward.
Suppose now that $\xi \neq \eta,\eta'$. Locally on $\Spec A$ the module $e_i \gL^{+}_{\xi \chi_i^{-1}}$ is by
definition a
  direct summand of $e_i \gL_{\xi \chi_i^{-1}}$.  Since $p_{i,\xi}$ is an isomorphism, the image of $p^+_{i,\xi}$ is locally on $\Spec A$ a direct summand of
$e_i\gL^+_{\xi}$. 
Under the assumption that~(\ref{eqn: explicit strong det condn slightly
  rewritten}) holds for $i$ and~ $\xi \chi_i^{-1}$, the condition~(\ref{eqn: explicit strong det condn slightly
  rewritten}) for $i$ and $\xi$ is therefore  equivalent to the surjectivity of
$p^+_{i,\xi}$.  The lemma follows upon  noting that $\chi_i$ is a generator
of the
group of characters $I(K'/K)\to \cO^{\times}$.  
\end{proof}

To conclude this section we describe the $\cO_{E'}$-points of $\cC^{\dd,\BT}$, for $E'/E$
a finite extension;\ recall that our convention is that a two-dimensional Galois
representation is Barsotti--Tate if all its labelled pairs of
Hodge--Tate weights are equal to~$\{0,1\}$ (and not just that all of
the labelled Hodge--Tate weights are equal to~$0$ or~$1$).
\begin{lem} 
\label{lem: O points of moduli stacks}Let $E'/E$ be a finite
extension. Then the $\Spf(\cO_{E'})$-points of $\cC^{\dd,\BT}$ correspond
precisely to the potentially Barsotti--Tate Galois representations
$G_K\to\GL_2(\cO_{E'})$ which become Barsotti--Tate over~$K'$;\ and the
$\Spf(\cO_{E'})$-points of $\cC^{\tau,\BT}$ correspond to those
representations which are potentially Barsotti--Tate of type~$\tau$.
\end{lem}
\begin{proof} 
In light of Lemma~\ref{lem:points-of-Cdd} and the first sentence of
Remark~\ref{rem: Kisin modules don't always have a type}, we are
reduced to checking that a Breuil--Kisin module of rank $2$ and height $1$ with $\cO_{E'}$-coefficients
and descent data corresponds to a potentially Barsotti--Tate
representation if and only if it satisfies the strong determinant
condition, as well as checking that the descent data on the Breuil--Kisin module
matches the type of the corresponding Galois representation. 

Let $\gM_{\cO_{E'}} \in\cC^{\dd,\BT}(\Spf(\cO_{E'}))$. Plainly
$\gM_{\cO_{E'}}$ satisfies the strong determinant condition if and
only if $\gM := \gM_{\cO_{E'}}[1/p]$ satisfies the strong determinant
condition (with the latter having the obvious meaning). 
Consider the filtration
\[
\mathrm{Fil}^1(\varphi^*(\gM_i)):=\{m\in \varphi^*(\gM_i)\mid \Phi_{\gM,i+1}(m)\in 
E(u)\gM_{i+1}\}\subset \varphi^*\gM_{i}
\]
inducing
\[
\mathrm{Fil}^1_i\subset \varphi^*(\gM_{i})/ E(u)\varphi^*(\gM_{i}).
\]
Note that $\varphi^*(\gM_{i})/E(u)\varphi^*(\gM_{i})$ 
is isomorphic to a free $K'\otimes_{W(k'),\sigma_i}E'$-module 
of rank $2$. Then $\gM$ corresponds to a Barsotti--Tate Galois representation
\[
G_{K'}\to \GL_{2}(E')
\]
if and only, if for every $i$, $\mathrm{Fil}^1_i$ is isomorphic to $K'\otimes_{W(k'),\sigma_i}E'$ as a 
$K'\otimes_{W(k'),\sigma_i}E'$-submodule of $\varphi^*\gM_{i}/E(u)\varphi^*\gM_{i}$. 
This follows, for example, by specialising the proof of~\cite[Cor. 2.6.2]{kisindefrings}
to the Barsotti--Tate case (taking care to note that the conventions
of \emph{loc.\ cit.} for Hodge--Tate weights and for Galois
representations associated to Breuil--Kisin modules are both dual to ours). 

Let $\xi: I(K'/K)\to\cO^\times$ be a character. Consider the filtration
\[
\mathrm{Fil}^1_{i,\xi}\subset \varphi^*(\gM_i)_{\xi}/E(u)\varphi^*(\gM_i)_{\xi}\simeq N^2\otimes_{W(k'),\sigma_i}E'
\]
induced by $\mathrm{Fil}^1_i$. The strong determinant condition on 
$(\im \Phi_{\gM,i+1}/E(u)\gM_{i+1})_{\xi}$
holds if and only if $\mathrm{Fil}^1_{i,\xi}$ 
is isomorphic to $N\otimes_{W(k'),\sigma_i}E'$. By~\cite[Lemma 5.10]{2015arXiv151007503C},
we have an isomorphism of $K'\otimes_{W(k'),\sigma_i}E'$-modules 
\[
\mathrm{Fil}^1_{i}\simeq K' \otimes_{N} \mathrm{Fil}^1_{i,\xi}. 
\] 
This, together with the previous paragraph, allows us to conclude. 
Note that, since $u$ acts invertibly when working with $E'$-coefficients
and after quotienting by $E(u)$, the argument is independent of the
choice of character $\xi$.  

For the statement about types, let $S_{K'_0}$ be Breuil's period ring
(see e.g.\ \cite[\S5.1]{MR1804530}) endowed with the evident action of
$\Gal(K'/K)$ compatible with the embedding $\gS \into S_{K'_0}$.
Here $K'_0$ is the maximal unramified extension in $K'$.  Recall that by \cite[Cor.\
3.2.3]{MR2388556} there is a canonical  $(\varphi,N)$-module isomorphism
\numequation\label{eq:compare-M-and-D}  S_{K'_0} \otimes_{\varphi,\gS} \gM \cong S_{K'_0} \otimes_{K'_0}
  D_{\mathrm{pcris}}(T(\gM)). \end{equation}
One sees from its construction 
that the 
isomorphism \eqref{eq:compare-M-and-D} is in fact equivariant for the
action of $I(K'/K)$, 
and
the claim follows by reducing modulo $u$ and all its divided
powers.
\end{proof}

\subsection{Change of extensions}\label{subsec: change of extension K
  prime over K}
We now discuss the behaviour of various of our constructions under
change of~$K'$. Let~$L'/K'$ be an extension of degree prime to~$p$ such that~$L'/K$ is
Galois. We suppose that we have fixed a uniformiser~$\pi''$ of~$L'$
such that
$(\pi'')^{e(L'/K')}=\pi'$. Let~$\gS'_A:=(W(l')\otimes_{\Zp}A)[[u]]$,
where~$l'$ is the residue field of~$L'$, and let~$\Gal(L'/K)$ and~$\varphi$ act
on~$\gS'_A$ via the prescription of Section~\ref{subsec: kisin modules
  with dd} (with~$\pi''$ in place of~$\pi'$). 

There is a natural injective ring
homomorphism~$\cO_{K'}\otimes_{\Zp}A\to\cO_{L'}\otimes_{\Zp}A$, which
is equivariant for the action of~$\Gal(L'/K)$ (acting on the source via the
natural surjection~$\Gal(L'/K)\to\Gal(K'/K)$). There is also an
obvious injective ring homomorphism $\gS_A\to\gS'_A$
sending~$u\mapsto u^{e(L'/K')}$, which is equivariant for the actions
of~$\varphi$ and~$\Gal(L'/K)$;\ we have $(\gS'_{A})^{\Gal(L'/K')}=\gS_A$. 
If~$\tau$ is an inertial type for~$I(K'/K)$,
we write~$\tau'$ for the corresponding type for~$I(L'/K)$, obtained by
inflation.

For any $(\gL,\gL^+)\in\cM_{\loc}^{K'/K}$, we define 
$(\gL',(\gL')^+)\in\cM_{\loc}^{L'/K}$
by \[(\gL',(\gL')^+):=\cO_{L'}\otimes_{\cO_{K'}}(\gL,\gL^+),\]with the diagonal
action of~$\Gal(L'/K)$.  Similarly, for
any~$\gM\in\cC^{\dd}(A)$, we let~$\gM':=\gS'_A\otimes_{\gS_A}\gM$,
with~$\varphi$ and~$\Gal(L'/K)$ again acting diagonally. 

\begin{prop} 
  \label{prop: change of K'}\hfill
  \begin{enumerate}
  \item The assignments~$(\gL,\gL^+)\to(\gL',(\gL')^+)$ and 
  $\gM\mapsto\gM'$ induce compatible monomorphisms 
  $\cM_{\loc}^{K'/K}\to\cM_{\loc}^{L'/K}$
  and~$\cC^{\dd}_{K'}\to\cC^{\dd}_{L'}$, i.e., as functors they are 
  fully faithful. 

\item The monomorphism $\cC^{\dd}_{K'}\to\cC^{\dd}_{L'}$ induces an
  isomorphism~$\cC^{\tau}\to\cC^{\tau'}$, as well as a monomorphism
$\cC^{\dd,\BT}_{K'}\to\cC^{\dd,\BT}_{L'}$ 
and an isomorphism $\cC^{\tau,\BT}\to\cC^{\tau',\BT}$. 
  \end{enumerate}
\end{prop}
\begin{proof} 
(1) One checks easily that the  assignments~$(\gL,\gL^+)\to(\gL',(\gL')^+)$ and 
  $\gM\mapsto\gM'$ are compatible. For the rest of the claim, we consider the case of the functor $\gM\mapsto\gM'$;\ the
  arguments for the local models case are similar but slightly easier,
  and we leave them to the reader. Let $A$ be a $\varpi$-adically
  complete $\cO$-algebra. 
If~$\gN$ is a
rank~$d$ Breuil--Kisin module with descent data from~$L'$ to~$K$, consider
the Galois
invariants~$\gN^{\Gal(L'/K')}$. Since~$(\gS'_{A})^{\Gal(L'/K')}=\gS_A$,
these invariants are naturally an~$\gS_A$-module, and moreover they
naturally carry a Frobenius and descent data satisfying the conditions
required of a Breuil--Kisin module of height at most $h$. In general the
invariants need not be projective of rank $d$, and so need not be a
rank~$d$ Breuil--Kisin module with descent data from~$K'$ to~$K$. However, in
the case~$\gN=\gM'$ we
have \[(\gM')^{\Gal(L'/K')}=(\gS'_A\otimes_{\gS_A}\gM)^{\Gal(L'/K')}=(\gS'_A)^{\Gal(L'/K')}\otimes_{\gS_A}\gM=\gM.\]
Here the second equality holds e.g.\ because $\Gal(L'/K')$ has order
prime to $p$, so that taking $\Gal(L'/K')$-invariants is exact, and
indeed is given by
multiplication by an idempotent $i \in \gS'_A$ (use the decomposition
$\gS'_A = i \gS'_A \oplus (1-i) \gS'_A$ and note that $i$ kills the
latter summand).
It
follows immediately that the functor~$\gM\to\gM'$ is fully faithful,
so $\cC^{\dd,a}_{K'}\to\cC^{\dd,a}_{L'}$ is a monomorphism.

(2) Suppose now that~$\gN$ has type~$\tau'$. In view of what we have
proven so far, in order to prove that
$\cC^{\tau}\to\cC^{\tau'}$ is an isomorphism, it is enough to
show that $\gN^{\Gal(L'/K')}$ is a rank~$d$ Breuil--Kisin module of type~$\tau$, and that the natural
map of $\gS'_A$-modules
\numequation\label{eqn: tensoring up invariants}\gS'_{A}\otimes_{\gS_A}\gN^{\Gal(L'/K')}\to\gN\end{equation}
is an
isomorphism. For the remainder of this proof, for clarity we write $u_{K'}$,
$u_{L'}$ instead of $u$ for the variables of $\gS_A$ and $\gS'_A$ respectively.
Since the type~$\tau'$ of $\gN$ is inflated from $\tau$, the
action of $\Gal(L'/K')$ on $\gN/u_{L'} \gN$ factors through
$\Gal(l'/k')$;\ noting that $W(l')$ has a normal basis for $\Gal(l'/k')$ over
$W(k')$,  we obtain an isomorphism
\numequation\label{eq:moduL'} 
 W(l')  \otimes_{W(k')}  (\gN/u_{L'}\gN)^{\Gal(L'/K')} \toisom
  \gN/u_{L'} \gN.\end{equation}
In particular the $W(k')\otimes_{\Zp} A$-module $(\gN/u_{L'} \gN)^{\Gal(L'/K')}$ is projective of rank $d$.

Observe however that $(\gN/u_{K'}\gN)^{\Gal(L'/K')} =
(\gN/u_{L'}\gN)^{\Gal(L'/K')}$. To see this, by the exactness of
taking $\Gal(L'/K')$ invariants it suffices to check that $u_{L'}^i
\gN/u_{L'}^{i+1}$ has trivial $\Gal(L'/K')$-invariants  for $0 < i <
e(L'/K')$. Multiplication by $u_{L'}^i$ gives an isomorphism 
$\gN/u_{L'} \gN \cong u_{L'}^i
\gN/u_{L'}^{i+1}$, so that for $i$ in the above range, the inertia group 
$I(L'/K')$ acts linearly on $ u_{L'}^i \gN/u_{L'}^{i+1}$ through a twist of
$\tau'$ by a nontrivial character;\ so there are no
$I(L'/K')$-invariants, and thus no $\Gal(L'/K')$-invariants either.

It
follows that the isomorphism \eqref{eq:moduL'} is the map \eqref{eqn: tensoring up
  invariants} modulo~$u_{L'}$. By 
Nakayama's lemma it follows that \eqref{eqn: tensoring up
  invariants} is surjective. Since~$\gN$ is projective,  the surjection~\eqref{eqn: tensoring up
  invariants} is split, and is therefore an isomorphism, since it is
an isomorphism modulo~$u_{L'}$. This isomorphism
exhibits~$\gN^{\Gal(L'/K')}$ as a direct summand (as
an~$\gS_A$-module) of the projective module~$\gN$, so it is also
projective;\ and it is projective of rank~$d$, since this holds modulo~$u_{K'}$. 

Finally, we need to check the compatibility of these maps with the
strong determinant condition.  By Corollary~\ref{cor: BT C is union of
  tau C}, it is enough to prove this for the case of the morphism
$\cC^{\tau}\to\cC^{\tau'}$ for some~$\tau$;\ by the compatibility in
part (1), this comes down to the same for the corresponding map of
local model stacks $\cM_{\loc}^{\tau}\to\cM_{\loc}^{\tau'}$. If
$(\gL,\gL^+) \in \cM_{\loc}^{\tau}$, it therefore suffices to show 
to show that conditions (1) and ($2'$) of  Lemma~\ref{lem: local model of type only have to check the strong
      determinant condition on that type} for $(\gL,\gL^+)$ and
    $(\gL', (\gL')^+) := \cO_{L'} \otimes_{\cO_{K'}} (\gL,\gL^+)$ are
    equivalent. This is immediate for condition ($2'$), since we have
    $(\gL')^+/\pi'' (\gL')^+ \cong l' \otimes_{k'} (\gL^+/\pi' \gL^+)$
 as $I(L'/K)$-representations.

Writing~$\tau=\eta\oplus\eta'$, it remains to relate the
    strong determinant conditions on the~$\eta,\eta'$-parts over
    both~$K'$ and~$L'$. Unwinding the definitions using Remark~\ref{rem: explicit version of determinant condition
      local models version}, one finds that the condition over $L'$ is
    a product of $[l':k']$ copies (with different sets of variables)
    of the condition over~$K'$. Thus the strong determinant condition
    over $K'$ implies the condition over $L'$, while the condition
    over $L'$ implies the condition over $K'$ up to an $[l':k']$th
    root of unity. Comparing the terms involving only copies of 
    $X_{0,\sigma_i}$'s shows that this root of unity must be $1$.
\end{proof}

\begin{remark}
 The morphism of local model stacks $\cM^{\tau}_{\loc} \to
 \cM^{\tau'}_{\loc}$ is not an isomorphism (provided that the extension
 $L'/K'$ is nontrivial). The issue is that, as we observed in the
 preceding proof, local models
 $(\gL',(\gL')^+)$ in the image
 of the morphism $\cM^{\tau}_{\loc} \to
 \cM^{\tau'}_{\loc}$ can have  $((\gL')^+/\pi'' (\gL')^+)_{\xi} \neq
 0$ only for characters $\xi : I(L'/K) \to \cO^{\times}$ that are
 inflated from $I(K'/K)$. However, one does obtain an isomorphism
 from the substack of $\cM^{\tau}_{\loc}$ of pairs $(\gL,\gL^+)$ satisfying condition~(2)
 of Lemma~\ref{lem: local model of type only have to check the strong
      determinant condition on that type} to the analogous substack of
    $\cM^{\tau'}_{\loc}$;\ therefore the induced map
    $\cM^{\tau,\BT}\rightarrow \cM^{\tau',\BT}$ \emph{will} also be an
    isomorphism.  Analogous remarks will apply to the maps 
    of local model stacks in \S\ref{subsec: explicit local models}.
\end{remark}

\subsection{Explicit local models}\label{subsec: explicit local models}
We now explain the connection between the moduli stacks $\cC^\tau$
and local models for Shimura varieties at hyperspecial and Iwahori level. 
This idea has been developed in~\cite{2015arXiv151007503C} for
Breuil--Kisin modules of arbitrary rank with tame principal series descent data, inspired by~\cite{kis04}, which analyses the case without descent data. 

The results
of~\cite{2015arXiv151007503C} relate the moduli stacks $\cC^{\tau}$(in
the case that~$\tau$ is a principal series type) via
a local model diagram to a mixed-characteristic deformation of the affine flag variety,
introduced in this generality by Pappas and Zhu~\cite{pappas-zhu}. The local models in~\cite[\S 6]{pappas-zhu}
are defined in terms of Schubert cells in the generic fibre of this mixed-characteristic deformation,
by taking the Zariski closure of these Schubert cells. The disadvantage of this
approach is that it does not give a direct moduli-theoretic interpretation of the special
fibre of the local model. Therefore, it is hard to check directly that our stack $\cC^{\tau, \BT}$, 
which has a moduli-theoretic definition,
corresponds to the desired local model under the diagram introduced in~\cite[Cor. 4.11]{2015arXiv151007503C} 
\footnote{One should be able to check this by adapting the ideas
in~\cite[\S 2.1]{haines-ngo} and~\cite[Prop. 6.2]{pappas-zhu} to $\Res_{K/\Q_p}\GL_n$
where $K/\Q_p$ can be ramified.}. 

In our rank~2 setting, the local models admit a much more explicit condition, using
lattice chains and Kottwitz's determinant condition, and in the cases
of nonscalar types, we will relate our local models to 
the naive local model at Iwahori level
for the Weil restriction of $\GL_2$, in the sense of~\cite[\S
2.4]{MR3135437}.

We begin with the simpler case of scalar inertial types. Suppose  that~$K'/K$ is totally ramified, and that~$\tau$ is a scalar
inertial type, say~$\tau=\eta\oplus\eta$. In this case we define the
local model stack~$\cM_{\loc,\hyp}$ (``hyp'' for ``hyperspecial'') to be the \emph{fppf} stack
over~$\Spf\cO$ (in fact, a $p$-adic formal algebraic stack),
which to each $p$-adically complete
$\cO$-algebra $A$ associates the groupoid~$\cM_{\loc,\hyp}(A)$
consisting of pairs~$(\gL,\gL^+)$,
where
\begin{itemize}
\item $\gL$ is a rank~$2$ projective~$\cO_K\otimes_{\Zp}A$-module, and
\item $\gL^+$ is an $\cO_K\otimes_{\Zp}A$-submodule of $\gL$,
	which is locally on~$\Spec A$ a direct summand of~$\gL$ as an $A$-module (or equivalently, for which the quotient $\gL/\gL^+$ is projective
	as an $A$-module).
\end{itemize}
We let~$\cM^{\BT}_{\loc,\hyp}$ be the substack of
pairs~$(\gL,\gL^+)$ with the property that for all $a\in\cO_K$,
we have \numequation\label{eqn: strong det condn for scalar local model} \det{\!}_A(a|\gL^+)
=\prod_{\psi:K\into E}\psi(a)
  \end{equation}
  as polynomial functions on $\cO_K$.

\begin{lem}
  \label{lem: local model scalar type is equivalent to what you
    expect}
The functor
  $(\gL',(\gL')^+) \mapsto ((\gL')_{\eta},(\gL')^+_{\eta})$
  defines a morphism
$\cM^\tau_{\loc}\to\cM_{\loc,\hyp}$ which induces an isomorphism
  $\cM^{\tau,\BT}_{\loc} \to \cM^{\BT}_{\loc,\hyp}$. 
{\em (}We remind the reader that $K'/K$ is assumed totally ramified,
and that $\tau$ is assumed to be a scalar inertial type associated
to the character $\eta$.{\em )}
\end{lem}
\begin{proof}
If $(\gL',(\gL')^+)$ is an object of $\cM^{\tau}_{\loc}$,
the proof that $((\gL')_\eta,(\gL')^+_\eta)$
  is indeed an object of~$\cM_{\loc,\hyp}(A)$ is
  very similar to the proof of Proposition~\ref{prop: change of K'},
  and is left to the reader.
  Similarly, the reader may verify that the functor
  \[(\gL,\gL^+) \mapsto
(\gL',(\gL')^+):=\cO_{K'}\otimes_{\cO_K}(\gL,\gL^+),\]
where the
action of~$\Gal(K'/K)$ is given by the tensor product of the natural
action on~$\cO_{K'}$ with the action on~$(\gL,\gL^+)$ given by the character~$\eta$,
defines a morphism
$\cM_{\loc,\hyp} \to \cM^\tau_{\loc}$.

 The composition
$\cM_{\loc,\hyp} \to \cM^\tau_{\loc}\to\cM_{\loc,\hyp}$ is evidently
the identity. The composition in the other order is not, in general, naturally equivalent to the
identity morphism, because for $(\gL,\gL^+) \in
\cM^{\tau}_{\loc}(A)$ one cannot necessarily recover $\gL^+$ from the
projection to its $\eta$-isotypic part. However, this 
\emph{will} hold if $\gL^+$ satisfies condition (2) of  Lemma~\ref{lem:
      local model of type only have to check the strong determinant
      condition on that type} (and so in particular will hold after
    imposing the strong determinant condition). 

Indeed, suppose $(\gL,\gL^+) \in
\cM^{\tau}_{\loc}(A)$. Then there is a natural
$\Gal(K'/K)$-equivariant map of projective $\cO_{K'}\otimes_{\Zp} A$-modules
\numequation\label{eq:recover-L} \cO_{K'} \otimes_{\cO_K} \gL_{\eta} \to \gL\end{equation}
of the same rank (in which $\Gal(K'/K)$ acts by $\eta$ on $\gL_{\eta}$). This map is
surjective because it is surjective on $\eta$-parts and the maps $p_{i,\xi}$ are surjective for all
$\xi\neq \eta$;\ therefore it is an isomorphism. One further has a 
a natural
$\Gal(K'/K)$-equivariant map of $\cO_{K'}\otimes_{\Zp} A$-modules
\numequation\label{eq:recover-L+} \cO_{K'} \otimes_{\cO_K}
\gL^+_{\eta} \to \gL^+ \end{equation}
that is injective because locally on $\Spec(A)$ it is a direct summand of the isomorphism \eqref{eq:recover-L}.
If one further assumes that $\gL^+$ satisfies condition (2) of  Lemma~\ref{lem:
      local model of type only have to check the strong determinant
      condition on that type} then \eqref{eq:recover-L+} is an
    isomorphism, as claimed.

  It remains to check the compatibility of these maps with the strong
  determinant condition. If $(\gL,\gL^+) \in \cM_{\loc,\hyp}$ then 
certainly condition (2) of  Lemma~\ref{lem:
      local model of type only have to check the strong determinant
      condition on that type} holds for
    $(\gL',(\gL')^+):=\cO_{K'}\otimes_{\cO_K}(\gL,\gL^+) \in
    \cM^\tau_{\loc}$.  By Lemma~\ref{lem:
      local model of type only have to check the strong determinant
      condition on that type} the strong determinant condition holds
    for~$(\gL',(\gL')^+)$ if and only if~(\ref{eqn: strong det condn}) holds
    for~$\gL'$ with~$\xi=\eta$;\ but this is exactly the
    condition~(\ref{eqn: strong det condn for scalar local model}) for
    $\gL$, as required.
\end{proof}

Next, we consider the case of principal series types. We suppose
that~$K'/K$ is totally ramified. We begin by
defining a $p$-adic formal algebraic stack $\cM_{\loc,\Iw}$
over~$\Spf\cO$ (``Iw'' for Iwahori). For each complete $\cO$-algebra $A$, we
let~$\cM_{\loc,\Iw}(A)$ be the groupoid of tuples~$(\gL_1,\gL^+_1,\gL_2,\gL^+_2,f_1,f_2)$,
where
\begin{itemize}
\item $\gL_1$, $\gL_2$ are rank~$2$
  projective~$\cO_K\otimes_{\Zp}A$-modules, 
\item $f_1:\gL_1\to\gL_2$, $f_2:\gL_2\to\gL_1$ are morphisms of
  $\cO_K\otimes_{\Zp}A$-modules, satisfying~$f_1\circ f_2=f_2\circ f_1=\pi$,
\item both $\coker f_1$ and~$\coker f_2$ are rank one projective $k\otimes_{\Zp}A$-modules,
\item $\gL^+_1$, $\gL^+_2$ are $\cO_K\otimes_{\Zp}A$-submodules of $\gL_1$, $\gL_2$, which are locally
  on~$\Spec A$ direct summands  as $A$-modules
  (or equivalently, for which the quotients
  $\gL_i/\gL^+_i$ ($i = 1,2$) are projective $A$-modules), and moreover
  for which
  the morphisms~$f_1,f_2$ restrict to morphisms~$f_1:\gL^+_1\to \gL^+_2$, $f_2:\gL^+_2\to \gL^+_1$.
\end{itemize}
We let~$\cM^{\BT}_{\loc,\Iw}$ be the substack of
tuples with the property that for all $a\in\cO_K$ and $i=1,2$,
we have \numequation\label{eqn: strong det condn for principal series local model} \det{\!}_A(a|\gL^+_i)
=\prod_{\psi:K\into E}\psi(a)
  \end{equation}
  as polynomial functions on $\cO_K$.

Write~$\tau=\eta\oplus\eta'$ with~$\eta\ne\eta'$. 
Recall that the character $h:\Gal(K'/K)=I(K'/K)\to W(k)^\times$ is
given by~$h(g)=g(\pi')/\pi'$. 
Since we are assuming
that~$\eta\ne\eta'$, for each embedding~$\sigma:k\into\F$ (which we
also think of as $\sigma:W(k)\into\cO$) there are
integers $0<a_\sigma,b_\sigma<e(K'/K)$ with the properties that
$\eta'/\eta=\sigma\circ h^{a_\sigma}$, $\eta/\eta'=\sigma\circ
h^{b_\sigma}$;\ in particular, $a_\sigma+b_\sigma=e(K'/K)$. Recalling
that~$e_\sigma\in W(k)\otimes_{\Zp}\cO\subset \cO_K$ is the idempotent
corresponding to~$\sigma$, we
set~$\pi_1=\sum_\sigma(\pi')^{a_\sigma}e_\sigma$,
$\pi_2=\sum_\sigma(\pi')^{b_\sigma}e_\sigma$;\ so we have
~$\pi_1\pi_2=\pi$, and~$\pi_1\in(\cO_{K'}\otimes_{\Zp}\cO)^{I(K'/K)=\eta'/\eta}$,
$\pi_2\in(\cO_{K'}\otimes_{\Zp}\cO)^{I(K'/K)=\eta/\eta'}$.

We define a
morphism $\cM^\tau_{\loc}\to\cM_{\loc,\Iw}$ as follows. Given a
pair~$(\gL,\gL^+)\in\cM^\tau_{\loc}(A)$, we
set~$(\gL_1,\gL^+_1)=(\gL_\eta,\gL^+_\eta)$,
$(\gL_2,\gL^+_2)=(\gL_{\eta'},\gL^+_{\eta'})$, and we let~$f_1$, $f_2$ be given by
multiplication by~$\pi_1$, $\pi_2$ respectively. The only point that
is perhaps not immediately obvious is to check that~$\coker f_1$ and~$\coker f_2$
are  rank one projective $k\otimes_{\Zp}A$-modules. To see this, note
that by~\cite[\href{http://stacks.math.columbia.edu/tag/05BU}{Tag
  05BU}]{stacks-project}, we can lift $(\gL/\pi'\gL)_{\eta'}$ to a
rank one summand~$U$ of the
projective~$\cO_{K'}\otimes_{\Zp}A$-module~$\gL$. Let~$U_\eta$
be the direct summand of~$\gL_\eta$ obtained by projection of~$U$ to
the $\eta$-eigenspace~$\gL_\eta$;\ then the projective $\cO_K
\otimes_{\Zp} A$-module $U_\eta$ has rank one, as can be checked
modulo $\pi$.
Similarly we
may lift $(\gL/\pi'\gL)_{\eta'}$ to a rank one
summand~$V$ of~$\gL$, and we let~$V_{\eta'}$ be the 
projection to the $\eta'$-part.

The natural map 
\numequation\label{eqn: decomposing principal series type over powers
  of pi
  prime}
\cO_{K'} \otimes_{\cO_K} (U_\eta \oplus V_{\eta'}) \to \gL 
\end{equation}
is an isomorphism,  since both sides are projective
$\cO_{K'}\otimes_{\Zp}A$-modules of rank two, and the given
map is an isomorphism modulo~$\pi'$.
It follows immediately from~(\ref{eqn: decomposing principal series type over powers
  of pi prime}) that~$\gL_\eta=U_\eta\oplus\pi_2V_{\eta'}$
and~$\gL_{\eta'}=\pi_1U_\eta\oplus V_{\eta'}$, so that~$\coker f_1$ and~$\coker f_2$ are projective of rank
one, as claimed.

\begin{prop}
  \label{prop: principal series case local model isomorphic to Iwahori
  local model}The morphism $\cM^\tau_{\loc}\to\cM_{\loc,\Iw}$ 
induces an isomorphism $\cM^{\tau,\BT}_{\loc}\to\cM^{\BT}_{\loc,\Iw}$.
{\em (}We remind the reader that $K'/K$ is assumed totally ramified,
and that $\tau$ is assumed to be a principal series inertial type.{\em )}
\end{prop} 
\begin{proof} 
We begin by constructing a morphism
  $\cM_{\loc,\Iw}\to\cM^\tau_{\loc}$, inspired by the arguments
  of~\cite[App.\ A]{rapoport-zink}. 
We define an~$\cO_K\otimes_{\Zp}A$-module~$\gL$ by
 \[\gL= \oplus_\sigma\left(\oplus_{i=0}^{a_\sigma-1}e_\sigma
    \gL_1^{(i)} \oplus_{j=0}^{b_\sigma-1}e_\sigma
    \gL_2^{(j)}\right), \] where the $\gL_1^{(i)}$'s and
  $\gL_2^{(j)}$'s are copies of~$\gL_1$, $\gL_2$ respectively. 

We can
  upgrade the ~$\cO_K\otimes_{\Zp}A$-module structure on~$\gL$ to that
  of an~$\cO_{K'}\otimes_{\Zp}A$-module by specifying how~$\pi'$
  acts.  If $i<a_\sigma-1$, then we let~$\pi':e_\sigma
    \gL_1^{(i)}\to e_\sigma
    \gL_1^{(i+1)}$ be the map induced by the identity
    on~$\gL_1$, and if~$j<b_\sigma-1$, then we let~$\pi':e_\sigma
    \gL_2^{(j)}\to e_\sigma
    \gL_2^{(j+1)}$ be the map induced by the identity
    on~$\gL_2$. We let $\pi':e_\sigma
    \gL_1^{(a_\sigma-1)}\to e_\sigma
    \gL_2^{(0)}$ be the map induced by~$f_1:\gL_1\to\gL_2$, and we let $\pi':e_\sigma
    \gL_2^{(b_\sigma-1)}\to e_\sigma
    \gL_1^{(0)}$ be the map induced by~$f_2:\gL_2\to\gL_1$. That this indeed gives~$\gL$ the structure of an
$\cO_{K'}\otimes_{\Zp}A$-module follows from our assumption that
$f_1\circ f_2=f_2\circ f_1=\pi$. We give~$\gL$ a semilinear action of~$\Gal(K'/K)=I(K'/K)$ 
by letting
it act via~$(\sigma \circ h)^{i} \cdot \eta$ on each~$e_\sigma
    \gL_1^{(i)} $ and via~$(\sigma \circ h)^j  \cdot \eta'$ on each~$e_\sigma
    \gL_2^{(j)} $.

We claim that~$\gL$ is a rank~$2$ projective
$\cO_{K'}\otimes_{\Zp}A$-module.  Since~$\coker f_2$ is projective by
  assumption, we can choose a section to the $k\otimes_{\Zp}A$-linear
  morphism $\gL_1/\pi\to\coker f_2$, with image~$\overline{U}_\eta$,
  say. Similarly we choose a section to $\gL_2/\pi\to\coker f_1$ with
  image~$\overline{V}_{\eta'}$. We choose lifts~$U_\eta$, $V_{\eta'}$
  of~$\overline{U}_\eta$, $\overline{V}_{\eta'}$ to direct summands of the
  $\cO_K\otimes_{\Zp}A$-modules $\gL_1$, $\gL_2$ respectively. There
  is a map of $\cO_{K'} \otimes_{\Zp} A$-modules
\numequation\label{eqn:
    reconstructing
    L} \lambda: \cO_{K'} \otimes_{\cO_K} (U_\eta \oplus V_{\eta'})  \to \gL
\end{equation}
induced by the map identifying $U_\eta, V_{\eta'}$ with their copies in $\gL_1^{(0)}$ and
$\gL_2^{(0)}$ respectively.  The map $\lambda$ is surjective modulo $\pi'$ by
construction, hence surjective by Nakayama's lemma. Regarding $\lambda$ as a map
of projective $\cO_K \otimes_{\Zp} A$-modules of equal rank $2e(K'/K)$, we
deduce that $\lambda$ is an isomorphism. Since the source of $\lambda$
is a projective $\cO_{K'}\otimes_{\Zp} A$-module, the claim follows.

We now set \[\gL^+=
  \oplus_\sigma\left(\oplus_{i=0}^{a_\sigma-1} e_\sigma
    (\gL_1^{(i)})^+ \oplus_{j=0}^{b_\sigma-1} e_\sigma
    (\gL_2^{(j)})^+\right) \subset \gL\]
It is
immediate from the construction that $\gL^+$ is preserved by
$\Gal(K'/K)$. The hypothesis that $f_1,f_2$ and preserve $\gL_1^+,\gL_2^+$ implies
that $\gL^+$ is an $\cO_{K'}\otimes_{\Zp} A$-submodule of $\gL$, while
the hypothesis that each $\gL_i/\gL_i^+$  is a projective $A$-module
implies the same for $\gL/\gL^+$.  This completes the construction of
our morphism 
  $\cM_{\loc,\Iw}\to\cM^\tau_{\loc}$.

Just as in the proof of Proposition~\ref{lem: local model scalar type is equivalent to what you
    expect}, the morphism $\cM_{\loc,\Iw}\to
\cM^{\tau}_{\loc}$ followed by our morphism
$\cM^{\tau}_{\loc}\to\cM_{\loc,\Iw}$ is the identity, while the composition in the other order is not, in general, naturally equivalent to the
identity morphism. 
However, it follows immediately from Lemma~\ref{lem: local model of type only have to check the strong
      determinant condition on that type} and the construction of $\gL^+$ that our morphisms $\cM_{\loc,\Iw}\to
\cM^{\tau}_{\loc}$ and
$\cM^{\tau}_{\loc}\to\cM_{\loc,\Iw}$ respect the strong determinant
condition, and so induce maps  $\cM_{\loc,\Iw}^{\BT} \to
\cM^{\tau,\BT}_{\loc}$ and
$\cM^{\tau,\BT}_{\loc} \to\cM_{\loc,\Iw}^{\BT}$.  To see that the
composite $\cM^{\tau,\BT}_{\loc} \to\cM_{\loc,\Iw}^{\BT} \to
\cM^{\tau,\BT}_{\loc}$ is naturally equivalent to the identity,
suppose that 
$(\gL,\gL^+) \in \cM^\tau_{\loc}(A)$ and observe that
there is a natural $\Gal(K'/K)$-equivariant isomorphism of $\cO_K\otimes_A\Zp$-modules
\numequation \label{eq:ps-natural-isom} \oplus_\sigma\left(\oplus_{i=0}^{a_\sigma-1}e_\sigma
    \gL_\eta^{(i)} \oplus_{j=0}^{b_\sigma-1}e_\sigma
    \gL_{\eta'}^{(j)}\right) \toisom \gL \end{equation}
induced by the maps $\gL_\eta^{(i)} \longnto{(\pi')^i} (\pi')^i
\gL_\eta$ and $\gL_{\eta'}^{(j)} \longnto{(\pi')^j} (\pi')^j
\gL_{\eta'}$. The commutativity of the diagram
\[ 
\xymatrix{ e_\sigma \gL_{\eta}^{(a_\sigma-1)} \ar[d]_{\pi_1}
  \ar[r]^{(\pi')^{a_\sigma-1}} & e_\sigma (\pi')^{a_\sigma-1} \gL_{\eta}
  \ar[d]^{\pi'} \\
e_\sigma \gL^{(0)}_{\eta'} \ar[r]_{\textrm{id}} & e_\sigma\gL_{\eta'} }\]
implies that the map in \eqref{eq:ps-natural-isom} is in fact an
$\cO_{K'}\otimes_{\Zp} A$-module isomorphism. The map
\eqref{eq:ps-natural-isom} induces an inclusion
\[ \oplus_\sigma\left(\oplus_{i=0}^{a_\sigma-1}e_\sigma
    (\gL_\eta^{(i)})^+ \oplus_{j=0}^{b_\sigma-1}e_\sigma
    (\gL_{\eta'}^{(j)})^+ \right) \rightarrow \gL^+.\]
If furthermore $(\gL,\gL^+) \in \cM^{\tau,\BT}_{\loc}$ then this is an isomorphism because $\gL^+$ satisfies condition (2) of Lemma~\ref{lem: local model of type only have to check the strong
      determinant condition on that type}.
  \end{proof}

Finally, we turn to the case of a cuspidal type. Let~$L$ as usual be a
quadratic unramified extension of~$K$, and
set~$K'=L(\pi^{1/(p^{2f}-1)})$. The field $N$ continues to denote the
maximal unramified extension of $K$ in $K'$, so that $N=L$.  Let~$\tau$ be a cuspidal type, so that
~$\tau=\eta\oplus\eta'$, where~$\eta\ne\eta'$ but~$\eta'=\eta^{p^f}$.

\begin{prop}
  \label{prop: cuspidal local model formally smooth over Iwahori}There
  is a morphism 
  $\cM_{\loc}^{\tau,\BT}\to\cM_{\loc,\Iw}^{\BT}$
which is representable by algebraic spaces and smooth.
{\em (}We remind the reader that $\tau$ is now assumed to be a cuspidal
inertial type.{\em )}
\end{prop}
\begin{proof}Let~$\tau'$ be the type~$\tau$, considered as a (principal series) type for the
  totally ramified extension~$K'/N$. Let~$c\in\Gal(K'/K)$ be the
  unique element which fixes~$\pi^{1/(p^{2f}-1)}$ but acts
  nontrivially on~$N$. For any map $\alpha : X\to Y$ of
  $\cO_N$-modules we write $\alpha^c$ for the twist  $1\otimes \alpha :
  \cO_N \otimes_{\cO_N,c} X \to \cO_N \otimes_{\cO_N,c} Y$.

We may think of an
  object~$(\gL,\gL^+)$ of~$\cM_{\loc}^{\tau,\BT}$ as an object~$(\gL',(\gL')^+)$
  of~$\cM_{\loc}^{\tau',\BT}$ equipped with the additional data of an
  isomorphism of $\cO_{K'}\otimes_{\Zp}A$-modules
  $\theta:\cO_{K'}\otimes_{\cO_{K'},c}\gL'\to\gL'$ which is compatible
  with~$(\gL')^+$, which satisfies~$\theta\circ\theta^c=\id$, and which is compatible with the
  action of~$\Gal(K'/N)=I(K'/N)$ in the sense that $\theta \circ (1
  \otimes g) = g^{p^f} \circ \theta$.

  Employing the isomorphism of Proposition~\ref{prop: principal series case local model isomorphic to Iwahori
  local model}, we think of~$(\gL',(\gL')^+)$ as a
tuple~$(\gL'_1, (\gL'_1)^+,\gL'_2,(\gL'_2)^+,f_1,f_2)$, where the~$\gL'_i,(\gL'_i)^+$ are
$\cO_N\otimes_{\Zp}A$-modules;\ by construction, the map~$\theta$
induces  isomorphisms $\theta_1:\cO_N\otimes_{\cO_N,c}\gL'_1\isoto\gL_2'$,
$\theta_2:\cO_N\otimes_{\cO_N,c}\gL'_2\isoto\gL'_1$, which are
compatible with~$(\gL'_1)^+,(\gL'_2)^+$ and~$f_1,f_2$, and
satisfy~$\theta_1\circ\theta_2^c=\id$.

Choose for each embedding~$\sigma:k\into\F$ an extension to an
embedding~$\sigma^{(1)}:k'\into\F$, 
set~$e_1=\sum_\sigma e_{\sigma^{(1)}}$, and write~$e_2=1-e_1$. Then
the map~$\theta_1$  induces isomorphisms
$\theta_{11}:e_1\gL'_1\isoto e_2\gL_2'$ and
$\theta_{12}:e_2\gL'_1\isoto e_1\gL_2'$, while~$\theta_2$  induces isomorphisms
$\theta_{21}:e_1\gL'_2\isoto e_2\gL_1'$ and
$\theta_{22}:e_2\gL'_2\isoto e_1\gL_1'$. The condition
that~$\theta_1\circ\theta_2^c=\id$ translates to $\theta_{22} =
\theta_{11}^{-1}$ and $\theta_{21}=\theta_{12}^{-1}$, and
compatibility with~$(\gL'_1)^+,(\gL'_2)^+$ implies that
$\theta_{11},\theta_{21}$ induce isomorphisms $e_1 (\gL'_1)^+ \isoto
e_2(\gL_2')^+$ and $e_1 (\gL'_2)^+ \isoto
e_2(\gL_1')^+$ respectively.

Furthermore $f_1,f_2$ induce maps
$e_1f_1:e_1\gL'_1\to e_1\gL_2'$, $e_1g:e_1\gL_2'\to e_1\gL_1'$.
 It
follows that there is a map $\cM_{\loc}^{\tau,\BT}\to\cM_{\loc,\Iw}$,
sending~$(\gL,\gL^+)$ to the tuple
$(e_1\gL_1',e_1(\gL_1')^+,e_1\gL_2',e_1(\gL_2')^+,e_1f_1,e_1f_2)$.
To see that it respects the strong determinant condition, 
one has to check that the
    conditions on~$\gL^+_1$, $\gL^+_2$ imply those for~$e_1(\gL'_1)^+$,
    $e_2(\gL'_2)^+$
    coincide, and this follows from the definitions (via Remark~\ref{rem: explicit version of determinant condition
      local models version}). We therefore obtain a map $\cM_{\loc}^{\tau,\BT}\to\cM^{\BT}_{\loc,\Iw}$,

Since this morphism is
given by forgetting the data of $e_2 \gL'_1,e_2 \gL'_2$ and the pair of isomorphisms~$\theta_{11},\theta_{21}$, it is evidently formally
smooth.
It is also a morphism between $\varpi$-adic formal algebraic stacks that
are locally of finite presentation, and so is representable by algebraic spaces
(by \cite[Lem.~7.10]{Emertonformalstacks}) and locally of finite presentation
(by \cite[Cor.~2.1.8]{EGstacktheoreticimages}
and~\cite[\href{http://stacks.math.columbia.edu/tag/06CX}{Tag
  06CX}]{stacks-project}).  Thus this morphism is in fact smooth.
\end{proof}

\subsection{Local models:\ local geometry}\label{subsec: local models
  relation to explicit PR ones}We now deduce our main results on the
local structure of our moduli stacks from results in the literature on local models for Shimura
varieties.

  \begin{prop}\label{prop: explicit description of local models stack
      as naive local model} 
  We can identify $\cM_{\loc,\Iw}^{\BT}$ with the quotient 
  of \emph{(}the $p$-adic formal completion of\emph{)} the naive local model for $\Res_{K/\Q_p}\GL_2$
  \emph{(}as defined in~\cite[\S 2.4]{MR3135437}\emph{)} by a smooth
  group scheme over $\cO$.
  \end{prop}
  
  \begin{proof} Let $\widetilde{\cM}_{\loc,\Iw}^{\BT}$ be the $p$-adic formal
    completion of the naive local model for $\Res_{K/\Q_p}\GL_2$
  corresponding to a standard lattice chain~ $\cL$, as defined in~\cite[\S 2.4]{MR3135437}. 
  By~\cite[Prop. A.4]{rapoport-zink}, the automorphisms of the standard lattice chain $\cL$ 
  are represented by a smooth group scheme $\cP_{\cL}$ over $\cO$. (This is in fact a parahoric subgroup scheme
  of $\mathrm{Res}_{\cO_K/\Z_p}\GL_2$, and in particular it is affine.) 
  Also by~\emph{loc.\ cit.}, every lattice chain of type $(\cL)$ is Zariski locally isomorphic to $\cL$. 
  By comparing the two moduli problems, we see
  that $\widetilde{\cM}_{\loc,\Iw}^{\BT}$ is a $\cP_{\cL}$-torsor over
  $\cM_{\loc,\Iw}^{\BT}$ 
  for the Zariski topology
  and the proposition follows. 
  \end{proof}
  
The following theorem describes the various regularity properties of local models.
Since we are working in the context of formal algebraic stacks,
we use the terminology developed in~\cite[\S 8]{Emertonformalstacks}
(see in particular~\cite[Rem.~8.21]{Emertonformalstacks}
and~\cite[Def.~8.35]{Emertonformalstacks}).

\begin{thm} 
  \label{thm: local consequences of local models}Suppose that~$d=2$
  and that $\tau$
  is a tame inertial type. Then
  \begin{enumerate}
  \item $\cM_{\loc}^{\tau,\BT}$ is residually Jacobson and analytically normal,
 and Cohen--Macaulay.
  \item The special fibre $\cM_{\loc}^{\tau,\BT,1}$ is reduced.
  \item  $\cM_{\loc}^{\tau,\BT}$ is flat over~$\cO$.
  \end{enumerate}

\end{thm}
\begin{proof}
For scalar types, this follows from~\cite[Prop.\
  2.2.2]{kis04} by Lemma~\ref{lem: local model scalar type is equivalent to what you
    expect}, 
and so we turn to studying the case of a non-scalar type.
The properties in question can be checked smooth locally
(see \cite[\S 8]{Emertonformalstacks} for~(1),
and \cite[\href{http://stacks.math.columbia.edu/tag/04YH}{Tag 04YH}]{stacks-project} for~(2);
for~(3), note that morphisms that are representable by algebraic spaces and smooth
are also flat, and take into account the final statement
of~\cite[Lem.~8.34]{Emertonformalstacks}),
and so 
  by Propositions~\ref{prop: principal series case local model isomorphic to Iwahori
  local model} and~\ref{prop: cuspidal local model formally smooth over
  Iwahori}
we reduce to checking the assertions of the theorem
for $\cM^{\BT}_{\loc,\Iw}.$  
Proposition~\ref{prop: explicit description of local models stack
    as naive local model}
then reduces us to checking these assertions for the $\varpi$-adic completion
of the naive local model at Iwahori level for $\Res_{K/\Q_p}\GL_2$.

Since this naive local model is a scheme of finite presentation over $\cO$,
its special fibre (i.e.\ its base-change to $\F$) is Jacobson,
and it is excellent;\ thus its $\varpi$-adic completion satisfies the properties
of~(1) if and only if the naive local model itself is normal and Cohen--Macaulay.
The special fibre of its $\varpi$-adic completion is of course just equal
to its own special fibre, and so verifying~(2) for the first of these special fibres
is equivalent to verifying it for the second.  Finally, the $\varpi$-adic completion
of an $\cO$-flat Noetherian ring is again $\cO$-flat, and so the $\varpi$-adic
completion of the naive local model will be $\cO$-flat if the naive local
model itself is. 

All the properties of the naive local model that are described in the preceding
paragraph,
other than the Cohen--Macaulay property, 
 are contained 
  in~\cite[Thm.\ 2.4.3]{MR3135437},
  if we can identify the naive local models at Iwahori level with the vertical local
  models at Iwahori level, in the sense of~\cite[\S 8]{pappas-rapoport1} 
  (see also the discussion above~\emph{loc.\ cit.}\
  in~\cite{MR3135437}).
  
  The vertical local models are obtained
  by intersecting the preimages at Iwahori level of the flat local models at hyperspecial level. Although naive local models 
  are not flat in general,  for the Weil restriction of~$\GL_2$
  at hyperspecial level they are, by a special case
  of~\cite[Cor.\ 4.3]{pappas-rapoport1}. Thus the naive local models at Iwahori level
  are identified with the vertical ones and~\cite[Thm.\
  2.4.3]{MR3135437} applies to them directly.  (To be precise, the results
  of~\cite{pappas-rapoport1} apply to restrictions of scalars
  $\Res_{F/F_0}\GL_2$ with $F/F_0$ totally ramified. However, thanks
  to the decomposition $\cO_K \otimes_{\Zp} A \cong \oplus_{\sigma :
    W(k) \into \cO} \cO_K \otimes_{W(k),\sigma} A$ the local
  model for
  $\Res_{K/\Qp}\GL_2$ decomposes as a product of local models for totally
  ramified extensions.) 

  Finally, Cohen--Macaulayness can be proved as in~\cite[Prop.\
  4.24]{MR1871975} and the discussion immediately following. We thank
  U.\ G\"ortz for explaining this argument to us.  As in the
  previous paragraph we reduce to the case of local models at Iwahori
  level for $\Res_{F/F_0}\GL_2$ with $F/F_0$ totally ramified of
  degree~$e$. In this
  setting the admissible set $M$ for the coweight $\mu = (e,0)$ has precisely
   one element of length $0$ and two elements of each length between $1$
   and $e$.  
Moreover, for elements $x,y \in M$ we have
$x<y$ in the Bruhat order if and only if $\ell(x)<\ell(y)$. One checks easily
that $M$ is $e$-Cohen--Macaulay in the sense of~\cite[Def.\
  4.23]{MR1871975}, and we conclude by~\cite[Prop.\
  4.24]{MR1871975}. Alternatively, this also follows from the much more
  general recent results of Haines--Richarz \cite{HainesRicharz}.
\end{proof}

\begin{cor} \label{cor: Kisin moduli consequences of local models}Suppose that~$d=2$
  and that $\tau$
  is a tame inertial type. Then
  \begin{enumerate}
  \item $\cC^{\tau,\BT}$ is analytically normal, and Cohen--Macaulay.
  \item The special fibre $\cC^{\tau,\BT,1}$ is reduced.
  \item  $\cC^{\tau,\BT}$ is flat over~$\cO$.
  \end{enumerate}
  \end{cor}
\begin{proof}
This follows from Theorems~\ref{thm: map from Kisin to general local model is
    formally smooth} and~\ref{thm: local consequences of local models},
since all of these properties can be verified smooth locally
(as was already noted in the proof of the second of these theorems).
\end{proof}

\begin{remark}
  There is another structural result about vertical local models
that is proved in~\cite[\S 8]{pappas-rapoport1} but which we haven't incorporated
into our preceding results, namely, that the irreducible components of
the special fibre are normal.  Since the notion of irreducible component
is not \'etale local (and so in particular not smooth local), this statement
does not imply the corresponding statement for the special fibres of
$\cM_{\loc}^{\tau,\BT}$ or $\cC^{\tau,\BT}$.  Rather, it implies the weaker,
and somewhat
more technical, statement that each of the analytic branches passing through each
closed point of the special fibre of these $\varpi$-adic formal algebraic
stacks is normal.  We won't discuss this further here, since we don't 
need this result.
\end{remark}

\subsection{The Dieudonn\'e stack}\label{subsec: Dieudonne
  stack}

In this subsection we assume that
$\eta\ne\eta'$.  As previously mentioned, the construction described in this
section will not be needed elsewhere in the present paper;\ however, it is
convenient to include it here because of its reliance on the discussion of
determinant conditions in Section~\ref{subsec: local
  models rank two}.

Let $\gM$ be a Breuil--Kisin module with $A$-coefficients and descent data of type
$\tau$ and height at most $1$, and let $D := \gM/u\gM$ be
its corresponding Dieudonn\'e module as in Definition~\ref{def: Dieudonne module formulas}. 
If we write
$D_i:=e_iD$, then this Dieudonn\'e module is given by rank two
projective modules
 $D_j$ over $A$ ($j = 0,\ldots, f'-1$) with linear maps $F: D_j \to D_{j+1}$
and $V: D_j \to D_{j-1}$ (subscripts understood modulo $f'$) 
such that $FV = VF = p$.

Now,
  $I(K'/K)$ is abelian of order prime to $p$, so we can write
  $D=D_\eta\oplus D_{\eta'}$, where $D_\eta$ is the submodule on which
  $I(K'/K)$ acts via~$\eta$. Since $\gM_\eta$ is obtained from the projective $\gS_A$-module $\gM$
by applying a projector, each $D_{\eta,j}$ is an invertible
$A$-module, and $F,V$ induce linear
maps $F:D_{\eta,j}\to D_{\eta,j+1}$ and $V: D_{\eta,j+1} \to D_{\eta,j}$ 
such that $FV = VF = p$. 

We can of course apply the same construction with $\eta'$ in the place
of $\eta$, obtaining a Dieudonn\'e module $D_{\eta'}$.
We now prove some lemmas relating these various Dieudonn\'e
modules. We will need to make use of a variant of the strong
determinant condition, so we begin by discussing this and its relationship to
the strong determinant condition of Subsection~\ref{subsec: local
  models rank two}. 

\begin{defn}
  \label{defn: determinant condition over L} 
Let~$(\gL,\gL^+)$ be a
  pair consisting of a rank two projective
  $\cO_{K'}\otimes_{\Zp}A$-module $\gL$, and an
  $\cO_{K'}\otimes_{\Zp}A$-submodule $\gL^+ \subset \gL$, such that
  Zariski locally on~$\Spec A$, $\gL^+$ is a direct summand of~$\gL$
  as an $A$-module. 

Then we say that the pair~$(\gL,\gL^+)$ \emph{satisfies the Kottwitz
  determinant condition over~$K'$} if for all~$a\in\cO_{K'}$, we
have \[\det{\!}_A(a|\gL^+)=\prod_{\psi:K'\into E}\psi(a) \]as polynomial functions on $\cO_{K'}$ in the sense of~\cite[\S
  5]{MR1124982}.
\end{defn}

There is a finite type stack~$\cM_{K',\det}$ over~$\Spec \cO$, with
$\cM_{K',\det}(\Spec A)$ being the groupoid of pairs~$(\gL,\gL^+)$ as
above which satisfy the Kottwitz determinant condition over~$K'$. As we
have seen above, by a
result of Pappas--Rapoport, this
stack is flat over~$\Spec\cO$ (see~\cite[Prop.\ 2.2.2]{kis04}).

\begin{lem}
  \label{lem: determinant condition  generic fibre}If $A$ is an
  $E$-algebra, then a pair~$(\gL,\gL^+)$ as in Definition~{\em \ref{defn: determinant condition over L}} satisfies the Kottwitz
  determinant condition over~$K'$ if and only if~$\gL^+$ is a rank one
  projective $\cO_{K'}\otimes_{\Zp}A$-module.
\end{lem}
\begin{proof}
  We may write
  $\cO_{K'}\otimes_{\Zp}A=K'\otimes_{\Qp}A\cong\prod_{\psi:K'\into E}A$, where
  the embedding $\psi:K'\into E$ corresponds to an idempotent
  $e_\psi\in K'\otimes_{\Qp}A$. Decomposing $\gL^+$ as $\oplus_{\psi}
  e_{\psi} \gL^+$, the left-hand side of the Kottwitz determinant condition
becomes $\prod_{\psi} \det{\!}_A(a|e_\psi \gL^+) = \prod_{\psi}
\psi(a)^{\mathrm{rk}_A e_{\psi} \gL^+}$. It follows that the Kottwitz
  determinant condition is satisfied if and only if the projective
  $A$-module $e_\psi
  \gL^+$ has rank one for all $\psi$, which is equivalent to~$\gL^+$ being a
  rank one projective~$K'\otimes_{\Qp}A$-module, as required.
\end{proof}

\begin{prop}
  \label{prop: for us strong det implies det}If~$\gM$ is an object
  of~$\cC^{\tau,\BT}(A)$, 
  then
  the pair \[(\gM/E(u)\gM,\im\Phi_{\gM}/E(u)\gM)\] satisfies the
  Kottwitz determinant condition for~$K'$.
\end{prop}
\begin{proof}Let~$\cC^{\tau,\BT'}$ be the closed substack
  of~$\cC^{\tau}$ consisting of those~$\gM$ for which the
  pair~$(\gM/E(u)\gM,\im\Phi_{\gM}/E(u)\gM)$ satisfies the Kottwitz
  determinant condition for~$K'$. We need to show
  that~$\cC^{\tau,\BT}$ is a closed substack
  of~$\cC^{\tau,\BT'}$. Since~$\cC^{\tau,\BT}$ is flat over~$\Spf\cO$
  by Corollary~\ref{cor: Kisin moduli consequences of local models},
  it is enough to show that if~$A^\circ$ is a flat $\cO$-algebra, then
  $\cC^{\tau,\BT}(A^\circ)=\cC^{\tau,\BT'}(A^\circ)$. 
  
To see this, let~$\gM^\circ$ be an object
of~$\cC^{\tau}(A^\circ)$. Since~$A^\circ$ is flat over~$\cO$, the Kottwitz determinant conditions
for~$\gM$ to be an object of~$\cC^{\tau,\BT'}(A^\circ)$ or
of~$\cC^{\tau,\BT}(A^\circ)$ can be checked after
inverting~$p$. Accordingly, we write $A=A^\circ[1/p]$, and write $\gM$ for
the base change of~$\gM^\circ$ to~$A$.  By Lemma~\ref{lem: determinant condition
  generic fibre}, $\gM^\circ$ is an object of~$\cC^{\tau,\BT'}(A^\circ)$ if and
only if $\im\Phi_{\gM}/E(u)\gM$ is a rank one
projective~$K'\otimes_{\Qp}A$-module. Similarly,  $\gM^\circ$ is an object of~$\cC^{\tau,\BT}(A^\circ)$ if and
only if for each~$\xi$, $(\im\Phi_{\gM})_{\xi}/E(u)\gM_{\xi}$ is a rank one
projective~$N\otimes_{\Qp}A$-module. Since \[\im\Phi_{\gM}/E(u)\gM=\oplus_{\xi}(\im\Phi_{\gM})_{\xi}/E(u)\gM_{\xi},\]the
equivalence of these two conditions is clear.
\end{proof}

\begin{lemma}
\label{lem:dets for local models}
If $(\gL,\gL^+)$ is an object of $\cM_{K',\det}(A)$ {\em (}i.e.\
satisfies the Kottwitz determinant condition over $K'${\em )}, 
then the morphism $\bigwedge^2_{\cO_{K'}\otimes_{\Z_p} A} \gL^+ \to
\bigwedge^2_{\cO_{K'}\otimes_{\Z_p} A} \gL$ induced
by the inclusion $\gL^+ \subset \gL$ is identically zero.
\end{lemma}

\begin{remark}
Note that, although $\gL^+$ need not be locally free over
$\cO_{K'}\otimes_{\Z_p} A$, its exterior square is nevertheless defined,
so that the statement of the lemma makes sense.
\end{remark}

\begin{proof}[Proof of Lemma~{\ref{lem:dets for local models}}]
Since~$\cM_{K',\det}$ is $\cO$-flat, 
it is enough to treat the case that $A$ is $\cO$-flat. In this case $\gL$, and thus
also $\bigwedge^2 \gL$, are $\cO$-flat.  Given this additional assumption,
it suffices to prove that the morphism of the lemma becomes zero after tensoring
with $\Q_p$ over $\Z_p$.   This morphism may naturally be identified with
the morphism
$\bigwedge^2_{K'\otimes_{\Z_p} A} \gL^+ \to \bigwedge^2_{K'\otimes_{\Z_p} A} \gL$
induced by the injection
$\Q_p\otimes_{\Z_p} \gL^+ \hookrightarrow \Q_p\otimes_{\Z_p} \gL$.
Locally on $\Spec A$, this is the embedding of a free $K'\otimes_{\Z_p} A$-module
of rank one as a direct summand of a free $K'\otimes_{\Z_p} A$-module of rank
two.  Thus $\bigwedge^2$ of the source in fact vanishes, and hence so does $\bigwedge^2$ of the embedding.
\end{proof}

\begin{lemma}
\label{lem:BT condition gives control of det Phi}
If $\gM$ is an object of $\cC^{\tau,\BT}(A)$, then
$\bigwedge^2 \Phi_{\gM}: \bigwedge^2 \varphi^*\gM  \hookrightarrow \bigwedge^2 \gM$
is exactly divisible by $E(u)$, i.e.\ can be written as $E(u)$ times an
isomorphism of $\gS_A$-modules.
\end{lemma}
\begin{proof}
It follows from Proposition~\ref{prop: for us strong det implies det} and
Lemma~\ref{lem:dets for local models} that the reduction of
$\bigwedge^2 \Phi_{\gM}$ modulo $E(u)$ vanishes, so we can think
of~$\bigwedge^2\Phi_{\gM}$ as a morphism
$\bigwedge^2\varphi^*\gM\to E(u)\bigwedge^2\gM$.
We need to show that the cokernel~$X$ of this
morphism vanishes. Since~$\im\Phi_{\gM}\supseteq E(u)\gM$, $X$ is a finitely generated $A$-module, so that in order to prove
that it vanishes, it is enough to prove that~$X/pX=0$. 

Since the formation of cokernels is compatible with
base change, this means that we can (and do) assume that~$A$ is an
$\F$-algebra. Since the special fibre~$\cC^{\tau,\BT}$ is of finite
type over~$\F$, we can and do assume that~$A$ is furthermore of finite
type over~$\F$. 
The special fibre of~$\cC^{\tau,\BT}$ is reduced by
Corollary~\ref{cor: Kisin moduli consequences of local models}, 
so we
may assume that~$A$ is reduced, 
and it is therefore enough to prove that~$X$ vanishes
modulo each maximal ideal of~$A$. Since the residue fields at such
maximal ideals are finite, 
we are reduced to the case that~$A$ is a finite
field, when the result follows from~\cite[Lem.\ 2.5.1]{kis04}. 
\end{proof}

\begin{lemma}
\label{lem:Dieudonne iso}
There is a canonical isomorphism
$$\text{``}(F\otimes F)/p\text{''}
:D_{\eta,j} \otimes_A
D_{\eta',j} \iso D_{\eta,j+1} \otimes_A D_{\eta',j+1},$$
characterised by the fact that it is compatible with change of scalars,
and that $$p \cdot \text{``}(F\otimes F)/p\text{''}= F\otimes F.$$
\end{lemma}
\begin{proof}
Since $\cC^{\BT}$ is flat over $\cO$, we see that in the universal case,
the formula $$p \cdot \text{``}(F\otimes F)/p\text{''}= F\otimes F$$
uniquely determines the isomorphism
$\text{``}(F\otimes F)/p\text{''}$ (if it exists).  Since any 
Breuil--Kisin module with descent data is obtained from the universal case by
change of scalars, we see that the isomorphism
$\text{``}(F\otimes F)/p\text{''}$ is indeed characterised by the properties
stated in the lemma, provided that it exists. 

To check that the isomorphism exists, we can again consider the universal case,
and hence assume that $A$ is a flat $\cO$-algebra. 
In this case, it suffices to check that the morphism
$F\otimes F: D_{\eta,j}\otimes_A D_{\eta',j} \to D_{\eta,j+1}\otimes_A D_{\eta,
j+1}$ is divisible by $p$, and that the formula
$(F\otimes F)/p$ is indeed an isomorphism.
Noting that the direct sum over $j = 0,\ldots, f'-1$ of these morphisms may
be identified with the reduction modulo $u$ of the morphism
$\bigwedge^2 \Phi_{\gM}: \bigwedge^2 \varphi^*\gM \to \bigwedge^2
\gM$, 
this follows from Lemma~\ref{lem:BT condition gives control of det Phi}.
\end{proof}

The isomorphism $\text{``}(F\otimes F)/p \text{''}$ of the preceding lemma
may be rewritten as an isomorphism of invertible $A$-modules
\numequation
\label{eqn:Dieudonne hom iso}
\Hom_A(D_{\eta,j}, D_{\eta,j+1}) \iso \Hom_A(D_{\eta',j+1},D_{\eta',j}).
\end{equation}

\begin{lem}
\label{lem:swapping F and V}
The isomorphism~{\em (\ref{eqn:Dieudonne hom iso})} takes $F$ to $V$.
\end{lem}
\begin{proof}
The claim of the lemma is equivalent to showing that
the composite 
$$D_{\eta,j} \otimes_A D_{\eta',j+1} \buildrel
\id\otimes V \over \longrightarrow 
D_{\eta,j}\otimes_A D_{\eta',j} \buildrel \text{``}(F\otimes F)/p \text{''}
\over \longrightarrow D_{\eta,j+1}\otimes_A D_{\eta',j+1}$$
coincides with the morphism $F \otimes \id.$  
It suffices to check this in the universal case, and thus we may assume that $p$ is a non-zero divisor in $A$, and hence verify the required identity of morphisms after multiplying each of them by $p$.
The identity to be verified then becomes
$$(F\otimes F) \circ (\id \otimes V) \buildrel ? \over = p (F\otimes \id),$$
which follows immediately from the formula $F V = p$.
\end{proof}

We now consider the moduli stacks classifying the Dieudonn\'e modules
with the properties we have just established, and the maps from the
moduli stacks of Breuil--Kisin modules to these stacks.   In this
discussion we specialise the choice of $K'$ in the following way. Choose a
tame inertial type $\tau=\eta\oplus\eta'$ with $\eta \neq \eta'$.
Fix a uniformiser $\pi$ of~$K$. If $\tau$ is a tame
principal series type, we take $K'=K(\pi^{1/(p^f-1)})$, while
if~$\tau$ is a tame cuspidal type, we let $L$ be an unramified
quadratic extension of~$K$, and set $K'=L(\pi^{1/(p^{2f}-1)})$. Let
$N$ be the maximal unramified extension of $K$ in $K'$. In
either case $K'/K$ is a Galois extension;\ in the principal series
case, we have $e'=(p^f-1)e$, $f'=f$, and in the cuspidal case we have
$e'=(p^{2f}-1)e$, $f'=2f$.

Suppose first that we are in the principal series case. Then there is a moduli stack classifying the data of the $D_{\eta,j}$
together with the
$F$ and $V$, namely the stack 
$$\cD_{\eta} :=
\Big[
\bigl( \Spec W(k)[X_0,Y_0,\ldots,X_{f-1},Y_{f-1}]/(X_j Y_j - p)_{j = 0,\ldots, f-1}) \bigr) /
\mathbb G_m^f  \big],$$
where the $f$ copies of $\mathbb G_m$ act as follows:
$$(u_0,\ldots,u_{f-1}) \cdot (X_j,Y_j) \mapsto (u_j u_{j+1}^{-1} X_j, u_{j+1} u_j^{-1} Y_j).$$
To see this, recall that the stack
$$[\text{point}/\Gm]$$ classifies line bundles, so the $f$ copies of $\mathbb G_m$
in $\cD_\eta$ correspond to $f$ line bundles, which are the line bundles $D_{\eta,j}$
($j = 0,\ldots,f-1$).  If we locally trivialise these line bundles, then the
maps $F: D_{\eta,j} \to D_{\eta,j+1}$ and $V:D_{\eta,j+1} \to D_{\eta,j}$ act by
scalars, which we denote by $X_j$ and $Y_j$ respectively.  The $f$ copies of
$\mathbb G_m$ are then encoding possible changes of trivialisation, by units
$u_j$, which induce the indicated changes on the $X_j$'s and $Y_j$'s.

There is then a natural map
\[\cC^{\tau,\BT}
 \to \cD_{\eta},\]
classifying the Dieudonn\'e modules
underlying the Breuil--Kisin modules with descent data. 

There is a more geometric way to think about what $\cD_{\eta}$ 
classifies.
To begin with, we just rephrase what we've already indicated:
it represents the functor which associates to a $W(k)$-scheme the
groupoid whose objects are $f$-tuples of line bundles
$(D_{\eta,j})_{j = 0,\ldots,f-1}$ equipped with morphisms $X_j: D_{\eta,j} \to D_{\eta,j+1}$
and $Y_j:D_{\eta,j+1} \to D_{\eta,j}$ such that $Y_j X_j = p$.
(Morphisms in the groupoid are just isomorphisms between collections of such data.)
Equivalently, we can think of this as giving the line bundle
$D_{\eta,0}$, and then the $f$ line bundles $\cD_j
:= D_{\eta,j+1}\otimes D_{\eta,j}^{-1}$,
equipped with sections $X_j \in \cD_j$ and $Y_j \in \cD_j^{-1}$ whose
product in $\cD_j \otimes \cD_j^{-1} = \cO$
(the trivial line bundle) is equal to the element~$p$.
Note that it superficially looks like we are remembering $f+1$ line bundles,
rather than $f$, but this is illusory, since in fact
$\cD_0 \otimes \cdots \otimes \cD_{f-1}$ is trivial;\ indeed,
the isomorphism 
$\cD_0 \otimes \cdots \otimes \cD_{f-1} \iso \cO$ is part of the data we should remember.

We now turn to the case that~$\tau$ is a cuspidal type. In this case
our Dieudonn\'e modules have unramified as well as inertial descent
data;\ accordingly, we let $\varphi^f$ denote the element of
$\Gal(K'/K)$ which acts trivially on $\pi^{1/(p^{2f}-1)}$ and
non-trivially on $L$. Then the descent data of~ $\varphi^f$ induces
isomorphisms $D_j\isoto D_{j+f}$, which are compatible with the $F,V$,
and which identify $D_{\eta,j}$ with $D_{\eta',f+j}$. 

If we choose local trivialisations of the line bundles
$D_{\eta,0},\dots,D_{\eta,f}$, then the maps $F : D_{\eta,j} \to
D_{\eta,j+1}$ and $V : D_{\eta,j+1} \to D_{\eta,j}$ for $0 \le j \le
f-1$ are given by scalars $X_j$ and $Y_j$ respectively. The identification of $D_{\eta,j}$ and $D_{\eta',f+j}$ given
by~$\varphi^f$ identifies $D_{\eta,j}\otimes D_{\eta,j+1}^{-1}$ with
$D_{\eta',f+j}\otimes D_{\eta',f+j+1}^{-1}$, which via the isomorphsim
\eqref{eqn:Dieudonne hom iso} is identified with
$D_{\eta,f+j+1}\otimes D_{\eta,f+j}^{-1}$. It follows that for $0\le j\le f-2$ the
data of $D_{\eta,j}$, $D_{\eta,j+1}$ and $D_{\eta,f+j}$ recursively
determines $D_{\eta,f+j+1}$. From Lemma~\ref{lem:swapping F and V}
we see, again recursively  for $0\le j\le f-2$, that there are unique trivialisations
of $D_{\eta,f+1},\dots,D_{\eta,2f-1}$ such that
$F:D_{\eta,f+j}\to D_{\eta,f+j+1}$ is given by~$Y_{j}$, and
$V:D_{\eta,f+j+1}\to D_{\eta,f+j}$ is given by~$X_j$. Furthermore,
there is some unit~$\alpha$ such that $F:D_{\eta,2f-1}\to D_{\eta,0}$
is given by~$\alpha Y_{f-1}$, and
$V:D_{\eta,0}\to D_{\eta,2f-1}$ is given by~$\alpha^{-1}
X_{f-1}$. Note that the map $F^{2f} : D_{\eta,0} \to D_{\eta,0}$ is
precisely $ p^f\alpha$.

Consequently, we see that the data of the  $D_{\eta,j}$ (together with
the $F,V$) is classified by the stack  $$\cD_{\eta} :=
\Big[
\bigl( \Spec W(k)[X_0,Y_0,\ldots,X_{f-1},Y_{f-1}]/(X_j Y_j - p)_{j = 0,\ldots, f-1}) \times\Gm\bigr) /
\mathbb G_m^{f+1}  \big],$$
where the $f+1$ copies of $\mathbb G_m$ act as follows:
$$(u_0,\ldots,u_{f-1},u_f) \cdot ((X_j,Y_j),\alpha) \mapsto ((u_j u_{j+1}^{-1}
X_j, u_{j+1} u_j^{-1} Y_j),\alpha ). $$
One again there is a natural map
\[\cC^{\tau,\BT}
 \to \cD_{\eta},\]
classifying the Dieudonn\'e modules
underlying the Breuil--Kisin modules with descent data.  

\section{Galois moduli stacks}\label{sec:galois-moduli}

We are finally ready to discuss the algebraic stacks $\cZ^{\dd,a}$ and
$\cZ^{\tau,a}$, which are defined in Section~\ref{subsec:results from
  EG and PR on finite flat moduli} as the scheme-theoretic images of
the morphisms $\cC^{\dd,\BT,a}\to\cR^{\dd,a}$ and
$\cC^{\tau,\BT,a}\to\cR^{\dd,a}$ respectively. We informally think of
these stacks as being moduli stacks of two-dimensional representations
of~$G_K$.

In
Section~\ref{subsec:Galois deformation rings} we relate the versal
rings of the stacks $\cZ^{\tau}$ to potentially crystalline
deformation rings, and use this to show that the stacks $\cZ^{\dd,a}$,
$\cZ^{\tau,a}$, and $\cC^{\tau,\BT,a}$ are equidimensional of
dimension $[K:\Qp]$ (see Propositions~\ref{prop: dimensions
  of the Z stacks} and~\ref{prop: C tau is equidimensional of the
  expected dimension}).

\subsection{Scheme-theoretic images}\label{subsec:results from EG
  and PR on finite flat moduli}
We continue to fix $d=2$, $h=1$, and we set
$K'=L(\pi^{1/{p^{2f-1}}})$, where $L/K$ is the unramified quadratic
extension, and~$\pi$ is a uniformiser of~$K$. (This is the choice of
$K'$ that we made before in the cuspidal case, and contains the choice
of $K'$ that we made in the principal series case;\ since the category
of Breuil--Kisin
modules with descent data for the smaller extension is by Proposition~\ref{prop: change of K'} naturally a
full subcategory of the category of Breuil--Kisin modules with descent data
for the larger extension, we can consider both principal series and
cuspidal types in this setting.)

\begin{df}
\label{def:defining the Z's}
For each $a\ge 1$ we write $\cZ^{\dd,a}$ and $\cZ^{\tau,a}$ for the scheme-theoretic
images (in the sense of~\cite[Defn.\ 3.2.8]{EGstacktheoreticimages}) of the morphisms
$\cC^{\dd,\BT,a}\to\cR^{\dd,a}$ and
$\cC^{\tau,\BT,a}\to\cR^{\dd,a}$ respectively. We write $\ocZ$,
$\ocZ^\tau$ for $\cZ^1, \cZ^{\tau,1}$ respectively. 
\end{df}

The following theorem records some basic properties of these
scheme-theoretic images. (See Definition~\ref{defn: tres ram} for the
notion of a tr\`es ramifi\'ee representation.) 
\begin{thm} 
  \label{thm: existence of picture with descent data and its basic
    properties}\leavevmode
  \begin{enumerate}
  \item  For each $a\ge 1$, $\cZ^{\dd,a}$ is an algebraic
    stack of finite presentation over 
    $\cO/\varpi^a$, and is a closed substack
    of~$\cR^{\dd,a}$. In turn, each $\cZ^{\tau,a}$ is
    a closed substack of~$\cZ^{\dd,a}$, and thus in particular is an
    algebraic stack of finite presentation over $\cO/\varpi^a$;\ and $\cZ^{\dd,a}$ is the union of
    the $\cZ^{\tau,a}$. 
  \item The morphism $\cC^{\dd,\BT,a}\to\cR^{\dd,a}$ factors through a morphism
$\cC^{\dd,\BT,a}\to \cZ^{\dd,a}$
which is representable by algebraic spaces, scheme-theoretically dominant, and proper.
Similarly, the morphism
    $\cC^{\tau,\BT,a}\to\cR^{\dd,a}$ factors through a morphism
    $\cC^{\tau,\BT,a}\to\cZ^{\tau,a}$
which is representable by algebraic spaces, scheme-theoretically dominant, and proper.
  \item The $\Fpbar$-points of $\ocZ$ are naturally in bijection with the
    continuous representations $\rbar:G_K\to\GL_2(\Fpbar)$ which are not a
    twist of a tr\`es ramifi\'ee extension of the trivial character by the
    mod~$p$ cyclotomic character. Similarly, the $\Fpbar$-points of $\ocZ^\tau$ are naturally in bijection with the
    continuous representations $\rbar:G_K\to\GL_2(\Fpbar)$ which have
    a potentially Barsotti--Tate lift of type~$\tau$. 
    \end{enumerate}

\end{thm}

\begin{proof} 

Part~(1) follows easily from Theorem~\ref{thm: R^dd_h is an algebraic
  stack}. Indeed, by \cite[Prop.\ 3.2.31]{EGstacktheoreticimages} we may think of $\cZ^{\dd,a}$ as the scheme-theoretic image of the
proper morphism of algebraic stacks~$\cC^{\dd,\BT,a}\to\cR^{\dd,a}_1$ and
similarly for each $\cZ^{\tau,a}$. The
existence of the factorisations in~(2) is then formal.

By~\cite[Lem.\ 3.2.14]{EGstacktheoreticimages}, for each finite
extension $\F'/\F$, 
the $\F'$-points of $\ocZ$ (respectively $\ocZ^\tau$) correspond
to the \'etale $\varphi$-modules with descent data of the form~$\gM[1/u]$,
where~$\gM$ is a Breuil--Kisin module of rank $2$ with descent data and
$\F$-coefficients which satisfies the strong
determinant condition (respectively, which satisfies the strong determinant
condition and is of type~$\tau$). By Lemma~\ref{lem: O points of
  moduli stacks} and Corollary~\ref{cor: Kisin moduli consequences of local models}, these precisely
correspond to the Galois representations $\rbar:G_K\to\GL_2(\F)$ which
admit potentially Barsotti--Tate lifts of some tame type (respectively, of
type~$\tau$). The result  now follows from Lemma~\ref{lem: list of things
  we need to know about Serre weights}.
\end{proof}

The thickenings 
$\cC^{\dd, \BT, a} \hookrightarrow \cC^{\dd, \BT, a+1}$
and $\cR^{\dd, a} \hookrightarrow \cR^{\dd, a+1}$
induce closed immersions
$\cZ^{\dd, a} \hookrightarrow \cZ^{\dd, a+1}$.
Similarly, the thickenings $\cC^{\tau,\BT,a}\hookrightarrow \cC^{\tau,\BT, a+1}$
give rise to closed immersions $\cZ^{\tau,a} \hookrightarrow \cZ^{\tau,a+1}.$

\begin{lemma}
\label{lem:Z thickenings}
Fix $a \geq 1$.
Then the morphism
$\cZ^{\dd,a} \hookrightarrow \cZ^{\dd,a+1}$ is a thickening,
and for each tame type~$\tau,$
the morphism
$\cZ^{\tau,a} \hookrightarrow \cZ^{\tau,a+1}$ is a thickening.
\end{lemma}
\begin{proof}
In each case, the claim of the lemma follows from the following more general
statement:
if $$\xymatrix{ \cX \ar[r] \ar[d] & \cX' \ar[d] \\ \cY \ar[r] & \cY'}$$
is a diagram of morphisms of algebraic stacks in which the upper horizontal arrow
is a thickening, the lower horizontal arrow is a closed immersion,
and each of the vertical arrows is representable by algebraic spaces, quasi-compact,
and scheme-theoretically
dominant, then the lower horizontal arrow is also a thickening.

Since the property of being a thickening may be checked smooth locally, and since
scheme-theoretic dominance of quasi-compact morphisms is preserved 
by flat base-change, we may show this
after pulling the entire diagram back over a smooth surjective morphism
$V' \to \cY$  whose source is a scheme, and thus reduce to the case in which
the lower arrow is a morphism of schemes, and the upper arrow is a morphism
of algebraic spaces.  A surjecive \'etale morphism is also scheme-theoretically
dominant, and so pulling back the top arrow over a surjective \'etale morphism
$U' \to V'\times_{\cY'} \cX'$ whose source is a scheme, we finally reduce to
considering a diagram of morphisms of schemes
$$\xymatrix{U \ar[r]\ar[d] & U' \ar[d] \\ V \ar[r] & V'}$$
in which the top arrow is a thickening,
the vertical arrows are quasi-compact and scheme-theoretically dominant,
and the bottom arrow is a closed immersion.

Pulling back over an affine open subscheme of $V'$, and then pulling back the top arrow
over the disjoint union of the members of a finite affine open
cover of the preimage of this affine open in $U'$ (note that this preimage is quasi-compact),
we further reduce
to the case when all the schemes involved are affine.  That is, we have a diagram
of ring morphisms
$$\xymatrix{A' \ar[r] \ar[d] & A \ar[d] \\ B' \ar[r] & B}$$
in which the vertical arrows are injective, the horizontal arrows are surjective,
and the bottom arrow has nilpotent kernel.  One immediately verifies that
the top arrow has nilpotent kernel as well.
\end{proof}

We write $\cC^{\dd, \BT} := \varinjlim_a\cC^{\dd,\BT,a}$
and $\cZ^{\dd}:=\varinjlim_a\cZ^{\dd,a}$;\ 
we then
have evident morphisms of 
Ind-algebraic stacks
$$ {\cC^{\dd ,\BT}} \to {\cZ^{\dd}} \to
{\cR^{\dd}}$$
lying over $\Spf \cO$,
both representable by algebraic spaces, with the first being furthermore
proper and scheme-theoretically dominant in the sense of~\cite[Def.\ 6.13]{Emertonformalstacks}, 
and the second being a closed immersion.

Similarly, for each choice of tame type $\tau$,
we set
$\cC^{\tau,\BT}=\varinjlim_a\cC^{\tau,a}$, and $\cZ^{\tau}:=\varinjlim_a\cZ^{\tau,a}$.
We again have morphisms 
$$ {\cC^{\tau ,\BT}} \to {\cZ^{\tau}} \to
{\cR^{\dd}}$$
of Ind-algebraic stacks over $\Spf \cO$,
both being representable by algebraic spaces,
the first being proper and scheme-theoretically dominant, and the second
being a closed immersion. Note that by Corollary~\ref{cor: BT C is
  union of tau C}, $\cC^{\dd,\BT}$ is the disjoint union of the
~$\cC^{\tau,\BT}$, so it follows that~$\cZ^{\dd,\BT}$ is the union
(but \emph{not} the disjoint union) of
the~$\cZ^{\tau}$.

Proposition~\ref{prop:C tau BT is algebraic} shows that $\cC^{\tau,\BT} $ is a $\varpi$-adic formal algebraic stack
of finite presentation over $\Spf \cO$.
Each $\cZ^{\tau,a}$ is an algebraic stack of finite presentation
over $\Spec \cO/\varpi^a$ by Theorem~\ref{thm: existence of picture with descent data and its basic
    properties}. 
Analogous remarks apply in the case of $\cC^{\dd,\BT}$ and~$\cZ^{\dd}$.

\begin{prop} 
	\label{prop: Z is a p-adic formal stack}
	$\cZ^{\dd}$, and each $\cZ^{\tau}$, are $\varpi$-adic formal algebraic stacks,
	of finite presentation over $\Spf \cO$.
\end{prop}
\begin{proof}We give the argument for~$\cZ^\dd$, the argument
  for~$\cZ^{\tau}$ being identical. (Alternatively,
  this latter case follows from the former and the fact
  that the canonical morphism $\cZ^{\tau} \hookrightarrow \cZ^{\dd}$
  is a closed immersion.) 
	It follows from~\cite[Prop.\ A.10]{EGmoduli}
	that $\cZ^{\dd}$ is a formal algebraic stack,
	and by construction it is locally Ind-finite type over $\Spec \cO$.
	Since each $\cZ^{\dd,a}$ is 
        of finite presentation over $\cO/\varpi^a$,
         it follows from~\cite[Prop.\ A.21]{EGmoduli} that $\cZ^{\dd}$ is a $\varpi$-adic formal
  algebraic stack.

	The isomorphism $\cZ^{\dd} \iso \varinjlim_a \cZ^{\dd,a}$ 
	induces an isomorphism
	$$\cZ^{\dd}\times_{\cO} \cO/\varpi^b \iso
	\varinjlim_a \cZ^{\dd,a}\times_{\cO} \cO/\varpi^b,$$
	for any fixed $b \geq 1$.
	Since $\cZ^{\dd}$ is quasi-compact and quasi-separated,
	so is $\cZ^{\dd}\times_{\cO} \cO/\varpi^b,$ and thus this isomorphism
	factors through $\cZ^{\dd,a}\times_{\cO} \cO/\varpi^b$ for some $a$.
	Thus the directed system $\cZ^{\dd,a}\times_{\cO}
        \cO/\varpi^b$ in $a$ 
	eventually stabilises, and so we see that
	$$\cZ^{\dd}\times_{\cO} \cO/\varpi^b \iso \cZ^{\dd,a} \times_{\cO}
	\cO/\varpi^b$$
	for sufficiently large values of $a$.  Since $\cZ^{\dd,a}$
	is of finite presentation over $\cO/\varpi^a$,
	we find that $\cZ^{\dd}\times_{\cO} \cO/\varpi^b$
	is of finite presentation over $\cO/\varpi^b$.  Consequently,
	we conclude that $\cZ^{\dd}$ is of finite presentation
	over $\Spf \cO$, as claimed.
\end{proof}

\begin{remark}
The thickening $\cZ^{\dd,a} \hookrightarrow \cZ^{\dd}\times_{\cO}
\cO/\varpi^a$ need not be an isomorphism {\em a priori},
and we have no reason to expect that it is.  
Nevertheless, in \cite{cegsA} 
we prove that this thickening
is {\em generically} an isomorphism for {\em every} value of $a \geq 1$,
and we will furthermore show that each 
$(\cZ^{\dd,a})_{/\F}$ is generically reduced. 
The proof of this result 
depends on the detailed analysis of the irreducible
	components of the algebraic stacks $\cC^{\tau,a}$ 
	and $\cZ^{\tau,a}$ that we  make in \cite{cegsC}.
\end{remark}

We conclude this subsection
by establishing some basic lemmas about the reduced substacks underlying 
each of $\cC^{\tau, \BT}$ and $\cZ^{\tau}$.

\begin{lem}
    \label{lem: underlying reduced substack mod p}Let $\cX$ be an
    algebraic stack over $\cO/\varpi^a$, and let $\cX_{\red}$ be
    the underlying reduced substack of~$\cX$. Then $\cX_{\red}$ is a
    closed substack of $\cX_{/\F}:=\cX\times_{\cO/\varpi^a}\F$.
  \end{lem}
  \begin{proof}The structural morphism $\cX\to\Spec\cO/\varpi^a$
    induces a natural morphism $\cX_{\red}\to
    (\Spec\cO/\varpi^a)_{\red}=\Spec\F$, so the natural morphism
    $\cX_{\red}\to \cX$ factors through $\cX_{/\F}$.  Since the
    morphisms $\cX_{\red}\to\cX$ and $\cX_{/\F}\to\cX$ are both closed
    immersions, so is the morphism $\cX_{\red}\to\cX_{/\F}$.
      \end{proof}
\begin{lem}
  \label{lem: scheme-theoretic image commutes with passing to
    underlying reduced substacks}If $f:\cX\to\cY$ is a quasi-compact morphism of
  algebraic stacks, and $\cW$ is the scheme-theoretic image of $f$,
  then the scheme-theoretic image of the induced morphism of
  underlying reduced substacks $f_{\red}:\cX_{\red}\to\cY_{\red}$ is the
  underlying reduced substack $\cW_{\red}$. 
\end{lem}
\begin{proof} 
  Since the definitions of
  the scheme-theoretic image and of the underlying reduced substack
  are both smooth local (in the former case
  see~\cite[Rem. 3.1.5(3)]{EGstacktheoreticimages}, and in the latter
  case it follows immediately from the construction in~\cite[\href{http://stacks.math.columbia.edu/tag/0509}{Tag 0509}]{stacks-project}), we immediately reduce to the case of
  schemes, which follows from~\cite[\href{http://stacks.math.columbia.edu/tag/056B}{Tag 056B}]{stacks-project}. 
  \end{proof}

\begin{lem}
  \label{lem: C 1 and Z 1 are the underlying reduced substacks}
For each $a\ge 1$, $\cC^{\dd,\BT,1}$ is the underlying reduced
  substack of~$\cC^{\dd,\BT,a}$, and $\cZ^{\dd,1}$ is the
  underlying reduced substack of~$\cZ^{\dd,a}$;\ consequently,
 $\cC^{\dd,\BT,1}$ is the underlying reduced substack of $\cC^{\dd,\BT}$,
 and $\cZ^{\dd,1}$ is the underlying reduced substack of~$\cZ^{\dd}$.
 Similarly, for each tame type $\tau,$
  $\cC^{\tau,\BT,1}$ is the underlying reduced
  substack of each~$\cC^{\tau,\BT,a}$, and of~$\cC^{\tau,\BT}$,
  while $\cZ^{\tau,1}$ is the
  underlying reduced substack of each ~$\cZ^{\tau,a}$, and of~$\cZ^{\tau}$.
\end{lem}
\begin{proof} The statements for the $\varpi$-adic formal algebraic stacks
	follow directly from the corresponding statements for the various
	algebraic stacks modulo $\varpi^a$, and so we focus on proving 
	these latter statements,
	beginning with the case of~$\cC^{\tau,\BT,a}$. Note that
  $\cC^{\tau,\BT,1}=\cC^{\tau,\BT,a}\times_{\cO/\varpi^a}\F$ is
  reduced by Corollary~\ref{cor: Kisin moduli consequences of
    local models}, so
  $\cC^{\dd,\BT,1}=\cC^{\dd,\BT,a}\times_{\cO/\varpi^a}\F$ is
  also reduced by Corollary~\ref{cor:
  BT C is union of tau C}.  The claim follows for $\cC^{\dd,\BT,a}$
and $\cC^{\tau,\BT,a}$  from Lemma~\ref{lem:
    underlying reduced substack mod p}. 
  
The claims for $\cZ^{\tau,a}$  and $\cZ^{\dd,a}$ are then immediate from Lemma~\ref{lem: scheme-theoretic image commutes with passing to
    underlying reduced substacks}, applied to the morphisms
  $\cC^{\tau,\BT,a}\to\cZ^{\tau,a}$ and  $\cC^{\dd,\BT,a}\to\cZ^{\dd,a}$.
\end{proof}

\subsection{Versal rings and equidimensionality}\label{subsec:Galois
  deformation rings}
We now show that $\cC^{\dd,\BT}$ and $\cZ^{\dd,\BT}$ (and their
substacks $\cC^{\tau,\BT}$, $\cZ^{\tau}$) are equidimensional, and
compute their dimensions, by making use of their versal
rings.  In~\cite[\S 5]{EGstacktheoreticimages} these versal rings were
constructed in a more general setting in terms of liftings of \'etale
$\varphi$-modules;\ in our particular setting, we will find it
convenient to interpret them as Galois deformation rings.

We remark that this equidimensionality is not immediate from the local model results of Section~\ref{subsec: local models
  relation to explicit PR ones}, because we have not shown that the fibres of the local model map have constant dimensions;\ establishing this would be closely related to the arguments with versal rings that follow.

  Fix a finite type point $x:\Spec\F'\to\cZ^{\tau,a}$, where $\F'/\F$
is a finite extension;\ we also denote the induced finite type point of 
$\cR^{\dd,a}$ 
by~$x$.  Let $\rbar:G_K\to\GL_2(\F')$ be the
 Galois representation corresponding to~$x$ by Theorem~\ref{thm: existence of picture with descent data and its basic
    properties}~(3). Let~$E'$ be the compositum of~$E$
  and~$W(\F')[1/p]$, with ring of integers~$\cO_{E'}$ and residue
  field~$\F'$. 

  We let $R^\square_{\rbar}$ denote the universal framed deformation
  $\cO_{E'}$-algebra of $\rbar$, and we let $R^\square_{\rbar,0,\tau}$ be the reduced
  and $p$-torsion free quotient of $R^\square_{\rbar}$ whose
  $\Qpbar$-points correspond to the potentially Barsotti--Tate lifts
  of~$\rbar$ of type~$\tau$. In this section we will denote
  $R^\square_{\rbar,0,\tau}$ by the more suggestive name
  $R^{\tau,\BT}_{\rbar}$. We recall, for instance from
  \cite[Thm.~3.3.8]{BellovinGee}, that the ring $R^{\tau,\BT}_{\rbar}[1/p]$
  is regular.

As in Section~\ref{subsec: etale phi modules
  and Galois representations}, we write $R_{\rbar|_{G_{K_\infty}}}$
for the universal framed deformation
$\cO_{E'}$-algebra 
for~$\rbar|_{G_{K_\infty}}$. By Lemma~\ref{lem: Galois rep is a
  functor if A is actually finite local}, we have a natural morphism
\numequation
\label{eqn:versal morphism to R}
\Spf R_{\rbar|_{G_{K_\infty}}}\to\cR^{\dd}.
\end{equation}

\begin{lemma}
\label{lem:versal morphism to R}
The morphism~{\em \eqref{eqn:versal morphism to R}} 
is versal {\em (}at~$x${\em )}.
\end{lemma}
\begin{proof}
By definition, it suffices to show that if  
$\rho: G_{K_\infty} \to \GL_d(A)$ is a representation with $A$ a
finite Artinian $\cO_{E'}$-algebra, and if $\rho_B: G_{K_\infty} \to \GL_d(B)$
is a second representation, with $B$ a finite Artinian $\cO_{E'}$-algebra
admitting a surjection onto $A$, such that the base change
$\rho_A$ of $\rho_B$ to $A$ is isomorphic to $\rho$
(more concretely, so that there exists $M \in \GL_d(A)$ with
$\rho = M \rho_A M^{-1}$),
then we may find $\rho': G_{K_\infty} \to \GL_d(B)$ which lifts $\rho$,
and is isomorphic to $\rho_B$.  This is straightforward:\ the natural
morphism $\GL_d(B) \to \GL_d(A)$ is surjective, 
and so if $M'$ is any lift of $M$ to an element of $\GL_d(B)$,
then we may set $\rho' = M' \rho_B (M')^{-1}$.
\end{proof}

\begin{df}
\label{def:GL_2-hat}
For any pro-Artinian $\cO_{E'}$-algebra $R$ with residue field $\F'$ we let  $\widehat{\GL_2}_{/R}$ denote the completion
of~$(\GL_2)_{/R}$ along the closed
subgroup of its special fibre given by the centraliser of
$\rbar|_{G_{K_\infty}}$. 
\end{df}

\begin{remark}\label{rem:pullback-GL_2}
For $R$ as above we have $\widehat{\GL_2}_{/R} = \Spf R \times_{\cO_{E'}} \widehat{\GL_2}_{/\cO_{E'}}$.
Indeed, if $R \cong \varprojlim_i A_i,$ then
$\Spf R \times_{\cO_{E'}} \widehat{\GL_2}_{/\cO_{E'}} \cong \varinjlim_i
\Spec A_i \times_{\cO_{E'}} \widehat{\GL_2}_{/\cO_{E'}}$,
and 
$\Spec A_i \times_{\cO_{E'}} \widehat{\GL_2}_{/\cO_{E'}}$ agrees with the 
completion of $(\GL_2)_{/A_i}$ because $A_i$ is a finite $\cO_{E'}$-module.

It follows from this that $\widehat{\GL_2}_{/R}$ has nice base-change properties
more generally:
if $R \to S$ is a morphism of pro-Artinian $\cO_{E'}$-algebras each
with residue field~$\F'$, then there is an isomorphism $\widehat{\GL_2}_{/S}
\cong \Spf S \times_{\Spf R} \widehat{\GL_2}_{/R}.$
We apply this fact without further comment in various arguments below.
\end{remark}

There is a pair of morphisms
$\widehat{\GL_2}_{/R_{\rbar|_{G_{K_\infty}}}} \rightrightarrows
\Spf R_{\rbar|_{G_{K_\infty}}},$
the first being simply the projection to~$\Spf R_{\rbar|_{G_{K_\infty}}},$
and the second being given by ``change of framing''. 
Composing such changes of framing endows $\widehat{\GL_2}_{/R_{\rbar|_{G_{K_\infty}}}}$ with
the structure of a groupoid over~$\Spf R_{\rbar|_{G_{K_\infty}}}.$
Note that the two morphisms 
$$\widehat{\GL_2}_{/R_{\rbar|_{G_{K_\infty}}}} \rightrightarrows \Spf R_{\rbar|_{G_{K_\infty}}}
\buildrel \text{\eqref{eqn:versal morphism to R}} \over \longrightarrow
\cR^{\dd}$$
coincide, since changing the framing does not change the isomorphism
class (as a Galois representation) of a deformation of~$\rbar|_{G_{K_\infty}}$.
Thus there is an induced morphism of groupoids over~$\Spf R_{\rbar|_{G_{K_\infty}}}$
\numequation
\label{eqn:GL_2 hat equivariance}
\widehat{\GL_2}_{/R_{\rbar|_{G_{K_\infty}}}} \to
\Spf R_{\rbar|_{G_{K_\infty}}}\times_{\cR^{\dd}}\Spf R_{\rbar|_{G_{K_\infty}}}.
\end{equation}

\begin{lemma}
\label{lem:GL_2 hat equivariance}
The morphism~{\em \eqref{eqn:GL_2 hat equivariance}} is an isomorphism.
\end{lemma}
\begin{proof}
If $A$ is an Artinian $\cO_{E'}$-algebra, with residue field $\F'$,
then a pair of $A$-valued points of 
$\Spf R_{\rbar|_{G_{K_\infty}}}$ map to the same point of $\cR^{\dd}$
if and only if they give rise to isomorphic deformations of $\rhobar$,
once we forget the framings.  But this precisely means that the second of them is obtained
from the first by changing the framing via an $A$-valued point of $\widehat{\GL_2}$.
\end{proof}

It follows from Lemma~\ref{lem:versal morphism to R} that,
for each~$a\ge 1$, the quotient $R_{\rbar|_{G_{K_\infty}}}/\varpi^a$ is a
(non-Noetherian) versal ring for~$\cR^{\dd,a}$ at~$x$. By~\cite[Lem.\
3.2.16]{EGstacktheoreticimages}, for each~$a\ge 1$ a versal ring $R^{\tau,a}$ for
$\cZ^{\tau,a}$ at~$x$ is given by the
scheme-theoretic image of the
morphism
\numequation
\label{eqn:pull-back morphism mod a}
\cC^{\tau,\BT,a}\times_{\cR^{\dd,a}}\Spf R_{\rbar|_{G_{K_\infty}}}/\varpi^a \to\Spf
R_{\rbar|_{G_{K_\infty}}}/\varpi^a,
\end{equation}
in the sense that we now explain.

In general, the notion of scheme-theoretic image for morphisms 
of formal algebraic stacks can be problematic;\ at the very least it should
be handled with care.  But in this particular context,
a definition is given in~\cite[Def.\ 3.2.15]{EGstacktheoreticimages}:\ 
we write $R_{\rbar|_{G_{K_\infty}}}/\varpi^a$ as an inverse limit of Artinian local
rings $A$, form the corresponding scheme-theoretic images
of the induced morphisms 
$\cC^{\tau,\BT,a}\times_{\cR^{\dd,a}}\Spec A \to\Spec A,$
and then take the inductive limit of these scheme-theoretic images;
this is a formal scheme, which is in fact of the form $\Spf R^{\tau,a}$
for some quotient $R^{\tau,a}$ of $R_{\rbar_{| G_{K_{\infty}}}}/\varpi^a$
(where {\em quotient} should be understood in the sense of topological rings),
and is by definition the scheme-theoretic
image in question.

The closed
immersions $\cC^{\tau,\BT,a}\into \cC^{\tau,\BT,a+1}$ induce
corresponding closed immersions
$$
\cC^{\tau,\BT,a}\times_{\cR^{\dd,a}}\Spf R_{\rbar|_{G_{K_\infty}}}/\varpi^a \to
\cC^{\tau,\BT,a+1}\times_{\cR^{\dd,a+1}}\Spf R_{\rbar|_{G_{K_\infty}}}/\varpi^{a+1},$$
and hence closed immersions of scheme-theoretic images
$\Spf R^{\tau,a}\to \Spf R^{\tau,a+1}$, corresponding to surjections
$R^{\tau,a+1} \to R^{\tau,a}$.
(Here we are using the fact that
an projective limit of surjections of finite Artin
rings is surjective.)
Thus we may form the pro-Artinian ring 
$\varprojlim_a R^{\tau,a}.$
This projective limit is a quotient (again in the sense of topological rings)
of $R_{\rbar_{| G_{K_{\infty}}}}$, 
and the closed formal subscheme $\Spf (\varprojlim R^{\tau,a})$ 
of $\Spf R_{\rbar_{| G_{K_{\infty}}}}$ 
is the scheme-theoretic image (computed in the sense 
described above) of the projection
\numequation
\label{eqn:pull-back morphism}
{\cC^{\tau,\BT}}\times_{{\cR^{\dd}}}\Spf
R_{\rbar|_{G_{K_\infty}}}\to \Spf R_{\rbar|_{G_{K_\infty}}}
\end{equation}
(This is a formal consequence of the construction
of the $\Spf R^{\tau,a}$ as scheme-theoretic images,
since any discrete Artinian quotient of
$R_{\rbar |_{G_{K_{\infty}}}}$
is a discrete Artinian quotient of 
$R_{\rbar |_{G_{K_{\infty}}}}/\varpi^a$,
for some~$a~\geq~1$.)
It also follows formally (for example, by the
same argument as in the proof of~\cite[Lem.\
4.2.14]{EGstacktheoreticimages})  that $\varprojlim R^{\tau,a}$ is a versal
ring to ${\cZ^\tau}$ at~$x$.

\begin{remark}
The same method constructs a versal ring $R^{\dd,a}$ for $\cZ^{\dd,a}$ at $x$, and each $R^{\tau,a}$ is a quotient of $R^{\dd,a}$.
\end{remark}

Our next aim is to identify 
 $\varprojlim R^{\tau,a}$
with $R^{\tau,\BT}_{\rbar}$.  Before we do this, we have to establish some preliminary facts 
related to the various objects and morphisms we have just introduced.

\begin{lemma}  
\label{lem:R-tau-a properties}\leavevmode
\begin{enumerate}
\item
Each of the rings $R^{\tau,a}$ is a complete local Noetherian ring,
endowed with its $\mathfrak m$-adic topology,
and the same is true of the inverse limit $\varprojlim_a R^{\tau,a}$.
\item
For each $a \geq 1,$ the morphism $\Spf R^{\tau,a} \to \Spf R_{\rbar_{|G_{K_{\infty}}}}$
induces an isomorphism
$$
\cC^{\tau,\BT,a}\times_{\cR^{\dd,a}}\Spf R^{\tau,a} \iso
\cC^{\tau,\BT,a}\times_{\cR^{\dd,a}}\Spf R_{\rbar|_{G_{K_\infty}}}/\varpi^a.
$$
\item
For each $a \geq 1$,
the morphism $\Spf R^{\tau,a} \to \cR^{\dd,a}$ is effective, i.e.\
may be promoted {\em (}in a unique manner{\em )} to a morphism
$\Spec R^{\tau,a} \to \cR^{\dd,a}$,
and the induced morphism $$\cC^{\tau,a} \times_{\cR^{\dd,a}} \Spec R^{\tau,a} 
\to \Spec R^{\tau,a}$$ is proper and scheme-theoretically dominant.
\item
Each transition morphism $\Spec R^{\tau,a} \hookrightarrow \Spec R^{\tau,a+1}$
is a thickening.
\end{enumerate}
\end{lemma}
\begin{proof}
Recall that in Section~\ref{subsubsec: deformation rings and Kisin
  modules} we defined a Noetherian quotient $R_{\rbar|_{G_{K_\infty}}}^{\le 1}$
of $R_{\rbar|_{G_{K_\infty}}}$, which is naturally identified with the framed
deformation ring $R_{\rbar}^{[0,1]}$ by Proposition~\ref{prop:
  restriction gives an equivalence of height at most 1
  deformation}. It follows from~\cite[Lem.\
5.4.15]{EGstacktheoreticimages}  (via an argument almost identical to
the one in the proof of~\cite[Prop.\ 5.4.17]{EGstacktheoreticimages}) that the morphism
$\Spf \varprojlim R^{\tau,a}\into\Spf R_{\rbar|_{G_{K_\infty}}}$ factors through
$\Spf R_{\rbar|_{G_{K_\infty}}}^{\le 1}=\Spf R_{\rbar}^{[0,1]}$, and indeed
that $\varprojlim R^{\tau,a}$ is a quotient of $\Spf R_{\rbar}^{[0,1]};$
this proves~(1).

It follows by the very construction of the $R^{\tau,a}$ that the
morphism~\eqref{eqn:pull-back morphism mod a}
factors through the closed subscheme $\Spf R^{\tau,a}$
of $\Spf R_{\rbar_{|G_{K_{\infty}}}}/\varpi^a$.
The claim of~(2) is a formal consequence of this.

We have already observed that the morphism $\Spf R^{\tau,a} \to \cR^{\dd,a}$
factors through $\cZ^{\tau,a}$.  This latter stack is algebraic, and of finite
type over~$\cO/\varpi^a$.  
    It follows
from~\cite[\href{https://stacks.math.columbia.edu/tag/07X8}{Tag
  07X8}]{stacks-project}
that the morphism $\Spf R^{\tau,a} \to \cZ^{\tau,a}$ is effective.
Taking into account part~(1) of the present lemma,
we deduce from the theorem on formal functions that
the formal completion of the scheme-theoretic image of the projection
$$\cC^{\tau,a}\times_{\cR^{\dd,a}} \Spec R^{\tau,a} \to \Spec R^{\tau,a}$$
at the closed point of $\Spec R^{\tau,a}$ coincides
with the scheme-theoretic image of the morphism
$$\cC^{\tau,a}\times_{\cR^{\dd,a}} \Spf R^{\tau,a} \to \Spf R^{\tau,a}.$$
Taking into account~(2), we see that this latter scheme-theoretic image coincides with
$\Spf R^{\tau,a}$ itself.  This completes the proof of~(3).

The claim of~(4) follows from a consideration of the diagram
$$\xymatrix{\cC^{\tau,a} \times_{\cR^{\dd,a}} \Spec R^{\tau,a} \ar[r] \ar[d]
& \cC^{\tau,a+1} \times_{\cR^{\dd,a+1}} \Spec R^{\tau,a+1}  \ar[d] \\
\Spec R^{\tau, a} \ar[r] & \Spec R^{\tau,a+1} }$$
just as in the proof of Lemma~\ref{lem:Z thickenings}.
\end{proof}

\begin{lemma}
\label{lem:algebraization}\leavevmode
\begin{enumerate}
\item The projection
$\cC^{\tau,\BT} \times_{\cR^{\dd}} 
\Spf R_{\rbar_{|G_{K_{\infty}}}}
\to
\Spf R_{\rbar_{|G_{K_{\infty}}}}$
factors through a morphism
$\cC^{\tau,\BT} \times_{\cR^{\dd}} 
\Spf R_{\rbar_{|G_{K_{\infty}}}}
\to \Spf(\varprojlim R^{\tau,a})$,
which is scheme-theoretically dominant
in the sense that its scheme-theoretic image {\em (}computed in the manner
described above{\em )} is equal to its target.
\item There is a projective
morphism of schemes $X_{\rbar} \to \Spec (\varprojlim R^{\tau,a})$,
which is uniquely 
determined, up to unique isomorphism, by the requirement that its
$\mathfrak m$-adic completion {\em (}where $\mathfrak m$ denotes
the maximal ideal of $\varprojlim R^{\tau,a}${\em )} may be identified with
the morphism
$\cC^{\tau,\BT} \times_{\cR^{\dd}} 
\Spf R_{\rbar_{|G_{K_{\infty}}}}
\to \Spf(\varprojlim R^{\tau,a})$
of~(1).
\end{enumerate}
\end{lemma}
\begin{proof}
Part~(1) follows formally from the various constructions and definitions
of the objects involved (just like part~(2) of Lemma~\ref{lem:R-tau-a properties}).

We now consider the morphism
\[{\cC^{\tau,\BT}}\times_{{\cR^{\dd}}}\Spf R_{\rbar|_{G_{K_\infty}}}
\to \Spf (\varprojlim R^{\tau,a}).\]  Once we recall
that
$\varprojlim R^{\tau,a}$
is Noetherian, by Lemma~\ref{lem:R-tau-a properties}~(1),
it follows exactly as in the proof
of~\cite[Prop.\ 2.1.10]{kis04} (which treats the case that
$\tau$ is the trivial type), 
via an application of formal GAGA~\cite[Thm.\ 5.4.5]{MR217085},
that this morphism arises as the formal completion along the maximal ideal of
$\varprojlim R^{\tau,a}$
of a projective morphism
$X_{\rbar}\to \Spec (\varprojlim R^{\tau,a})$
(and $X_{\rbar}$ is unique up to unique isomorphism, by~\cite[Thm.\ 5.4.1]{MR217085}). 
\end{proof}

We next establish various properties of the scheme $X_{\rbar}$
constructed in the previous lemma.
To ease notation going forward, we write $\widehat{X}_{\rbar}$ to
denote the fibre product 
${\cC^{\tau,\BT}}\times_{{\cR^{\dd}}}\Spf R_{\rbar|_{G_{K_\infty}}}$
(which is reasonable, since this fibre product is isomorphic to the
formal completion of $X_{\rbar}$).

\begin{lemma}
\label{lem:X is normal and flat}
The scheme $X_{\rbar}$ is Noetherian, normal, and flat over $\cO_{E'}$.
\end{lemma}
\begin{proof}
Since $X_{\rbar}$ is projective over the Noetherian ring $\varprojlim R^{\tau,a}$,
it is Noetherian.   The other claimed properties of~$X_{\rbar}$
will be deduced from the corresponding properties of $\cC^{\tau,\BT}$ that
are proved in 
Corollary~\ref{cor: Kisin moduli consequences of local models}.

To this end, we first note that,
since the morphism $\cC^{\tau,\BT} \to \cR^{\dd}$ factors through~$\cZ^{\tau}$,
it follows (for example as in the proof of~\cite[Lem.\
3.2.16]{EGstacktheoreticimages}) that we have isomorphisms
\begin{multline*}
\cC^{\tau,\BT}\times_{\cZ^{\tau}}\Spf (\varprojlim R^{\tau,a})
\iso
\cC^{\tau,\BT}\times_{\cZ^{\tau}} \cZ^{\tau}\times_{\cR^{\dd}}
\Spf R_{\rbar|_{G_{K_\infty}}}
\\
\iso
\cC^{\tau,\BT}\times_{\cR^{\dd}} 
\Spf R_{\rbar|_{G_{K_\infty}}} =: \widehat{X}_{\rbar}.
\end{multline*}
In summary, we may identify
$\widehat{X}_{\rbar}$ with the fibre product
$\cC^{\tau,\BT}\times_{\cZ^{\tau}}\Spf (\varprojlim R^{\tau,a})$.

We now show that
$\widehat{X}_{\rbar}$ 
is analytically normal.   To see this,
let $\Spf B \to \widehat{X}_{\rbar}$
be a morphism whose source is a Noetherian affine formal algebraic
space, which is
representable by algebraic spaces and smooth. 
We must show
that the completion $\widehat{B}_{\mathfrak n}$ is normal,
for each maximal ideal $\mathfrak n$ of $B$.
In fact, it suffices to verify this for some collection of such $\Spf B$
which cover $\widehat{X}_{\rbar}$, and so without loss of generality
we may choose our $B$ as follows:
first, choose a collection of morphisms 
$\Spf A \to \cC^{\tau,\BT}$
whose sources are Noetherian affine formal algebraic spaces, and which are
representable by algebraic spaces and smooth, which, taken together, cover~$\cC^{\tau,\BT}$.
Next, for each such $A$, choose a collection of morphisms
$$\Spf B \to \Spf_A \times_{\cC^{\tau,\BT}} \widehat{X}_{\rbar}$$
whose sources are Noetherian affine formal algebraic spaces, and which are
representable by algebraic spaces and smooth, which, taken together, cover the
fibre product.  Altogether (considering all such $B$ associated to all such $A$),
the composite morphisms
$$\Spf B \to \Spf_A \times_{\cC^{\tau,\BT}} \widehat{X}_{\rbar} \to \widehat{X}_{\rbar}$$
are representable by algebraic spaces and smooth,
and cover~$\widehat{X}_{\rbar}$. 

Now, let $\mathfrak n$ be a maximal ideal in one of these rings $B$,
lying over a maximal ideal $\mathfrak m$ in the corresponding ring $A$.
The extension of residue fields $A/\mathfrak m \to B/\mathfrak n$ is
finite, and each of these fields is finite over $\F'$. Enlarging $\F'$ sufficiently,
we may assume that in fact each of these residue fields coincides with $\F'$.
(On the level of rings, this amounts to forming various tensor products of the
form $\text{--}\otimes_{W(\F')} W(\F'')$, which doesn't affect the question
of normality.)  
The morphism $\Spf B_{\mathfrak n} \to \Spf A_{\mathfrak n}$ is then
seen to be smooth
in the sense of~\cite[\href{https://stacks.math.columbia.edu/tag/06HG}{Tag
  06HG}]{stacks-project}, i.e., it satisfies the infinitesimal lifting property for 
finite Artinian $\cO'$-algebras with residue field $\F'$:\ this follows from
the identification of $\widehat{X}_{\rbar}$ above as a fibre product,
and the fact that $\Spf (\varprojlim R^{\tau,a}) \to \cZ^{\tau}$ is versal at
the closed point~$x$.
Thus $\Spf B_{\mathfrak n}$ is a formal power series ring
over $\Spf A_{\mathfrak m}$, 
by~\cite[\href{https://stacks.math.columbia.edu/tag/06HL}{Tag
  06HL}]{stacks-project}, 
and hence $\Spf B_{\mathfrak n}$ is indeed normal,
since $\Spf A_{\mathfrak m}$ is so,
by Corollary~\ref{cor: Kisin moduli consequences of local models}.
By Lemma~\ref{lem:GAGA facts} below,
this implies that the algebraization $X_{\rbar}$ of $\widehat{X}_{\rbar}$ is normal.

We next claim that the morphism
\numequation
\label{eqn:versal morphism to Z}
\Spf (\varprojlim R^{\tau,a}) \to \cZ^{\tau}
\end{equation}
is a flat morphism of formal algebraic stacks,
in the sense of \cite[Def.~8.35]{Emertonformalstacks}.
Given this, we find that the base-changed morphism
$\widehat{X}_{\rbar} \to \cC^{\tau,\BT}$
is also flat.
Since
Corollary~\ref{cor: Kisin moduli consequences of local models} shows that
$\cC^{\tau,\BT}$ is flat over $\cO_{E'},$ 
we conclude that the same is true of~$\widehat{X}_{\rbar}$.
Again, by Lemma~\ref{lem:GAGA facts}, this implies that the algebraization $X_{\rbar}$
is also flat over~$\cO_{E'}$.

It remains to show the claimed flatness.  To this end, we note first
that for each $a \geq 1$, the morphism
\numequation
\label{eqn:versal morphism to Z-a}
\Spf R^{\tau,a} \to \cZ^{\tau,a}
\end{equation}
is a versal morphism from a complete Noetherian local ring to an
algebraic stack which is locally of finite type over $\cO/\varpi^a$.
We already observed in the proof of Lemma~\ref{lem:R-tau-a properties}~(3) 
that~\eqref{eqn:versal morphism to Z-a} is effective,
i.e.\ can be promoted to a morphism
$\Spec R^{\tau,a} \to \cZ^{\tau,a}$.  It then follows 
from~\cite[\href{https://stacks.math.columbia.edu/tag/0DR2}{Tag
  0DR2}]{stacks-project} that this latter morphism is flat,
and thus that~\eqref{eqn:versal morphism to Z-a}
is flat in the sense of \cite[Def.~8.35]{Emertonformalstacks}.
It follows easily that the morphism~\eqref{eqn:versal morphism to Z}
is also flat:\ use the fact that a morphism
 of $\varpi$-adically complete local Noetherian $\cO$-algebras
 which becomes flat upon reduction modulo $\varpi^a$, for each $a \geq 1$,
 is itself flat, which follows from (for example) \cite[\href{https://stacks.math.columbia.edu/tag/0523}{Tag 0523}]{stacks-project}.
\end{proof}

The following lemma is standard, and is presumably well-known.  We sketch the proof,
since we don't know a reference.

\begin{lemma}
\label{lem:GAGA facts}
If $S$ is a complete Noetherian local $\cO$-algebra and $Y \to \Spec S$ is a proper
morphism of schemes, then $Y$ is flat over $\Spec \cO$
{\em (}resp.\ normal{\em )} if and only
$\widehat{Y}$ {\em (}the $\mathfrak m_S$-adic completion of~$Y${\em )}
is flat over $\Spf \cO$ {\em (}resp.\ is analytically normal{\em )}.
\end{lemma}
\begin{proof}
The properties of $Y$ that are in question can be tested by considering the 
various local rings $\cO_{Y,y}$, as $y$ runs over the points of $Y$;\ namely,
we have to consider whether or not these rings are flat over~$\cO$,
or normal.  Since
any point $y$ specializes to a closed point $y_0$ of $Y$, so that 
$\cO_{Y,y}$ is a localization of $\cO_{Y,y_0}$, and thus $\cO$-flat (resp.\
normal) if $\cO_{Y,y_0}$ is, it suffices to consider the rings $\cO_{Y,y_0}$
for closed points $y_0$ of~$Y$.  
Note also that since $Y$ is proper over $\Spec S$, any closed point of $Y$ lies over
the closed point of $\Spec S$.  

Now let $\Spec A$ be an affine neighbourhood of a closed point $y_0$ of $Y$;
let $\mathfrak m$ be the corresponding maximal ideal of~$A$.  As we noted,
$\mathfrak m$ lies over~$\mathfrak m_S$, and so gives rise to a maximal
ideal $\widehat{\mathfrak m} := \mathfrak m \widehat{A}$ of~$\widehat{A}$,
the $\mathfrak m_S$-adic completion of~$A$;
and any
maximal ideal of $\widehat{A}$ contains $\mathfrak m_S \widehat{A}$,
and so arises from a maximal ideal of~$A$ in this manner (since $A/\mathfrak m_S \iso
\widehat{A}/\mathfrak m_S$).  Write $\widehat{A}_{\mathfrak m}$ to denote
the $\mathfrak m$-adic completion of $A$ (which maps isomorphically to
the $\widehat{\mathfrak m}$-adic completion of $\widehat{A}$).  Then $\widehat{A}$
is faithfully flat over the localization $A_{\mathfrak m} = \cO_{Y,y_0}$,
and hence $A_{\mathfrak m}$ is flat over $\cO$ if and only if
$\widehat{A}_{\mathfrak m}$ is. Consequently we see that $Y$ is flat over 
$\cO$ if and only if, for each affine open subset $\Spec A$ of $Y$, 
the corresponding $\mathfrak m_S$-adic completion $\widehat{A}$ becomes flat over
$\cO$ after completing at each of its maximal ideals.  Another application
of faithful flatness of completions of Noetherian local rings shows
that this holds if and only if each such $\widehat{A}$ is flat over $\cO$ after localizing
at each of its maximal ideals, which holds if and only each
such $\widehat{A}$ is flat over~$\cO$.  This is precisely what it means for
$\widehat{Y}$ to be flat over~$\cO$.

The proof that analytic normality of $\widehat{Y}$ implies that $Y$ is normal
is similar.  Indeed, analytic normality by definition means that the
completion of $\widehat{A}$ at each of its maximal ideals is normal.  This completion
is faithfully flat over the localization of $\Spec A$ at its corresponding 
maximal ideal,
and so~\cite[\href{https://stacks.math.columbia.edu/tag/033G}{Tag 033G}]{stacks-project}
implies that this localization is also normal.  The discussion of the first paragraph
then implies that $Y$ is normal.  For the converse direction, we have to deduce
normality of the completions $\widehat{A}_{\mathfrak m}$ from the normality of 
the corresponding localizations~$A_{\mathfrak m}$.  This follows from that fact
that  $Y$
is an excellent scheme (being of finite type over the complete local ring~$S$),
so that each $A$ is an excellent
ring~\cite[\href{https://stacks.math.columbia.edu/tag/0C23}{Tag
    0C23}]{stacks-project}.
\end{proof}

\begin{prop}\label{prop: projective morphism as in Kisin}The projective morphism $X_{\rbar}\to \Spec R_{\rbar}^{[0,1]}$ factors
  through a projective and scheme-theoretically dominant morphism
  \numequation
  \label{eqn:projective morphism as in Kisin}
  X_{\rbar}\to \Spec R^{\tau,\BT}_{\rbar}
  \end{equation}
  which becomes an isomorphism after inverting~$\varpi$.  \end{prop}
\begin{proof}
We begin by showing the existence of~\eqref{eqn:projective morphism as in Kisin},
and that it induces a bijection on closed points after inverting~$\varpi$.
    Since $X_{\rbar}$ 
    is $\cO$-flat, by Lemma~\ref{lem:X is normal and flat},
it suffices to show that the induced morphism
    $$\Spec E \times_{\cO} X_{\rbar} \to \Spec R^{[0,1]}_{\rbar}[1/\varpi]$$
    factors through a morphism
  \numequation
  \label{eqn:projective morphism as in Kisin p inverted}
    \Spec E\times_{\cO} X_{\rbar} \to \Spec R^{\tau,\BT}_{\rbar}[1/\varpi],
\end{equation}
    which induces a bijection on closed points.

  This can be proved in exactly the same way as~\cite[Prop.\
2.4.8]{kis04}, 
    which treats the case that $\tau$ is
    trivial. 
    Indeed,
the computation of the $D_{\cris}$ of a Galois
    representation in the proof
    of~\cite[Prop.\ 2.4.8]{kis04} goes over essentially unchanged to
    the case of a Galois representation coming from $\cC^{\tau,\BT}$, and
    finite type points of $\Spec R_{\rbar}^{\tau,\BT}[1/\varpi]$ yield $p$-divisible
    groups and thus Breuil--Kisin modules exactly as in the proof
    of~\cite[Prop.\ 2.4.8]{kis04} (bearing in mind Lemma~\ref{lem: O
        points of moduli stacks} above). 
    The tame descent data comes along
    for the ride.

The morphism~\eqref{eqn:projective morphism as in Kisin p inverted}
is a projective morphism whose target is Jacobson, and which induces a bijection
on closed points.  It is thus proper and quasi-finite, and hence finite.
Its source is reduced (being even normal, by Lemma~\ref{lem:X is normal and flat}),
and its target is normal (as it is even regular, as we noted above).
A finite morphism whose source is reduced,
whose target is normal and Noetherian, and which
induces a bijection on finite type points, is indeed an isomorphism.
(The connected components of a normal scheme are integral,
and so base-changing over the connected components of the target,
we may assume that the target is integral.  The source is a union of finitely many 
irreducible components, each of which has closed image in the target.
Since the morphism is surjective on finite type  points, it is surjective,
and thus one of these closed images coincides with the target.  The injectivity
on finite type points then shows that the source is also irreducible,
and thus integral, as it is reduced.   
It follows from~\cite[\href{https://stacks.math.columbia.edu/tag/0AB1}{Tag 0AB1}]{stacks-project} that the morphism is an isomorphism.)
Thus~\eqref{eqn:projective morphism as in Kisin p inverted}
is an isomorphism.
Finally, since $R^{\tau,\BT}_{\rbar}$ is also flat over $\cO$ (by its definition),
this implies
that~\eqref{eqn:projective morphism as in Kisin}
is scheme-theoretically dominant.
  \end{proof}

\begin{cor}
  \label{cor: R tau BT is a versal ring to Z-hat}$\varprojlim
  R^{\tau,a}=R^{\tau,\BT}_{\rbar}$;\ thus $R^{\tau,\BT}_{\rbar}$ is a
  versal ring to ${\cZ^\tau}$ at~$x$.
\end{cor}
\begin{proof}
The theorem on formal functions shows that if we write the
scheme-theoretic image
of~\eqref{eqn:projective morphism as in Kisin}
in the form $\Spec B$, for some quotient $B$ of $R_{\rbar_{|G_{K_{\infty}}}}$,
then the scheme-theoretic image of the morphism~\eqref{eqn:pull-back morphism}
coincides with $\Spf B$.
The corollary then follows from Proposition~\ref{prop: projective morphism
as in Kisin}, which shows that~\eqref{eqn:projective morphism as in Kisin}
is scheme-theoretically dominant.
  \end{proof}

\begin{prop}
  \label{prop: dimensions of the Z stacks}
The algebraic stacks $\cZ^{\dd,a}$ and  $\cZ^{\tau,a}$ are equidimensional of
dimension~$[K:\Qp]$.
\end{prop}
\begin{proof}
Let $x$ be a finite type point of $\cZ^{\tau,a}$, defined over some finite extension
$\F'$ of $\F$, and corresponding to a Galois representation $\rbar$
with coefficients in~$\F'$.
By Corollary~\ref{cor: R tau BT is a versal ring to Z-hat}
the ring $R^{\tau,\BT}_{\rbar}$ coincides with the versal ring $\varprojlim_a R^{\tau, a}$
at $x$ of the $\varpi$-adic formal
algebraic stack~$\cZ^{\tau}$,
and so $\Spf R^{\tau,a} \iso \Spf R^{\tau,\BT}_{\rbar} \times_{\cZ^{\tau}} \cZ^{\tau,a}.$
Since $\cZ^{\tau}$ is a $\varpi$-adic formal algebraic stack,
the natural morphism $\cZ^{\tau,1} \to \cZ^{\tau}\times_{\Spf \cO}\F$
is a thickening, and thus the same is true of the morphism
$\Spf R^{\tau,1} \to \Spf R^{\tau,\BT}_{\rbar}/\varpi$
obtained by pulling the former morphism back over~$\Spf R^{\tau,\BT}_{\rbar}/\varpi$.

  Since $R^{\tau,\BT}_{\rbar}$ is
  flat over $\cO_{E'}$ and 
  equidimensional of dimension $5+[K:\Qp]$, it follows that
  $R^{\tau,1}$ 
  is equidimensional of dimension~$4~+~[K:\Qp]$.
  The same is then true of each~$R^{\tau,a}$, since these are thickenings
  of~$R^{\tau,1}$, by Lemma~\ref{lem:R-tau-a properties}~(4).

We have a versal morphism $\Spf
  R^{\tau,a}\to\cZ^{\tau,a}$ at the finite type point~$x$ of $\cZ^{\tau,a}$. It follows from Lemma~\ref{lem:GL_2 hat
    equivariance} that 
\[\widehat{\GL_2}_{/\Spf R^{\tau,a}} \iso
\Spf R^{\tau,a}\times_{\cZ^{\tau,a}}\Spf R^{\tau,a}.\]

To find the dimension of $\cZ^{\tau,a}$ it
suffices to compute its dimension at finite type points (\emph{cf}.\ \cite[\href{https://stacks.math.columbia.edu/tag/0DRX}{Tag
  0DRX}]{stacks-project}, recalling the definition of the dimension of
an algebraic stack,
\cite[\href{https://stacks.math.columbia.edu/tag/0AFP}{Tag
  0AFP}]{stacks-project}).  It follows from~\cite[Lem.\ 2.40]{2017arXiv170407654E}
applied to the presentation
  $[\Spf R^{\tau,a} / \widehat{\GL_2}_{/\Spf R^{\tau,a}}]$ 
  of $\widehat{\cZ}^{\tau,a}_{x}$, together with Remark~\ref{rem:pullback-GL_2},
that $\cZ^{\tau,a}$ is
equidimensional of dimension $[K:\Qp]$. Since $\cZ^{\dd,a}$ is the
 union of the $\cZ^{\tau,a}$ by Theorem~\ref{thm: existence of picture with descent data and its basic
    properties},  $\cZ^{\dd,a}$ is also
equidimensional of dimension $[K:\Qp]$
by~\cite[\href{https://stacks.math.columbia.edu/tag/0DRZ}{Tag
  0DRZ}]{stacks-project}. 
  \end{proof}

\begin{prop}
  \label{prop: C tau is equidimensional of the expected dimension}The
  algebraic stacks $\cC^{\tau,\BT,a}$ are equidimensional of
  dimension $[K:\Qp]$.
\end{prop}
\begin{proof}
Let $x'$ be a finite type point of  $\cC^{\tau,\BT,a}$, defined
over some finite extension $\F'$ of $\F$,
lying over the
finite type point~$x$ of~$\cZ^{\tau,a}$.  Let $\rbar$ be the Galois
representation with coefficients in $\F'$ corresponding to~$x$,
and recall that $X_{\rbar}$ denotes a projective $\Spec R_{\rbar}^{\tau,\BT}$-scheme
whose pull-back $\widehat{X}_{\rbar}$
over $\Spf R_{\rbar_{|G_{K_{\infty}}}}$ is isomorphic to
${\cC^{\tau,\BT}}\times_{{\cR^{\dd}}}\Spf R_{\rbar|_{G_{K_\infty}}}.$
The point $x'$ gives rise to a closed point $\tx$ of $X_{\rbar}$ (of which $x'$
is the image under the morphism $X_{\rbar} \to \cC^{\tau,\BT}$).
Let
$\widehat{\cO}_{X_{\rbar},\tx}$ denote the complete local ring to
$X_{\rbar}$ at the point~$\tx$;\ then the natural morphism $\Spf
\widehat{\cO}_{X_{\rbar},\tx}\to {\cC^{\tau,\BT}}$
is versal at~$\tx$,
so that $\widehat{\cO}_{X_{\rbar},\tx}/\varpi^a$ is a versal ring for the point
$x'$ of $\cC^{\tau,\BT,a}$.

The isomorphism~\eqref{eqn:GL_2 hat equivariance} induces
(after pulling back over $\cC^{\tau,\BT}$)
an isomorphism 
$$\widehat{\GL_2}_{/\widehat{X}_{\rbar}} \iso
\widehat{X}_{\rbar}\times_{\cC^{\tau,\BT}} \widehat{X}_{\rbar},$$
and thence an isomorphism
$$\widehat{\GL_2}_{/\widehat{\cO}_{X_{\rbar},\tx}} \iso
\widehat{\cO}_{X_{\rbar},\tx}
\times_{\cC^{\tau,\BT}} 
\widehat{\cO}_{X_{\rbar},\tx}.$$

Since $R^{\tau,\BT}$ is  equidimensional of dimension $5+[K:\Qp]$,
it follows from Proposition~\ref{prop: projective morphism as in Kisin} that
$X_{\rbar}$ is equidimensional of dimension $5+[K:\Qp]$,
and thus (taking into account the flatness statement of
Lemma~\ref{lem:X is normal and flat}) that
$\widehat{\cO}_{X_{\rbar},\tx}/\varpi^a$ is equidimensional of
dimension $4+[K:\Qp]$.
As in the proof of 
  Proposition~\ref{prop: dimensions of the Z stacks},
an application of~\cite[Lem.\ 2.40]{2017arXiv170407654E}
shows that $\dim_{x'}{\cC^{\tau,\BT,a}}$ is equal to $[K:\Qp]$. 
Since~$x'$ was an arbitrary finite type point, the result follows.
\end{proof}

\renewcommand{\theequation}{\Alph{section}.\arabic{subsection}} 
\appendix

\section{Tr\`es ramifi\'ee representations and Serre weights}\label{sec: appendix on tame
  types}

To complete our characterization in Theorem~\ref{thm: existence of picture with descent data and its basic
    properties}(3) of the $\F$-points of the
stack $\cZ^{\dd,1}$, in this appendix we give a global argument that requires a brief  detour into the theory of Serre weights.

By definition, a \emph{Serre weight} is an irreducible
$\F$-representation of $\GL_2(k)$, where we recall that $k$ denotes the residue field of $K$. Concretely, such a
representation is of the form
\[\sigmabar_{\vec{t},\vec{s}}:=\otimes^{f-1}_{j=0}
(\det{\!}^{t_j}\Sym^{s_j}k^2) \otimes_{k,\sigma_{j}} \F,\]
where $0\le s_j,t_j\le p-1$ and not all $t_j$ are equal to
$p-1$. We say that a Serre weight is \emph{Steinberg} if $s_j=p-1$ for all $j$,
and \emph{non-Steinberg} otherwise.

We now recall the definition of the set of Serre weights~$W(\rbar)$
associated to a representation $\rbar:G_K\to\GL_2(\Fpbar)$. 

\begin{adefn}\label{def:de rham of type}
 We say that a crystalline representation $r : G_K \to
  \GL_2(\Qpbar)$ has \emph{type $\sigmabar_{\vec{t},\vec{s}}$}
  provided that for each embedding $\sigma_j : k \into \F$ there is an
  embedding $\widetilde{\sigma}_j : K \into \Qpbar$ lifting $\sigma_j$
  such that the $\widetilde{\sigma}_j$-labeled Hodge--Tate weights
  of~$r$ are $\{-s_j-t_j,1-t_j\}$, and the remaining $(e-1)f$ pairs of  Hodge--Tate weights
  of~$r$ are all $\{0,1\}$. \emph{(}In particular the representations of
   type $\sigmabar_{\vec{0},\vec{0}}$ \emph{(}the trivial weight\emph{)} are the same as those of Hodge type $0$.\emph{)}
\end{adefn}

\begin{adefn}\label{def:serre weights}
Given a representation $\rbar:G_K\to\GL_2(\Fpbar)$ we define $W(\rbar)$
to be the set of Serre weights $\sigmabar$ such that $\rbar$ has a
crystalline lift of type $\sigmabar$.
\end{adefn}

It follows easily from the formula
$\sigmabar_{\vec{t},\vec{s}}^\vee=\sigmabar_{-\vec{s}-\vec{t},\vec{s}}$ 
that $\sigmabar \in W(\rbar)$ if and only if $\sigmabar^\vee$ is in
the set of Serre weights associated to $\rbar^\vee$  
in~\cite[Defn.\ 4.1.3]{gls13}.  

Let $\sigma(\tau)$
    denote $\GL_2(\cO_K)$-representation associated to $\tau$ by
    Henniart's inertial local Langlands correspondence
    \cite{breuil-mezard}, and let $\sigmabar(\tau)$ denote the special fibre
    of any lattice in $\sigma(\tau)$, well-defined up to semi-simplification.

    \begin{remark}\label{rem:normalisations}
      There are several definitions of the set $W(\rbar)$ in the literature,
which by the papers~\cite{blggu2,geekisin,gls13} are known to be
equivalent (up to normalisation). However the normalisations of
Hodge--Tate weights and of inertial local Langlands used in
\cite{geekisin,gls13} are not all the same.  Our conventions for Hodge--Tate weights and
inertial types agree with those of~\cite{geekisin}, but our
representation~$\sigma(\tau)$ is the
representation~$\sigma(\tau^\vee)$ of~\cite{geekisin}
(where~$\tau^\vee=\eta^{-1}\oplus(\eta')^{-1}$);\ to see this, note the
dual in the definition of~$\sigma(\tau)$ in~\cite[Thm.\
2.1.3]{geekisin}.

On the other hand, our set of weights~$W(\rbar)$ is the
set of duals of the weights~$W(\rbar)$ considered
in~\cite{geekisin}. In turn, the paper~\cite{gls13} has the opposite
convention for the signs of Hodge--Tate weights to our convention (and
to the convention of~\cite{geekisin}), so we find that our set of
weights~$W(\rbar)$ is the set of duals of the weights~$W(\rbar^\vee)$
considered in~\cite{gls13}. 
    \end{remark}

    We have the following standard definition.
    \begin{adefn}\label{defn: tres ram}
      We say that~$\rbar$ is~\emph{tr\`es ramifi\'ee} if it is a twist
      of an extension of the trivial character by the mod~$p$
      cyclotomic character, and if furthermore the splitting field of
      its projective image is \emph{not} of the form
      $K(\alpha_1^{1/p},\dots,\alpha_s^{1/p})$ for some
      $\alpha_1,\dots,\alpha_s\in\cO_K^\times$. 
    \end{adefn}

    The following characterisation of tr\`es ramifi\'ee
    representations presumably admits a purely local proof, but we
    find it convenient to make a global argument.
\begin{alemma}\label{lem: list of things we need to know about Serre weights}\leavevmode
  \begin{enumerate}
  \item If ~$\tau$ is a tame type, then $\rbar$ has a potentially
    Barsotti--Tate lift of type~$\tau$ if and only
    if $W(\rbar)\cap\JH(\sigmabar(\tau))\ne 0$.

  \item The following conditions are equivalent:
    \begin{enumerate}
    \item $\rbar$ admits a potentially Barsotti--Tate lift of some tame type.
    \item $W(\rbar)$ contains a non-Steinberg Serre weight.
    \item $\rbar$ is not tr\`es    ramifi\'ee.
    \end{enumerate}
  \end{enumerate}
\end{alemma}
\begin{proof}
(1)   By the main result of~\cite{gls13}, and bearing in mind the
differences between our conventions and those of~\cite{geekisin} as
recalled in Remark~\ref{rem:normalisations}, we have
$\sigmabar\in W(\rbar)$ if and only if $\sigmabar^\vee\in \WBT(\rbar)$,
where $\WBT(\rbar)$ is the set of weights defined
  in~\cite[\S3]{geekisin}. 
  By~\cite[Cor.\ 3.5.6]{geekisin} (bearing in mind once again the
  differences between our conventions and those of~\cite{geekisin}),
  it follows that we have $W(\rbar)\cap\JH(\sigmabar(\tau))\ne 0$ if
  and only if $e(R^{\square,0,\tau}/\pi)\ne 0$ in the notation of
  \emph{loc.\ cit.},
  and 
 by definition $\rbar$ has a potentially Barsotti--Tate lift of
  type~$\tau$ if and only if $R^{\square,0,\tau}\ne 0$. It follows
  from~\cite[Prop.\ 4.1.2]{emertongeerefinedBM} 
that
  $R^{\square,0,\tau}\ne 0$ if and only if
  $e(R^{\square,0,\tau}/\pi)\ne 0$, as required.

  (2) By part~(1), condition~(a) is equivalent to $W(\rbar)$
  containing a Serre weight occurring as a Jordan--H\"older factor of
  ~$\sigmabar(\tau)$ for some tame type~$\tau$. It is easy to see (from explicit tables of the Serre weights occurring as Jordan--H\"older factors of tame types, as computed e.g.\ in \cite{MR2392355})
    that the Serre weights
  occurring as Jordan--H\"older factors of the~$\sigmabar(\tau)$ are precisely the non-Steinberg Serre weights,
  so~(a) and~(b) are equivalent.

  Suppose that~(a) holds;\ then~$\rbar$ becomes finite flat over a tame
  extension. However the restriction to a tame extension of a tr\`es
  ramifi\'ee representation is still tr\`es
  ramifi\'ee, and therefore not finite flat, so~(c) also
  holds. Conversely, suppose for the sake of contradiction that~(c)
  holds, but that~(b) does not hold, i.e.\ that $W(\rbar)$ consists of
  a single Steinberg weight.

Twisting, we
can without loss of generality assume that~$W(\rbar)=\{\sigmabar_{\vec{0},\vec{p-1}}\}$. By~\cite[Cor.\ A.5]{geekisin} we
can globalise~$\rbar$, and then the hypothesis that~$W(\rbar)$
contains $\sigmabar_{\vec{0},\vec{p-1}}$ implies that it has a semistable lift of
Hodge type~$0$. If this lift were in fact crystalline, then~$W(\rbar)$
would also contain the weight $\sigmabar_{\vec{0},\vec{0}}$ by~(1). So
this lift is not crystalline, and in particular the monodromy
operator~$N$ on the corresponding weakly admissible module is
nonzero. But then $\ker(N)$ is a free filtered submodule of rank $1$,
and since the lift has Hodge type~$0$, $\ker(N)$ is in fact a weakly
admissible submodule. It follows that the lift is an unramified twist of an
extension of $\varepsilon^{-1}$ by the  trivial
character, so that~$\rbar$ is an unramified twist of an extension of
$\varepsilonbar^{-1}$ by the trivial character. But we are assuming
that~(c) holds, so~$\rbar$ is finite flat, so that by~(1), $W(\rbar)$
contains the weight $\sigmabar_{\vec{0},\vec{0}}$, a contradiction.
\end{proof}

\bibliographystyle{amsalpha-custom}
\bibliography{dieudonnelattices}
\end{document}